# Group cohomology and the singularities of the Selberg zeta function associated to a Kleinian group

By Ulrich Bunke and Martin Olbrich

## Contents

## 1. Introduction

Let $G := \mathrm{SO}(1, n)_0$ denote the group of orientation-preserving isometries of the $n$-dimensional hyperbolic space $X := H^n$ equipped with the Riemannian metric of constant sectional curvature $-1$. We consider a convex cocompact, torsion-free discrete subgroup $\Gamma \subset G$. The quotient $Y := \Gamma \backslash X$ is a complete hyperbolic manifold, and we assume that $\mathrm{vol}(Y) = \infty$.

The Selberg zeta function $Z_S(s)$, $s \in \mathbf{C}$, associated to this geometric situation encodes the length spectrum of closed geodesics of $Y$ together with the eigenvalues of their Poincaré maps. Note that $Z_S(s)$ is given by an Euler product on some half-plane $\mathrm{Re}(s) > c$, and it has a meromorphic continuation to the complex plane.

The goal of the present paper is a description of the singularities of the Selberg zeta function in terms of the group cohomology of $\Gamma$ with coefficients in certain infinite dimensional representations. Such a relation was conjectured by Patterson [37].

### 1.1. *The Selberg zeta function.*

In order to fix our conventions we define $Z_S(s)$ in terms of group theory. Let $\mathbf{g} = \mathbf{k} \oplus \mathbf{p}$ be a Cartan decomposition of the Lie algebra of $G$, where $\mathbf{k}$ is the Lie algebra of a maximal compact subgroup $K \subset G$, $K \cong \mathrm{SO}(n)$. We fix a one-dimensional subspace $\mathbf{a} \subset \mathbf{p}$ and let $M \subset K$, $M \cong \mathrm{SO}(n-1)$, denote the centralizer of $\mathbf{a}$. The Riemannian metric of $X$ induces a metric on $\mathbf{a}$. We fix an isometry $\mathbf{a} \cong \mathbf{R}$. Let $\mathbf{a}_+$ denote the half-space corresponding to the ray $[0, \infty)$ and set $A := \exp(\mathbf{a})$, $A_+ := \exp(\mathbf{a}_+)$. By $\mathbf{n} \subset \mathbf{g}$ we denote the positive root space of $\mathbf{a}$ in $\mathbf{g}$. For $H \in \mathbf{a}$ we set $\rho(H) := \frac{1}{2}\mathrm{tr}\,\mathrm{ad}(H)_{|\mathbf{n}}$. The isometry $\mathbf{a}^* \cong \mathbf{R}$ identifies $\rho$ with $\frac{n-1}{2}$.

By the Cartan decomposition $G = KA_+K$, any element $g \in G$ can be written as $g = ha_g k$, $h, k \in K$, where $a_g \in A_+$ is uniquely determined. We have $a_g = \mathrm{e}^{\mathrm{dist}(\mathcal{O}, g\mathcal{O})}$, where $\mathcal{O} = K \in X$ is the origin of $X = G/K$, and dist denotes the hyperbolic distance. The basic quantity associated with $\Gamma$ is its exponent $\delta_\Gamma$ which measures the growth of $\Gamma$ at infinity.

*Definition* 1.1. The exponent $\delta_\Gamma \in \mathbf{a}^* \cong \mathbf{R}$ of $\Gamma$ is defined to be the smallest number such that the Poincaré series

$$(1) \qquad\qquad \sum_{g \in \Gamma} a_g^{-(s+\rho)}$$

converges for all $s > \delta_\Gamma$.

Although we do not use it in the present paper note that $\delta_\Gamma + \frac{n-1}{2}$ is equal to the Hausdorff dimension of the limit set of $\Gamma$ (see [35], [47]).

Any element $g \in \Gamma$ is conjugated in $G$ to an element of the form $m(g)a(g) \in MA_+$, where $a(g)$ is unique. By $l_g := \log(a(g))$ we denote the length of the closed geodesic of $Y$ corresponding to the conjugacy class of $g \in \Gamma$.

Let $C\Gamma$ denote the set of conjugacy classes $[g] \neq 1$ of $\Gamma$. If $[g] \in C\Gamma$, then $n_\Gamma(g) \in \mathbf{N}$ is the multiplicity of $[g]$, i.e., the largest number $k \in \mathbf{N}$ such that $[g] = [h^k]$ for some $[h] \in C\Gamma$. If $A : V \to V$ is a linear homomorphism of a complex vector space $V$, then by $S^k A : S^k V \to S^k V$ we denote its $k^{\text{th}}$ symmetric power.

*Definition* 1.2. The Selberg zeta function $Z_S(s)$, $s \in \mathbf{C}$, $\text{Re}(s) > \delta_\Gamma$, associated to $\Gamma$ is defined by the infinite product

$$(2) \qquad Z_S(s) := \prod_{[g] \in C\Gamma, n_\Gamma(g)=1} \prod_{k=0}^\infty \det\left(1 - \mathrm{e}^{-(s+\rho)l_g} S^k(\mathrm{Ad}(m(g)a(g))_{|\mathbf{n}}^{-1})\right) \ .$$

*Remark.* We defined the Selberg zeta function such that its critical line is $\{\text{Re}(s) = 0\}$. In the literature the convention is often such that the critical line is at $\rho = \frac{n-1}{2}$. The same applies to our definition of $\delta_\Gamma$ which differs by $\rho$ from the usual convention.

In [38] (see also [42]) it was shown that the infinite product converges for $\text{Re}(s) > \delta_\Gamma$, and that the Selberg zeta function has a meromorphic continuation to all of $\mathbf{C}$. In the special case of surfaces this was also proved in [17]. Partial results concerning the logarithmic derivative of the Selberg zeta function have been obtained in [34] for $\delta_\Gamma < 0$ and in [41] in the general case.

## 1.2. *Singularities and spectrum.*

The Selberg zeta function is a meromorphic function defined in terms of a classical Hamiltonian system, namely the geodesic flow on the unit sphere bundle $SY$ of $Y$. Philosophically, the singularities of the Selberg zeta function should be considered as quantum numbers of an associated quantum mechanical system. One way to quantize the geodesic flow is to take as the Hamiltonian the Laplace-Beltrami operator $\Delta_Y$ acting on functions on $Y$.

In order to explain this philosophy let $Y = \Gamma\backslash G/K$ for a moment be a compact locally symmetric space of rank one. Then the sphere bundle of $Y$ can be written as $SY = \Gamma\backslash G/M$. If $(\sigma, V_\sigma)$ is a finite-dimensional unitary representation of $M$, then we consider the bundle $V(\sigma) := \Gamma\backslash G \times_M V_\sigma$ over $SY$. The geodesic flow admits a lift to $V(\sigma)$ and gives rise to a more general Selberg zeta function $Z_S(s, \sigma)$ which also encodes the holonomy in $V(\sigma)$ of the flow along the closed geodesics. The Selberg zeta function $Z_S(s)$ defined above corresponds to the trivial representation of $M$. It was shown in [13] that $Z_S(s, \sigma)$ admits a meromorphic continuation to all of $\mathbf{C}$. In this generality a description of the singularities of $Z_S(s, \sigma)$ was first obtained by [21] (see also [48] for a closely related Selberg zeta function).

*Remark.* The case of a Riemann surface is classical. The Selberg zeta functions for general rank-one symmetric spaces and trivial $\sigma$ have been discussed in [14]. For a detailed account of the literature see [5].

The description of the singularities of $Z_S(s, \sigma)$ given in [21] corresponds to a different method of quantization of the geodesic flow (see subsection 1.3). The spectral description of the singularities of $Z_S(s, \sigma)$ uses differential operators acting on sections of bundles on $Y$.

One distinguishes between two types of singularities, so-called topological and spectral singularities. If $\sigma$ is trivial then the spectral singularities are connected with the eigenvalues of the Laplace-Beltrami operator $\Delta_Y$. The topological singularities depend on the spectrum of the Laplace-Beltrami operator on the compact dual symmetric space to $X$. For general $\sigma$ we found the corresponding quantum mechanical system in [5]; i.e., we determined the locally homogeneous vector bundle on $Y$ together with a corresponding locally invariant differential operator whose eigenvalues are responsible for the spectral singularities of $Z_S(s, \sigma)$. The analytic continuation of this operator to the compact dual symmetric space gives rise to the topological singularities.

We now return to our present case that $\Gamma \subset \mathrm{SO}(1, n)_0$ is convex cocompact, $Y$ is a noncompact hyperbolic manifold, and $\sigma$ is the trivial representation of $M$. It was shown in [28] that the spectrum of $\Delta_Y$ consists of finitely many isolated eigenvalues in the interval $[(\frac{n-1}{2})^2 - |\delta_\Gamma|\delta_\Gamma, (\frac{n-1}{2})^2)$. Moreover, in [29] the same authors show that the remaining spectrum of $\Delta_Y$ is the absolute continuous spectrum of infinite multiplicity in the interval $[(\frac{n-1}{2})^2, \infty)$. It turns out that the eigenvalues of $\Delta_Y$ are responsible for singularities of $Z_S(s)$ as in the cocompact case. This stage of understanding is not satisfactory. On the one hand $Z_S(s)$ may have more singularities. On the other hand the continuous spectrum was neglected.

A finer investigation of the continuous spectrum can be based on study of the resolvent kernel, i.e. the distributional kernel of the inverse $(\Delta_Y - (\frac{n-1}{2})^2 + \lambda^2)^{-1}$. It is initially defined for $\mathrm{Re}(\lambda) \gg 0$. A continuation of this kernel up to the imaginary axis implies absolute continuity of the essential spectrum by the limiting absorption principle (see [39] and for surfaces also [12], [10]). But this kernel behaves much better. It was shown in [32] that it has a meromorphic continuation to the whole complex plane (for surfaces see also [9], [1]). The poles of this continuation with positive real part correspond to the eigenvalues of $\Delta_Y$. The poles with nonpositive real part are called resonances.

Let us consider the resonances as sorts of eigenvalues associated to the continuous spectrum. Then they lead to singularities of $Z_S(s)$ in the same way as true eigenvalues. In detail, the spectral description of the singularities of $Z_S(s)$ was worked out in [38] for even dimensions $n$.

### 1.3. *Singularities and group cohomology.*

An important feature of the spectral description of the singularities of the Selberg zeta function is the distinction between spectral and topological singularities. By now there are two approaches to describe *all* singularities of the Selberg zeta function in a *uniform way*, avoiding a separation of topological and spectral singularities.

We will explain these approaches again for $Z_S(s, \sigma)$ in the case where $Y$ is a compact locally symmetric space of rank one.

The *first approach* was worked out in [21] (ideas can be traced back to [16], [36]) and corresponds to a quantization which is different from the one considered for the spectral description. Here one considers the cohomology of the $\sigma$-twisted tangential de Rham complex (called tangential cohomology), i.e. the restriction of the de Rham complex of $SY$ to the stable foliation twisted with $V(\sigma)$. This complex is equivariant with respect to the flow. The order of the singularity of $Z_S(s, \sigma)$ at $s = \lambda$ is related to the Euler characteristic of the $\lambda + \rho$-eigenspace of the flow generator on the tangential cohomology.

The tangential cohomology comes with a natural topology, and it is still an open problem to show that this topology is Hausdorff. Therefore, in [21] the result is phrased in terms of representation theory, in particular in terms of Lie algebra cohomology of $\mathbf{n}$ with coefficients in the Harish-Chandra modules of the unitary representations occurring in the decomposition of the right regular representation of $G$ on $L^2(\Gamma \backslash G)$.

The *second approach* was proposed by S. Patterson [37]. The parameters $\sigma$ and $\lambda \in \mathbf{C}$ fix a principal series representation $(\pi^{\sigma,\lambda}, H^{\sigma,\lambda})$ of $G$. The space $H^{\sigma,\lambda}$ can be realized as the space of sections of a homogeneous vector bundle $V(\sigma_\lambda) := G \times_P V_{\sigma_\lambda}$, where $P := MAN$, $N := \exp(\mathbf{n})$, and $\sigma_\lambda$ is the representation of $P$ on $V_\sigma$ given by $P = MAN \ni man \mapsto a^{\rho-\lambda}\sigma(m)$.

Taking distribution sections we obtain the distribution globalization of this principal series representation which we denote by $H^{\sigma,\lambda}_{-\infty}$. If $V$ is a complex representation of $\Gamma$, then $H^*(\Gamma, V)$ denotes group cohomology of $\Gamma$ with coefficients in $V$. The specialization of Patterson's conjecture to cocompact $\Gamma$ is:

(i) $\dim H^*(\Gamma, H^{\sigma,\lambda}_{-\infty}) < \infty$,

(ii) $\chi(\Gamma, H^{\sigma,\lambda}_{-\infty}) := \sum_{i=0}^{\infty}(-1)^i \dim \, H^i(\Gamma, H^{\sigma,\lambda}_{-\infty}) = 0$, and

(iii) $-\chi_1(\Gamma, H^{\sigma,\lambda}_{-\infty}) := -\sum_{i=0}^{\infty}(-1)^i i \dim \, H^i(\Gamma, H^{\sigma,\lambda}_{-\infty})$ is the order of $Z_S(s, \sigma)$ at $s = \lambda$ (a pole has negative and a zero has positive order).

This conjecture was proved in [4] and [7] (with a slight modification for $s = 0$) by clarifying the relation between both approaches. In (i)–(iii) one can replace $H^{\sigma,\lambda}_{-\infty}$ by the space $H^{\sigma,\lambda}_{-\omega}$ of hyperfunction sections.

One way to define group cohomology is to write down an explicit complex $(C^*, d)$ such that $H^*(\Gamma, V)$ is the cohomology of this complex. One can take e.g.

$$C^p := \{f : \underbrace{\Gamma \times \ldots \Gamma}_{p+1 \text{ times}} \to V | f(gg_0, \ldots, gg_p) = gf(g_0, \ldots, g_p)\}$$

and

$$(df)(g_0, \ldots, g_{p+1}) := \sum_{i=0}^{p+1} (-1)^i f(g_0, \ldots, \check{g}_i, \ldots, g_{p+1}) \ .$$

Alternatively one can define group cohomology as the right derived functor of the left exact functor from the category of complex representations of $\Gamma$ to complex vector spaces which takes in each representation the subspace of $\Gamma$-invariant vectors. By homological algebra one can compute group cohomology using acyclic resolutions. To find workable acyclic resolutions for the representations of interest is one of the main goals of the present paper.

Let now $\Gamma \subset G$ be convex cocompact such that $Y$ is a noncompact manifold. The space $G/P$ can be identified with the geodesic boundary $\partial X$ of $X$. There is a $\Gamma$-invariant partition $\partial X = \Omega \cup \Lambda$, where $\Lambda$ is the limit set, and $\Omega \neq \emptyset$ is the domain of discontinuity for $\Gamma$ with compact quotient $B := \Gamma \backslash \Omega$. According to the conjecture of Patterson for the convex cocompact case one should replace $H^{\sigma,\lambda}_{-\infty}$ by the $\Gamma$-submodule of distribution sections of $V(\sigma_\lambda)$ with support on the limit set $\Lambda$.

The main aim of this paper is to prove the conjecture of Patterson for convex cocompact $\Gamma \subset \mathrm{SO}(1,n)_0$ and trivial $\sigma$ up to two modifications which we will now describe.

The first modification is that we consider hyperfunctions instead of distributions. From the technical point of view hyperfunctions are more natural and easier to handle. In fact, in several places we use the flabbiness of the sheaf of hyperfunctions. Hyperfunctions also appear in a natural way as boundary values of eigenfunctions of $\Delta_X$. In order to obtain distribution boundary values one would have to require growth conditions. Guided by the experience with cocompact groups and by the fact that the spaces of $\Gamma$-invariant hyperfunctions and distributions with support on $\Lambda$ coincide, we believe that the cohomology groups are insensitive to replacing hyperfunctions by distributions (note that the situation is different if $\Gamma$ has parabolic elements; see [7]).

The second modification is in fact already necessary in the cocompact case to get things right at the point $s = 0$. Note that the principal series representations $H^{\sigma,\lambda}$ come as a holomorphic family parametrized by $\lambda \in \mathbf{C}$. In the conjecture we replace $H^{\sigma,\lambda}_{-\omega}$ by the representation of $\Gamma$ on the space of Taylor series of length $k > 0$ at $\lambda$ of holomorphic families $\mathbf{C} \ni \mu \mapsto f_\mu \in H^{\sigma,\mu}_{-\omega}$ such that $\mathrm{supp}(f_\mu) \subset \Lambda$ for all $\mu$ (this is the space $\mathcal{O}_{(\lambda,k)}C^{-\omega}(\Lambda)$ below).

### 1.4. *The main result.*

In the present paper we prove the conjecture of Patterson for $\Gamma \subset SO(1, n)_0$ a convex cocompact, non-cocompact and torsion-free subgroup and the trivial representation $\sigma$ of $M$. We restrict ourselves to this special case mainly because of the lack of information about the Selberg zeta function in the other cases.

Let us first define $\mathcal{O}_{(\lambda,k)}C^{-\omega}(\Lambda)$. The group $G$ acts on the geodesic boundary $\partial X \cong S^{n-1}$ by means of conformal automorphisms. Let $\partial X = \Omega \cup \Lambda$ be the decomposition of $\partial X$ into the limit set $\Lambda$ and the domain of discontinuity $\Omega$.

For any $\lambda \in \mathbf{C}$ let $V_\lambda$ be the representation of $P$ on $\mathbf{C}$ given by $man \mapsto a^{\rho-\lambda} := \mathrm{e}^{(\rho-\lambda)\log(a)}$. Let $V(\lambda) := G \times_P V_\lambda$ be the associated homogeneous line bundle. Note that $V(-\rho) \cong \Lambda^{n-1}T^*_{\mathbf{C}}\partial X$ is the complexified bundle of volume forms. Moreover, $V(\lambda) \cong (\Lambda^{n-1}T^*_{\mathbf{C}}\partial X)^{\frac{n-1-2\lambda}{2(n-1)}}$. If we choose a nowhere-vanishing volume form vol on $\partial X$, then $\mathrm{vol}^{\frac{n-1-2\lambda}{2(n-1)}}$ is a section trivializing $V(\lambda)$. Sections of $V(\lambda)$ can thus be viewed as functions which transform under $G$ according to a conformal weight related to $\lambda$.

The union $\bigcup_{\lambda \in \mathbf{C}} V(\lambda) \to \partial X$ has the structure of a holomorphic family of line bundles. Using the nowhere-vanishing volume form vol we define isomorphisms $\mathrm{vol}^{-\frac{\mu}{n-1}} : V(\lambda) \cong V(\lambda + \mu)$. Hence we can identify the space of sections of $V(\lambda)$ of a given regularity with the corresponding space of sections of a fixed bundle, e.g., of the trivial one $V(\rho)$. This allows us to speak of holomorphic families of sections or homomorphisms.

By $\pi^\lambda(g) : C^{-\omega}(\partial X, V(\lambda)) \to C^{-\omega}(\partial X, V(\lambda))$, $g \in G$, we denote the representation of $G$ on the space of hyperfunction sections of $V(\lambda)$. As a topological vector space $C^{-\omega}(\partial X, V(\lambda))$ is the space of continuous linear functionals on $C^\omega(\partial X, V(-\lambda))$. If $g \in G$ is fixed, then $\pi^\lambda(g)$ depends holomorphically on $\lambda$.

Since $\Lambda$ is $\Gamma$-invariant the space of hyperfunctions $C^{-\omega}(\Lambda, V(\lambda)) \subset C^{-\omega}(\partial X, V(\lambda))$ with support in $\Lambda$ carries a representation of $\Gamma$ induced by $\pi^\lambda$. Let $\mathcal{O}_\lambda C^{-\omega}(\Lambda)$ denote the space of germs at $\lambda$ of holomorphic families of sections $\mathbf{C} \ni \mu \mapsto f_\mu \in C^{-\omega}(\partial X, V(\mu))$ with $\mathrm{supp}(f_\mu) \subset \Lambda$. The representation of $\Gamma$ on that space is given by $(\pi(g)f)_\mu := \pi^\mu(g)f_\mu$, $g \in \Gamma$.

If $L_\lambda^k$ denotes the multiplication operator $L_\lambda^k : f_\mu \mapsto (\mu - \lambda)^k f_\mu$, $k \in \mathbf{N}$, then we have a short exact sequence

$$0 \to \mathcal{O}_\lambda C^{-\omega}(\Lambda) \overset{L_\lambda^k}{\to} \mathcal{O}_\lambda C^{-\omega}(\Lambda) \to \mathcal{O}_{(\lambda,k)}C^{-\omega}(\Lambda) \to 0$$

of $\Gamma$-modules defining $\mathcal{O}_{(\lambda,k)}C^{-\omega}(\Lambda)$. Note that $\mathcal{O}_{(\lambda,1)}C^{-\omega}(\Lambda) \cong C^{-\omega}(\Lambda, V(\lambda))$. Now we can formulate the main theorem of our paper.

THEOREM 1.3.   *For any* $\lambda \in \mathbf{C}$ *there is* $k(\lambda) \in \mathbf{N}_0$ *such that the following assertions hold*:

(i)   $\dim H^*(\Gamma, \mathcal{O}_{(\lambda,k)} C^{-\omega}(\Lambda)) < \infty$ *for all* $k$, $\dim H^*(\Gamma, \mathcal{O}_\lambda C^{-\omega}(\Lambda)) < \infty$.

(ii)   $\chi(\Gamma, \mathcal{O}_{(\lambda,k)} C^{-\omega}(\Lambda)) = 0$ *for all* $k$.

(iii)   *If* $k \geq k(\lambda)$, *then* $\dim H^*(\Gamma, \mathcal{O}_{(\lambda,k+1)} C^{-\omega}(\Lambda)) = \dim H^*(\Gamma, \mathcal{O}_{(\lambda,k)} C^{-\omega}(\Lambda))$.

(iv)   *If* $k \geq k(\lambda)$, *then the order of the Selberg zeta function at* $\lambda$ *is given by*

$$\begin{align}
(3) \qquad \text{ord}_{s=\lambda} Z_S(s) &= -\chi(\Gamma, \mathcal{O}_\lambda C^{-\omega}(\Lambda)) \\
(4) \qquad &= -\chi_1(\Gamma, \mathcal{O}_{(\lambda,k)} C^{-\omega}(\Lambda)) \ ,
\end{align}$$

*where for any* $\Gamma$-*module* $V$ *with* $\dim H^*(\Gamma, V) < \infty$ *its first derived Euler characteristic* $\chi_1(\Gamma, V)$ *is defined by*

$$\chi_1(\Gamma, V) := \sum_{p=1}^{n} p(-1)^p \dim H^p(\Gamma, V) \ .$$

It will be shown in Proposition 4.19 that one can take $k(\lambda) := \text{Ord}_{\mu=\lambda} \text{ext}_\mu + \varepsilon$, where $\varepsilon = 0$ if $\lambda \notin -\mathbf{N}_0 - \rho$, and $\epsilon = 1$ otherwise, and also where $\text{ext}_\mu$ is the extension map explained in the next subsection. In contrast to ord, $\text{Ord}_{\mu=\lambda}$ denotes the (positive) order of a pole at $\mu = \lambda$, if there is one, and it is zero otherwise.

For generic $\Gamma$ and most $\lambda$ one expects $k(\lambda) \leq 1$. For those $\lambda$ one can replace $\mathcal{O}_{(\lambda,k)} C^{-\omega}(\Lambda)$ by $C^{-\omega}(\Lambda, V(\lambda))$ and probably also by $C^{-\infty}(\Lambda, V(\lambda))$, the space appearing in Patterson's original conjecture.

In [38] the order of $Z_S(s)$ at $s = 0$ was not given explicitly. As a corollary of our computations we obtain:

COROLLARY 1.4.

$$\text{ord}_{s=0} Z_S(s) = \dim {}^\Gamma C^{-\omega}(\Lambda, V(0)) = \dim \ker(S_0 + \text{id}) \ ,$$

*where* $S_0$ *is the normalized scattering matrix* (*introduced in Section* 2) *at zero*.

The proof of Theorem 1.3 consists of three steps. The first step which occupies most of the paper is an explicit computation of the cohomology groups $H^*(\Gamma, \mathcal{O}_{(\lambda,k)} C^{-\omega}(\Lambda))$, $H^*(\Gamma, \mathcal{O}_\lambda C^{-\omega}(\Lambda))$. One of our main tools in these computations is the extension map which will be explained in the next subsection. The final results which are of interest in their own right, and much more detailed than needed for the proof of Theorem 1.3, are stated in subsections 4.3 and 4.4. They are generalizations of our results for the two-dimensional case [6].

In a second step we compare the result of this computation with the spectral description of the singularities of $Z_S$ given by Patterson-Perry [38] in case $n \equiv 0(2)$, $\delta_\Gamma < 0$. Finally we employ the embedding trick in order to drop this assumption. The second and third steps are performed in Section 5.

In order to compute $H^*(\Gamma, \mathcal{O}_\lambda C^{-\omega}(\Lambda))$ we use suitable acyclic resolutions by $\Gamma$-modules formed by germs of holomorphic families of hyperfunctions on $\partial X$ and $\Omega$. The proof of exactness and acyclicity of these resolutions is quite involved and uses some hyperfunction theory and analysis on the symmetric space $X$. It turned out that we need facts which hold true in much more general situations but have not been considered in the literature so far (up to our knowledge). This accounts for the length of subsections 4.1 and 4.2.

### 1.5. *The extension map.*

The present paper has a close companion [8] in which we consider the decomposition of the right regular representation of $G$ on $L^2(\Gamma \backslash G)$ into irreducible unitary representations in the case that $\Gamma$ is a convex cocompact subgroup of a simple Lie group of real rank one. The main ingredient of both papers is the extension map $\text{ext}_\lambda$.

Consider any $\Gamma$-invariant hyperfunction section $f \in {}^\Gamma C^{-\omega}(\Omega, V(\lambda))$. By the flabbiness of the sheaf of hyperfunction sections it can be extended across $\Lambda$; i.e., there is a hyperfunction section $\tilde{f} \in C^{-\omega}(\partial X, V(\lambda))$ which restricts to $f$. In general $\tilde{f}$ is not $\Gamma$-invariant. Our extension map solves the problem of finding a $\Gamma$-invariant extension $\text{ext}_\lambda(f)$. It turns out that the invariant extension exists and is unique for generic $\lambda$.

Now the maps $\text{ext}_\lambda$ form in fact a meromorphic family of maps with finite-dimensional singularities. The highest singular part of ext at the poles of ext gives invariant hyperfunction (in fact distribution) sections of $V(\lambda)$ with support on $\Lambda$. This can be considered as a generalization of the construction of the Patterson-Sullivan measure which is given by the residue of ext at $\lambda = \delta_\Gamma$.

In Section 3 we employ a version of Green's formula in order to get a hold on the spaces of invariant hyperfunction sections with support on $\Lambda$. In particular, it follows that there is a discrete set of $\lambda \in \mathbf{C}$, where $f \in {}^\Gamma C^{-\omega}(\Omega, V(\lambda))$ has to satisfy a finite number of nontrivial linear conditions in order to be invariantly extendable. To find a $\Gamma$-invariant extension of $f$ is a cohomological problem and $H^1(\Gamma, C^{-\omega}(\Lambda, V(\lambda))$ essentially appears as its obstruction group. This is the basic observation which enables us to compute these cohomology groups.

In order to provide a feeling for the extension problem let us discuss a toy example. We consider the extension of hyperfunctions $f$ on $\mathbf{R} \setminus \{0\}$, which transform as $f(rx) = r^\lambda f(x)$ for all $r \in \mathbf{R}_+^*$. Let $f_\lambda$ be given by $f_\lambda(x) := |x|^\lambda$. For $\text{Re}(\lambda) > -1$ its invariant extension as a distribution (and hence as a hyperfunction) to $\mathbf{R}$ is just given by

$$\langle \mathrm{ext}_\lambda(f_\lambda), \phi \rangle := \int_{\mathbf{R}} \phi(x)|x|^\lambda dx \ , \quad \phi \in C_c^\infty(\mathbf{R}) \ .$$

It is well-known that $\mathrm{ext}_\lambda(f_\lambda)$ has a meromorphic continuation to all of $\mathbf{C}$ with poles at negative integers. The residues of the continuation at these points are proportional to derivatives of delta distributions located at $\{0\}$.

The construction of the meromorphic continuation of the extension map is closely related to the meromorphic continuation of the scattering matrix

$$S_\lambda : {}^\Gamma C^{-\omega}(\Omega, V(\lambda)) \to {}^\Gamma C^{-\omega}(\Omega, V(-\lambda)) \ .$$

If $\Gamma$ is the trivial subgroup, then $\Omega = \partial X$ and $S_\lambda$ coincides with the (normalized) Knapp-Stein intertwining operator

$$J_\lambda : C^{-\omega}(\partial X, V(\lambda)) \to C^{-\omega}(\partial X, V(-\lambda)) \ .$$

The operators $J_\lambda$ form a meromorphic family and are well-studied in representation theory. In this paper and in [8] we approach the scattering matrix $S_\lambda$ starting from $J_\lambda$ and using the basic identity

$$S_\lambda = \mathrm{res}_{-\lambda} \circ J_\lambda \circ \mathrm{ext}_\lambda \ ,$$

where $\mathrm{res}_\lambda$ denotes the restriction $\mathrm{res}_\lambda : {}^\Gamma C^{-\omega}(\partial X, V(\lambda)) \to {}^\Gamma C^{-\omega}(\Omega, V(\lambda))$.

In the literature the scattering matrix is usually considered as a certain pseudodifferential operator which can be applied to smooth, resp. distribution sections. A meromorphic continuation of $S_\lambda$ was obtained in [35] for surfaces, and in [40], [30] in general. In [8] we develop the theory of the scattering matrix and the extension map in a smooth/distribution framework for general rank-one spaces and arbitrary $\sigma$. The main point in the present paper is the transition to the real analytic/hyperfunction framework.

The main difference from most previous papers is that our primary analysis concerns objects on the boundary $\partial X$. The spectral theory of $\Delta_Y$ is only needed in some very weak form. One the other hand one can deduce the spectral decomposition, the meromorphic continuation of Eisenstein series and the properties of the resolvent kernel from the theory on the boundary. To illustrate this consider, e.g., the Eisenstein series. Let $P_\lambda : C^{-\omega}(\partial X, V(\lambda)) \to C^\infty(X)$ denote the Poisson transform. For $b \in \Omega$ we define $f_b \in {}^\Gamma C^{-\omega}(\Omega, V(\lambda))$ by $f_b := \sum_{\gamma \in \Gamma} \pi^\lambda(\gamma)(\delta_b \mathrm{vol}^\mu)$, where $\delta_b \in C^{-\omega}(\Omega, V(-\rho))$ is the delta distribution located at $b$ and $\mu := -\frac{n-1+2\lambda}{2(n-1)}$. Then the Eisenstein series can be written as

$$E_\lambda(x, b) := (P_\lambda \circ \mathrm{ext}_\lambda(f_b))(x) \ .$$

These applications will be contained in [8] and its continuations.

*Acknowledgement.* We thank S. Patterson and P. Perry for keeping us informed about the progress of [38]. Moreover we thank A. Juhl for pointing out some wrong arguments in a previous version of this paper and for further useful remarks. This work was partially supported by the Sonderforschungsbereich 288, *Differentialgeometrie und Quantenphysik*.

## 2. Restriction, extension, and the scattering matrix

### 2.1. *Basic notions.*

The sheaf of hyperfunction sections of a real analytic vector bundle over a real analytic manifold is flabby. Thus the following sequence of $\Gamma$-modules

$$0 \to C^{-\omega}(\Lambda, V(\lambda)) \to C^{-\omega}(\partial X, V(\lambda)) \overset{\mathrm{res}_\Omega}{\to} C^{-\omega}(\Omega, V(\lambda)) \to 0$$

is exact, where $\mathrm{res}_\Omega$ is the restriction of sections to $\Omega$.

Let $V_B(\lambda) := \Gamma \backslash V(\lambda)_{|\Omega}$. If we identify $^\Gamma C^{-\omega}(\Omega, V(\lambda)) \cong C^{-\omega}(B, V_B(\lambda))$, then $\mathrm{res}_\Omega$ induces a map $\mathrm{res}_\lambda : {}^\Gamma C^{-\omega}(\partial X, V(\lambda)) \to C^{-\omega}(B, V_B(\lambda))$. Here $^\Gamma C^{-\omega}(., V(\lambda))$ denotes the subspace of $\Gamma$-invariant sections.

The main topic of this section is the construction of a meromorphic family of maps $\mathrm{ext}_\lambda : C^{-\omega}(B, V_B(\lambda)) \to {}^\Gamma C^{-\omega}(\partial X, V(\lambda))$ which are right inverse to $\mathrm{res}_\lambda$. The poles of $\mathrm{ext}_\lambda$ will correspond exactly to those points $\lambda \in \mathbf{C}$ where $\mathrm{res}_\lambda$ fails to be an isomorphism.

The strategy of the construction of $\mathrm{ext}_\lambda$ is the following. We first construct $\mathrm{ext}_\lambda$ for $\mathrm{Re}(\lambda) > \delta_\Gamma$. Then we introduce the scattering matrix $\hat{S}_\lambda$ : $C^{-\omega}(B, V_B(\lambda)) \to C^{-\omega}(B, V_B(-\lambda))$ by

$$(5) \qquad \hat{S}_\lambda := \mathrm{res}_{-\lambda} \circ \hat{J}_\lambda \circ \mathrm{ext}_\lambda \ ,$$

where $\hat{J}_\lambda : C^{-\omega}(\partial X, V(\lambda)) \to C^{-\omega}(\partial X, V(-\lambda))$ is the Knapp-Stein intertwining operator (see [26]) which we will introduce below. Assuming for a moment that $\delta_\Gamma < 0$ we obtain a meromorphic continuation of $\hat{S}_\lambda$ using results of Patterson [34] and [8]. Then we construct the meromorphic continuation of $\mathrm{ext}_\lambda$ by

$$(6) \qquad \mathrm{ext}_\lambda := J_{-\lambda} \circ \mathrm{ext}_{-\lambda} \circ S_\lambda, \quad \mathrm{Re}(\lambda) < 0 \ ,$$

where $J_{-\lambda}$ ($S_\lambda$) is the normalized intertwining operator (scattering matrix).

If $\delta_\Gamma \geq 0$, then we employ the embedding $SO(1, n)_0 \hookrightarrow SO(1, n + m)_0$, $m$ sufficiently large, in order to reduce to the case $\delta_\Gamma < 0$.

### 2.2. *Holomorphic functions to topological vector spaces.*

In order to carry out the program sketched above we need to consider holomorphic families of vectors in topological vector spaces. We will also consider holomorphic families of continuous linear maps between such vector spaces. Therefore we collect some preparatory material of a technical nature. All topological vector spaces appearing in this paper are Hausdorff and complete.

A holomorphic family of vectors in a locally convex topological vector space $\mathcal{F}$ defined on $U \subset \mathbf{C}$ is by definition a continuous function from $U$ to $\mathcal{F}$ which is weakly holomorphic. Using Cauchy's integral formula one can show an equivalent characterization of holomorphic families. A map $f : U \to \mathcal{F}$ is holomorphic, if and only if for any $z_0 \in U$ there is a neighbourhood $z_0 \in V \subset U$ and a sequence $\{f_k \in \mathcal{F}\}_{k \in \mathbf{N}_0}$ such that for $z \in V$ the sum $\sum_{k=0}^{\infty} f_k(z - z_0)^k$ converges and is equal to $f(z)$.

In order to speak of holomorphic families of homomorphisms from $\mathcal{F}$ to $\mathcal{G}$, where $\mathcal{G}$ is another locally convex topological vector space we equip $\mathrm{Hom}(\mathcal{F}, \mathcal{G})$ with the topology of uniform convergence on bounded sets.

Let $f : U \setminus \{z_0\} \to \mathrm{Hom}(\mathcal{F}, \mathcal{G})$ be holomorphic and $f(z) = \sum_{k=-N}^{\infty} f_k(z - z_0)^k$ for all $z \neq z_0$ close to $z_0$. Then we say that $f$ is meromorphic and has a pole of order $N$ at $z_0$. If $f_k$, $k = -N, \ldots, -1$, are finite-dimensional, then, by definition, $f$ has a finite-dimensional singularity.

Holomorphy of a map $f : U \to \mathrm{Hom}(\mathcal{F}, \mathcal{G})$ can be characterized in the following weak form. Let $\mathcal{G}'$ denote the dual space of $\mathcal{G}$ with its strong topology. We call a subset $A \subset \mathcal{F} \times \mathcal{G}'$ sufficiently large if for $B \in \mathrm{Hom}(\mathcal{F}, \mathcal{G})$ the condition $\langle \phi, B\psi \rangle = 0$, for all $(\psi, \phi) \in A$, implies $B = 0$.

LEMMA 2.1. *The following assertions are equivalent*:

(i) $f : U \to \mathrm{Hom}(\mathcal{F}, \mathcal{G})$ *is holomorphic.*

(ii) $f : U \to \mathrm{Hom}(\mathcal{F}, \mathcal{G})$ *is continuous, and there is a sufficiently large set $A \subset \mathcal{F} \times \mathcal{G}'$ such that for all $(\psi, \phi) \in A$ the function $U \ni z \mapsto \langle \phi, f(z)\psi \rangle$ is holomorphic.*

*Proof.* It is obvious that (i) implies (ii). We show that (i) follows from (ii). It is easy to see that $f$ is holomorphic at $z$ if and only if

$$f(z') = \frac{1}{2\pi\iota} \oint \frac{f(x)}{x - z'} dx$$

for $z'$ close to $z$, where the path of integration is a small circle surrounding $z$ counterclockwise. Let $f$ satisfy (ii). Then we form

$$d(z') := f(z') - \frac{1}{2\pi\iota} \oint \frac{f(x)}{x - z'} dx \ .$$

We must show that $d = 0$. It suffices to show that for all $(\psi, \phi) \in A$ we have $\langle \phi, d(z')\psi \rangle = 0$. But this follows from (ii) and Cauchy's integral formula. $\square$

LEMMA 2.2. *Let $f_i : U \to \mathrm{Hom}(\mathcal{F}, \mathcal{G})$ be a sequence of holomorphic maps. Moreover let $f : V \to \mathrm{Hom}(\mathcal{F}, \mathcal{G})$ be continuous such that for a sufficiently large set $A \subset \mathcal{F} \times \mathcal{G}'$ the functions $\langle \phi, f_i \psi \rangle$ converge locally uniformly in $U$ to $\langle \phi, f\psi \rangle$. Then $f$ is holomorphic, too.*

*Proof.* Since the holomorphic functions $\langle \phi, f_i \psi \rangle$ converge locally uniformly to $\langle \phi, f\psi \rangle$ we conclude that the latter function is holomorphic, too. The lemma is now a consequence of Lemma 2.1. $\qquad \square$

LEMMA 2.3. *Let* $f : U \to \mathrm{Hom}(\mathcal{F}, \mathcal{G})$ *be continuous. Then the adjoint* $f' : U \to \mathrm{Hom}(\mathcal{G}', \mathcal{F}')$ *is continuous. If* $f$ *is holomorphic, then so is* $f'$.

*Proof.* We first show that the adjoint $f' : U \to \mathrm{Hom}(\mathcal{G}', \mathcal{F}')$ is continuous at $z_0 \in U$. Let $B \subset \mathcal{G}'$ be a bounded set. Let $q$ be any continuous seminorm on $\mathcal{F}'$. Then we have to show that for any $\epsilon > 0$ there exists $\delta > 0$ such that if $|z - z_0| < \delta$, then $\sup_{\phi \in B} q(f'(z)\phi - f'(z_0)\phi) < \epsilon$. The strong topology of $\mathcal{F}'$ is generated by the seminorms $q_D$ associated to bounded subsets $D \subset \mathcal{F}$, where $q_D(\psi) := \sup_{\kappa \in D} |\langle \psi, \kappa \rangle|$. We have

$$
\begin{aligned}
\sup_{\phi \in B} q_D(f'(z)\phi - f'(z_0)\phi) &= \sup_{\phi \in B, \kappa \in D} |\langle f'(z)\phi - f'(z_0)\phi, \kappa \rangle| \\
&= \sup_{\phi \in B, \kappa \in D} |\langle \phi, f(z)\kappa - f(z_0)\kappa \rangle| \\
&= \sup_{\kappa \in D} q_B(f(z)\kappa - f(z_0)\kappa) \ .
\end{aligned}
$$

Here $q_B$ is a continuous seminorm on the bidual $\mathcal{G}''$. Since the embedding $\mathcal{G} \hookrightarrow \mathcal{G}''$ is continuous and $f$ is continuous at $z_0$ we can find $\delta > 0$ for any $\epsilon > 0$ as required.

If $f$ is holomorphic, then holomorphy of $f'$ follows from Lemma 2.1 when we take the sufficiently large set $\mathcal{F} \times \mathcal{G}'$, and use the fact that $\langle \phi, f\psi \rangle = \langle f'\phi, \psi \rangle$ is holomorphic for all $\psi \in \mathcal{F}$, $\phi \in \mathcal{G}'$. $\qquad \square$

A locally convex vector space is called Montel if its closed bounded subsets are compact.

LEMMA 2.4. *Let* $\mathcal{F}, \mathcal{G}, \mathcal{H}$ *be locally convex topological vector spaces and assume that* $\mathcal{F}$ *is a Montel space. If* $f : U \to \mathrm{Hom}(\mathcal{F}, \mathcal{G})$ *and* $f_1 : U \to \mathrm{Hom}(\mathcal{G}, \mathcal{H})$ *are continuous, then the composition* $f_1 \circ f : V \to \mathrm{Hom}(\mathcal{F}, \mathcal{H})$ *is continuous. If* $f$ *and* $f_1$ *are holomorphic, then so is the composition* $f_1 \circ f$.

*Proof.* We first prove continuity of the composition $f_1 \circ f$ at $z_0 \in U$. It is here where we need the assumption that $\mathcal{F}$ is Montel. Let $B \subset \mathcal{F}$ be a bounded set and $s$ be a seminorm of $\mathcal{H}$. We have to show that for all $\epsilon > 0$ there exists $\delta > 0$ such that if $|z - z_0| < \delta$, then $\sup_{\phi \in B} s((f_1 \circ f)(z)(\phi) - (f_1 \circ f)(z_0)(\phi)) < \epsilon$.

Let $W \subset U$ be a compact neighbourhood of $z_0$. The map $F : W \times B \to \mathcal{G}$ given by $F(v, \phi) := f(v)(\phi)$ is continuous. Since $\mathcal{F}$ is Montel, any bounded set $B$ is precompact. Thus the image $B_1 := F(W \times B)$ of the compact set $W \times B$ is precompact, too. In particular, $B_1$ is bounded. Since $f_1(z_0)$ is continuous there is a seminorm $t$ of $\mathcal{G}$ such that $s(f_1(z_0)(\psi)) < t(\psi)$, for all $\psi \in \mathcal{G}$. Now let $V \subset W$ be a neighbourhood of $z_0$ so small that $t(f(z)(\phi) - f(z_0)(\phi)) < \epsilon/2$, for all $z \in V$, for all $\phi \in B$, and $s(f_1(z)(\psi) - f_1(z_0)(\psi)) < \epsilon/2$, for all $\psi \in B_1$, for all $z \in V$. Then for all $z \in V$, $\phi \in B$,

$$
\begin{aligned}
s(f_1(z) \circ f(z)(\phi) - f_1(z_0) \circ f(z_0)(\phi)) \quad \leq \quad & s([f_1(z) - f_1(z_0)]f(z)(\phi)) \\
& + s(f_1(z_0)[f(z)(\phi) - f(z_0)(\phi)]) \\
< \quad & \epsilon \ .
\end{aligned}
$$

Thus we can find $\delta > 0$ for any $\epsilon > 0$ as required. This proves continuity of the composition $f_1 \circ f$.

We now show that this composition is holomorphic if $f, f_1$ are so. We have

$$
\frac{1}{2\pi \imath} \oint \frac{f_1(x) \circ f(x)}{x - z} dx = \frac{1}{(2\pi \imath)^2} \oint \oint \frac{f_1(x) \circ f(y)}{(x-z)(y-x)} dy \, dx \ .
$$

If we restrict this equation to a bounded set $B \subset \mathcal{F}$, we see as above that there exists some bounded set $B_1 \subset \mathcal{G}$ such that $\frac{f(y)}{(x-z)(y-x)}(B) \subset B_1$ for all $y, x$ in the domain of integration. Hence we can apply Fubini's theorem to the double integral and obtain

$$
\begin{aligned}
\frac{1}{2\pi \imath} \oint \frac{f_1(x) \circ f(x)}{x - z} dx \quad = \quad & \frac{1}{(2\pi \imath)^2} \oint \oint \frac{f_1(x) \circ f(y)}{(x-z)(y-x)} dx \, dy \\
= \quad & \frac{1}{2\pi \imath} \oint \frac{f_1(z) \circ f(y)}{y - z} dy \\
= \quad & (f_1 \circ f)(z) \ .
\end{aligned}
$$

This shows that $f_1 \circ f$ is holomorphic. □

Consider a real analytic vector bundle over a closed real analytic manifold. Then the spaces of real analytic, smooth, distribution, and hyperfunction sections of the bundle equipped with their natural locally convex topologies are Montel spaces.

### 2.3. *The push-down.*

In the present subsection we define a push-down map

$$
\pi_{*,-\lambda} : C^{\sharp}(\partial X, V(-\lambda)) \to C^{\sharp}(B, V_B(-\lambda)), \quad \sharp \in \{\omega, \infty\}, \ \mathrm{Re}(\lambda) > \delta_\Gamma \ .
$$

The extension $\mathrm{ext}_\lambda$ will then appear as its adjoint.

Using the identification $C^{\sharp}(B, V_B(-\lambda)) \cong {}^{\Gamma} C^{\sharp}(\Omega, V(-\lambda))$ we want to define $\pi_{*,-\lambda}$ by

$$
(7) \qquad \pi_{*,-\lambda}(f)(b) = \sum_{g \in \Gamma} (\pi^{-\lambda}(g)f)(b), \quad b \in \Omega, \quad f \in C^{\sharp}(\partial X, V(-\lambda))
$$

provided the sum converges. In order to prove the convergence for $\mathrm{Re}(\lambda) > \delta_\Gamma$ we need the following two geometric lemmas.

We adopt the following conventions about the notation for points of $X$ and $\partial X$. A point $x \in \partial X$ can equivalently be denoted by a subset $kM \subset K$ or $gP \subset G$ representing this point in $\partial X = K/M$ or $\partial X = G/P$. If $F \subset \partial X$, then $FM := \bigcup_{kM \in F} kM \subset K$. Analogously, we can denote a point $b \in X$ by a set $gK \subset G$, where $gK$ represents $b$ in $X = G/K$.

Adjoining the boundary at infinity we can consider $X \cup \partial X$ as a compact manifold with boundary carrying a smooth action of $G$. Let $\Gamma \subset G$ be a torsion-free convex cocompact subgroup. An equivalent characterization of being convex cocompact is that $\Gamma$ acts freely and cocompactly on $X \cup \Omega$.

LEMMA 2.5. *If $F \subset \Omega$ is compact, then $\sharp(\Gamma \cap (FM)A_+K) < \infty$.*

*Proof.* Note that $(FM)A_+K \cup F \subset X \cup \Omega$ is compact. Thus its intersection with the orbit $\Gamma K$ of the origin of $X$ is finite. □

Using the Iwasawa decomposition $G = KAN$ we write $g \in G$ as $g = \kappa(g)a(g)n(g) \in KAN$. The Iwasawa decomposition and, in particular, the maps $g \mapsto \kappa(g)$, $g \mapsto a(g)$, $g \mapsto n(g)$ are real analytic and extend to a complex neighbourhood of $G$ in its complexification $G_{\mathbf{C}}$. Let $K_{\mathbf{C}}, A_{\mathbf{C}}$ be the complexifications of $K, A$. We identify $A_{\mathbf{C}}$ with the multiplicative group $\mathbf{C}^*$ such that $A_+$ corresponds to $[1, \infty) \subset \mathbf{C}^*$. Any $g \in G$ has a Cartan decomposition $g = ha_gh' \in KA_+K$, where $a_g$ is uniquely determined.

The next lemma gives some control on the complex extension of the Iwasawa decomposition.

LEMMA 2.6. *Let $k_0M \in \partial X$. For any compact $W \subset (\partial X \setminus k_0M)M$ and complex neighbourhood $S_{\mathbf{C}} \subset K_{\mathbf{C}}$ of $K$ there are a complex neighbourhood $U_{\mathbf{C}} \subset K_{\mathbf{C}}$ of $k_0M$ and constants $c > 0$, $C < \infty$, such that for all $g = ha_gh' \in WA_+K$ the components $\kappa(g^{-1}k)$ and $a(g^{-1}k)$ extend holomorphically to $k \in U_{\mathbf{C}}$, and for all $k \in U_{\mathbf{C}}$*

(8) $$ca_g \leq |a(g^{-1}k)| \leq Ca_g,$$

(9) $$\kappa(g^{-1}k) \in S_{\mathbf{C}}.$$

*Proof.* The set $W^{-1}k_0M$ is compact and disjoint from $M$. Let $w \in N_K(M)$ represent the nontrivial element of the Weyl group of $(\mathbf{g}, \mathbf{a})$. Set $\bar{\mathbf{n}} := \theta(\mathbf{n})$, where $\theta$ is the Cartan involution of $G$ fixing $K$, and define $\bar{N} := \exp(\bar{\mathbf{n}})$. By the Bruhat decomposition $G = w\bar{N}P \cup P$ we have $K = w\kappa(\bar{N})M \cup M$. Thus there is a compact $V \subset \bar{N}$ such that $W^{-1}k_0M \subset w\kappa(V)M$. By enlarging $V$ we can assume that $V$ is $A_+$-invariant, where $A$ acts on $\bar{N}$ by $(a, \bar{n}) \to a\bar{n}a^{-1}$. There exists a complex compact $A_+$-invariant neighbourhood $V_{\mathbf{C}} \subset \bar{N}_{\mathbf{C}}$ of $V$ such that $\kappa(\bar{n}), a(\bar{n}), n(\bar{n})$ extend to $V_{\mathbf{C}}$ holomorphically. Moreover, there exists a complex neighbourhood $U^1_{\mathbf{C}} \subset K_{\mathbf{C}}$ of $k_0M$ such that $w^{-1}W^{-1}U^1_{\mathbf{C}}M \subset \kappa(V_{\mathbf{C}})M$. Let $k \in U^1_{\mathbf{C}}$ and $g = ha_gh' \in WA_+K$. Then $h^{-1}k = w\kappa(\bar{n})m$ for $\bar{n} \in V_{\mathbf{C}}$, $m \in M_{\mathbf{C}}$; i.e., we parametrize $h^{-1}U^1_{\mathbf{C}}$ by $V_{\mathbf{C}} \times M_{\mathbf{C}}$. Furthermore,

$$\begin{aligned}
a(g^{-1}k) &= a(h'^{-1}a_g^{-1}h^{-1}k) \\
&= a(a_g^{-1}w\kappa(\bar{n})m) \\
&= a(a_g\kappa(\bar{n})) \\
&= a(a_g\bar{n}n(\bar{n})^{-1}a(\bar{n})^{-1}) \\
&= a(a_g\bar{n}a_g^{-1})a(\bar{n})^{-1}a_g \ .
\end{aligned}$$

Now $a_g\bar{n}a_g^{-1} \in V_{\mathbf{C}}$. Thus $a(g^{-1}k)$ extends holomorphically to $k \in U_{\mathbf{C}}^1$. Set

$$\begin{aligned}
c &:= \inf_{\bar{n}\in V_{\mathbf{C}}}|a(\bar{n})| \inf_{\bar{n}\in V_{\mathbf{C}}}|a(\bar{n})^{-1}| \\
C &:= \sup_{\bar{n}\in V_{\mathbf{C}}}|a(\bar{n})| \sup_{\bar{n}\in V_{\mathbf{C}}}|a(\bar{n})^{-1}| \ .
\end{aligned}$$

Since $V_{\mathbf{C}}$ is compact we have $0 < c \le C < \infty$. Then

$$ca_g \le |a(g^{-1}k)| \le Ca_g \ .$$

Now considering $\kappa$ we have

$$\begin{aligned}
\kappa(g^{-1}k) &= \kappa(h'^{-1}a_g^{-1}h^{-1}k) \\
&= h'^{-1}\kappa(a_g^{-1}w\kappa(\bar{n}))m \\
&= h'^{-1}w\kappa(a_g\bar{n}n(\bar{n})^{-1}a(\bar{n})^{-1})m \\
&= h'^{-1}w\kappa(a_g\bar{n}a_g^{-1})m \ .
\end{aligned}$$

Since $a_g\bar{n}a_g^{-1} \in V_{\mathbf{C}}$ we see that $\kappa(g^{-1}k)$ extends holomorphically to $k \in U_C^1$. If we take $V_{\mathbf{C}}$ small enough we can also satisfy $h'^{-1}w\kappa(V_{\mathbf{C}}) \subset S_{\mathbf{C}}$. Thus for a smaller open subset $U_{\mathbf{C}} \subset U_{\mathbf{C}}^1$ we have $\kappa(g^{-1}k) \in S_{\mathbf{C}}$ for all $g \in WA_+K$ and $k \in U_{\mathbf{C}}$. $\qquad\square$

LEMMA 2.7. *If* $\operatorname{Re}(\lambda) > \delta_\Gamma$, *then the sum* (7) *converges and defines a holomorphic family of continuous maps*

$$\pi_{*,-\lambda} : C^\sharp(\partial X, V(-\lambda)) \to C^\sharp(B, V_B(-\lambda)), \quad \sharp \in \{\infty, \omega\} \ .$$

*Proof.* In case $\sharp = \infty$ the lemma was proved in [8]. Thus we assume $\sharp = \omega$. First we recall the definition of the topology on the spaces of real analytic sections of real analytic vector bundles.

We describe $C^\omega(\partial X)$ as a direct limit of Banach spaces. Let $S_{\mathbf{C}} \subset K_{\mathbf{C}}$ be a compact right $M$-invariant complex neighbourhood of $K$. Let $\mathcal{H}(S_{\mathbf{C}})$ denote the Banach space of bounded holomorphic functions on $S_{\mathbf{C}}$ equipped with the norm

$$\|f\| = \sup_{k\in S_{\mathbf{C}}}|f(k)|, \quad f \in \mathcal{H}(S_{\mathbf{C}}) \ .$$

If $S'_{\mathbf{C}} \subset S_{\mathbf{C}}$, then we have a continuous restriction $\mathcal{H}(S_{\mathbf{C}}) \hookrightarrow \mathcal{H}(S'_{\mathbf{C}})$. Then

$$C^\omega(K) = \varinjlim \mathcal{H}(S_{\mathbf{C}})$$

as a topological vector space, where the limit is taken over all compact right $M$-invariant complex neighbourhoods of $K$. The compact group $M$ acts continuously on $\mathcal{H}(S_{\mathbf{C}})$, and this action is compatible with the restriction to smaller $S'_{\mathbf{C}}$. Thus the natural right action of $M$ on $C^\omega(K)$ is continuous. We obtain $C^\omega(\partial X)$ as the closed subspace of right $M$-invariants in $C^\omega(K)$. Using the $K$-invariant volume form on $\partial X$ we identify $C^\omega(\partial X, V(\lambda)) \cong C^\omega(\partial X)$.

Next we describe the topological structure of $C^\omega(B, V_B(\lambda))$. Let $\{F_\alpha\}$ be a finite cover of $B$ by compact neighbourhoods which have diffeomorphic lifts $\tilde{F}_\alpha \subset \Omega$. Using again the $K$-invariant volume form we identify $C^\omega(\tilde{F}_\alpha, V(\lambda))$ with the direct limit of the Banach spaces $\mathcal{H}(U_{\mathbf{C}, \alpha})^M$ of $M$-invariant bounded holomorphic functions on $U_{\mathbf{C}, \alpha}$ where $U_{\mathbf{C}, \alpha} \subset K_{\mathbf{C}}$ are $M$-invariant compact complex neighbourhoods of $\tilde{F}_\alpha M$. Any section $f \in C^\omega(B, V_B(\lambda))$ defines sections $f_\alpha \in C^\omega(\tilde{F}_\alpha, V(\lambda))$ by restricting the lift of $f$. The correspondence $f \mapsto \bigoplus_\alpha f_\alpha$ defines an embedding

$$I : C^\omega(B, V_B(\lambda)) \hookrightarrow \bigoplus_\alpha C^\omega(\tilde{F}_\alpha, V(\lambda))$$

as a closed subspace and, hence, a topology on $C^\omega(B, V_B(\lambda))$.

The space $C^\omega(B, V_B(\lambda))$ depends on $\lambda$. Since $\Gamma$ acts on $\Omega$ preserving the orientation the manifold $B$ is orientable, and we can find an analytic real-valued nowhere-vanishing volume form $\mathrm{vol}_B \in C^\omega(B, V_B(-\rho))$. Multiplcation by suitable complex powers of $\mathrm{vol}_B$ provides isomorphisms between any two bundles $V_B(\lambda_1)$ and $V_B(\lambda_2)$. Thus we may identify $C^\omega(B, V_B(\lambda))$ with the fixed space $C^\omega(B, V_B(\rho)) = C^\omega(B)$. Using this identification we can speak of continuous or holomorphic families of continuous homomorphisms to or from $C^\omega(B, V_B(\lambda))$.

Let $\mathcal{G}$ be any locally convex topological vector space, $U \subset \mathbf{C}$ be open, and $U \ni \lambda \mapsto \phi_\lambda \in \mathrm{Hom}(\mathcal{G}, C^\omega(B, V_B(\lambda)))$ be a family of continuous maps. This family depends continuously (holomorphically) on $\lambda$, if and only if the composition $I \circ \phi_\lambda$ does.

Let $F \in \{F_\alpha\}$ and $\tilde{F} \subset \Omega$ be its lift. Let $\tilde{F} \subset U_1 \subset \Omega$ be an open neighbourhood of $\tilde{F}$ and $W := (\partial X \setminus U_1)M$. By Lemma 2.6 we can find for any complex $M$-invariant neighbourhood $S_{\mathbf{C}} \subset K_{\mathbf{C}}$ of $K$ a complex neighbourhood $U_{\mathbf{C}} \subset K_{\mathbf{C}}$ of $\tilde{F}M$ such that $a(g^{-1}k)$, $\kappa(g^{-1}k)$ extend to $k \in U_{\mathbf{C}}$ for all $g \in WA_+K$ and $\kappa(g^{-1}k) \in S_{\mathbf{C}}$. This allows us to define the map

$$\mathrm{res}_{U_{\mathbf{C}}} \circ \pi^\lambda(g) : \mathcal{H}(S_{\mathbf{C}}) \to \mathcal{H}(U_{\mathbf{C}})$$

for all $g \in WA_+K$. Here $(\pi^\lambda(g)f)(k) = a(g^{-1}k)^{\lambda-\rho}f(\kappa(g^{-1}k))$, $f \in \mathcal{H}(S_{\mathbf{C}})$, $k \in U_{\mathbf{C}}$.

In order to estimate the norm of $\mathrm{res}_{U_{\mathbf{C}}} \circ \pi^\lambda(g)$ we compute for $f \in \mathcal{H}(S_{\mathbf{C}})$ using (8),

$$
\begin{aligned}
\| \mathrm{res}_{U_{\mathbf{C}}} \circ \pi^\lambda(g) f \| &= \sup_{k \in U_{\mathbf{C}}} | \mathrm{res}_{U_{\mathbf{C}}} \circ \pi^\lambda(g)(f)(k) | \\
&= \sup_{k \in U_{\mathbf{C}}} | a(g^{-1} k)^{\lambda - \rho} f(\kappa(g^{-1} k)) | \\
&\leq C a_g^{\mathrm{Re}\lambda - \rho} \sup_{k \in S_{\mathbf{C}}} |f(k)| \\
&= C a_g^{\mathrm{Re}\lambda - \rho} \| f \| \,,
\end{aligned}
$$

where $C$ is independent of $k, g, f$.

By Lemma 2.5 we have $\sharp(\Gamma \setminus \Gamma \cap WA_+K) < \infty$. Thus for almost all $g \in \Gamma$ we have $\| \mathrm{res}_{U_{\mathbf{C}}} \circ \pi(g) \| \leq C a_g^{\mathrm{Re}\lambda - \rho}$. Fix $\epsilon > 0$. Then by the definition of $\delta_\Gamma$ the sum

$$
\mathcal{H}(S_{\mathbf{C}})^M \ni f \mapsto \sum_{g \in \Gamma \cap WA_+K} \mathrm{res}_{U_{\mathbf{C}}} \circ \pi^\lambda(g) f \in \mathcal{H}(U_{\mathbf{C}})
$$

converges uniformly for $\mathrm{Re}(\lambda) \in (-\infty, -\delta_\Gamma - \epsilon)$, and on the unit ball of the Banach space $\mathcal{H}(S_{\mathbf{C}})^M$. Going over to the direct limits we conclude that $\Pi_\lambda^F := \sum_{g \in \Gamma} \mathrm{res}_{\tilde{F}} \circ \pi^\lambda(g)$ is a convergent sum of holomorphic maps

$$
\{ \mathrm{Re}(\lambda) < -\delta_\Gamma - \epsilon \} \ni \lambda \mapsto \mathrm{res}_{\tilde{F}} \circ \pi^\lambda(g) \in \mathrm{Hom}(C^\omega(\partial X, V(\lambda)), C^\omega(\tilde{F}, V(\lambda))) \,,
$$

which is uniformly convergent on bounded subsets of $C^\omega(\partial X, V(\lambda))$. Thus $\Pi_\lambda^F$ is holomorphic with respect to $\lambda$ by Lemma 2.2. In order to finish the proof of Lemma 2.7 observe that $I \circ \pi_{*, \lambda} = \bigoplus_\alpha \Pi_\lambda^{F_\alpha}$. $\qquad \square$

For any $\lambda \in \mathbf{C}$ we have $V(\lambda) \otimes V(-\lambda) = \Lambda^{n-1} T^* \partial X \otimes \mathbf{C}$. Integration induces a $G$-invariant pairing of sections of $V(\lambda)$ with sections of $V(-\lambda)$. Therefore we define $C^{-\omega}(\partial X, V(\lambda)) := C^\omega(\partial X, V(-\lambda))'$. Similarly we proceed with $V_B(\lambda)$ and distribution sections of $V(\lambda)$ and $V_B(\lambda)$, respectively.

*Definition* 2.8. For $\mathrm{Re}(\lambda) > \delta_\Gamma$ and $\sharp \in \{\omega, \infty\}$ define the extension map $\mathrm{ext}_\lambda : C^{-\sharp}(B, V_B(\lambda)) \to C^{-\sharp}(\partial X, V(\lambda))$ by $\mathrm{ext}_\lambda := \pi'_{*, -\lambda}$.

2.4. *Elementary properties of* $\mathrm{ext}_\lambda$.

LEMMA 2.9. *For* $\mathrm{Re}(\lambda) > \delta_\Gamma$ *and* $\sharp \in \{\omega, \infty\}$ *the map*

$$
\mathrm{ext}_\lambda \in \mathrm{Hom}(C^{-\sharp}(B, V_B(\lambda)), C^{-\sharp}(\partial X, V(\lambda)))
$$

*depends holomorphically on* $\lambda$.

*Proof.* This follows from Lemmas 2.7 and 2.3. $\qquad \square$

LEMMA 2.10. *Let* $\mathrm{Re}(\lambda) > \delta_\Gamma$ *and* $\sharp \in \{\omega, \infty\}$. *Then the extension map* $\mathrm{ext}_\lambda$ *has values in* $^\Gamma C^{-\sharp}(\partial X, V(\lambda))$.

*Proof.* This follows immediately from the equation $\pi_{*,-\lambda}(\pi^{-\lambda}(g)f) = \pi_{*,-\lambda}(f)$, for all $f \in C^{\sharp}(\partial X, V(-\lambda))$. $\qquad\square$

Next we want to show that $\mathrm{ext}_\lambda$ is right inverse to $\mathrm{res}_\lambda$. For distributions this was shown in [8]. For hyperfunctions we want to argue by continuity. Therefore we must show that $\mathrm{res}_\lambda$ is continuous. This is not straightforward since we have defined $\mathrm{res}_\lambda$ as the composition of $\mathrm{res}_\Omega$ and an identification, and there is no natural topology on $C^{-\omega}(\Omega, V(\lambda))$.

LEMMA 2.11. *Let* $\lambda \in \mathbf{C}$. *Then the restriction map*
$$\mathrm{res}_\lambda : {}^\Gamma C^{-\omega}(\partial X, V(\lambda)) \to C^{-\omega}(B, V_B(\lambda))$$
*is continuous.*

*Proof.* We begin by giving a description of $\mathrm{res}_\lambda$ which is more appropriate for the present purpose. It is based on an alternative definition of the topological vector space $C^{-\omega}(\partial X, V(\lambda))$.

Let $\{F_\alpha\}$ be a finite cover of $B$ by compact neighbourhoods which have diffeomorphic lifts $\tilde{F}_\alpha \subset \Omega$. We consider the finite set $L := \{g \in \Gamma \mid \exists \alpha, \beta \text{ such that } g\tilde{F}_\alpha \cap \tilde{F}_\beta \neq \emptyset\}$. Set $W := \mathrm{clo}(\partial X \setminus \bigcup_{g \in L, \alpha} g\tilde{F}_\alpha)$. Then any hyperfunction $f \in C^{-\omega}(\partial X, V(\lambda))$ can be represented as a sum of hyperfunctions $f_{g,\alpha}$, $f_W \in C^{-\omega}(\partial X, V(\lambda))$, $g \in L$, where $\mathrm{supp}(f_{g,\alpha}) \subset g\tilde{F}_\alpha$, and $\mathrm{supp}(f_W) \subset W$.

If $f$ is $\Gamma$-invariant, then by lifting a corresponding image of $\mathrm{res}_\lambda(f)$ on $B$ to $\Omega$ we can choose the hyperfunctions $f_{g,\alpha}$ such that

(10) $$f_{g,\alpha} = \pi_\lambda(g)f_{1,\alpha}, \quad \text{for all} \quad \alpha, g \in L .$$

This choice made, we have $\mathrm{res}_\lambda(f) = \sum_\alpha f_{1,\alpha}$, where we view $f_{1,\alpha}$ as an element of $C^{-\omega}(B, V_B(\lambda))$ with support in $F_\alpha$.

We now argue that $\mathrm{res}_\lambda$ is continuous. Consider the Fréchet spaces
$$\mathcal{G} \;:=\; \bigoplus_{g \in L, \alpha} C^{-\omega}(g\tilde{F}_\alpha, V(\lambda)) \oplus C^{-\omega}(W, V(\lambda)) ,$$
$$^L\mathcal{G} \;:=\; \{f \in \mathcal{G} \mid f \text{ satisfies (10)}\} .$$

Here for any closed set $A \subset \partial X$ the space $C^{-\omega}(A, V(\lambda))$ is the space of hyperfunction sections of $V(\lambda)$ with support on $A$, i.e. the topological dual of the space of germs at $A$ of real analytic sections of $V(-\lambda)$. By the above we have a surjective continuous map $\phi : \mathcal{G} \to C^{-\omega}(\partial X, V(\lambda))$ with the property that $^\Gamma C^{-\omega}(\partial X, V(\lambda)) \subset \phi(^L\mathcal{G})$.

There is a continuous map $\tilde{\mathrm{res}} : \mathcal{G} \to C^{-\omega}(B, V_B(\lambda))$ given by $\tilde{\mathrm{res}}(f) := \sum_\alpha f_{1,\alpha}$. For $f \in {}^L\mathcal{G}$ we have $\tilde{\mathrm{res}}(f)_{|\mathrm{int}(F_\alpha)} = \phi(f)_{|\mathrm{int}(\tilde{F}_\alpha)}$. Thus $\tilde{\mathrm{res}}_{|\ker \phi \cap {}^L\mathcal{G}} = 0$, and $\tilde{\mathrm{res}}$ factorizes over the Fréchet space $\mathcal{F}$ defined by
$$\mathcal{F} := [\phi^{-1}(^\Gamma C^{-\omega}(\partial X, V(\lambda))) \cap {}^L\mathcal{G}]/[\ker(\phi) \cap {}^L\mathcal{G}] .$$
$\phi$ induces an isomorphism $\bar\phi : \mathcal{F} \to {}^\Gamma C^{-\omega}(\partial X, V(\lambda))$, and $\mathrm{res}_\lambda = \tilde{\mathrm{res}} \circ \bar\phi^{-1}$. Hence $\mathrm{res}_\lambda$ is continuous. $\qquad\square$

LEMMA 2.12. *For* $\mathrm{Re}(\lambda) > \delta_\Gamma$ *and* $\sharp \in \{\omega, \infty\}$ *the map* $\mathrm{ext}_\lambda$ *satisfies*

$$\mathrm{res}_\lambda \circ \mathrm{ext}_\lambda = \mathrm{id} \ .$$

*Proof.* For $\sharp = \infty$ this is as shown in [8, Lemma 4.5]. Since $C^{-\infty}(B, V_B(\lambda))$ is dense in $C^{-\omega}(B, V_B(\lambda))$, and $\mathrm{res}_\lambda \circ \mathrm{ext}_\lambda$ is continuous, the assertion for $\sharp = \omega$ follows. $\square$

If $^\Gamma C^{-\omega}(\Lambda, V(\lambda)) = 0$, then $\mathrm{res}_\lambda$ is injective. If in addition $\mathrm{Re}(\lambda) > \delta_\Gamma$, then we have $\mathrm{ext}_\lambda \circ \mathrm{res}_\lambda = \mathrm{id}$. In order to apply this observation we need vanishing results for $^\Gamma C^{-\omega}(\Lambda, V(\lambda))$.

LEMMA 2.13. *If* $\mathrm{Re}(\lambda) > 0$ *and* $\mathrm{Im}(\lambda) \neq 0$, *then* $^\Gamma C^{-\omega}(\Lambda, V(\lambda)) = 0$.

*Remark.* We will obtain much stronger vanishing results later.

*Proof.* The proof is based on the fact that the Laplace-Beltrami operator $\Delta_Y$ is self-adjoint. If $0 \neq \phi \in {}^\Gamma C^{-\omega}(\Lambda, V(\lambda))$, then, using the Poisson transform $P_\lambda$, we will construct a nontrivial $L^2$-eigenfunction $P_\lambda \phi$ of $\Delta_Y$ to an eigenvalue $\mu \notin \mathbf{R}$, obtaining a contradiction.

We now explain the details. Let $g = \kappa(g)a(g)n(g) \in KAN$ be the Iwasawa decomposition of $g \in G$.

*Definition* 2.14. The Poisson transform.

$$P_\lambda : C^{-\omega}(\partial X, V(\lambda)) \to C^\infty(X)$$

is defined by

$$(P_\lambda \phi)(gK) := \int_K a(g^{-1}k)^{-(\lambda+\rho)} \phi(k) dk \ .$$

Here $\phi \in C^{-\omega}(\partial X, V(\lambda))$ is viewed as a "function" on $G$ with values in $V_\lambda \cong \mathbf{C}$, and the integral is a formal notation meaning that the analytic functional $\phi$ has to be applied to the analytic integral kernel.

The Poisson transform $P_\lambda$ is continuous, $G$-equivariant, and $(\Delta_X - \rho^2 + \lambda^2)P_\lambda \phi = 0$. It is injective whenever $\lambda \notin -\rho - \mathbf{N}_0$. In this case $P_\lambda$ provides a topological isomorphism between $C^{-\omega}(\partial X, V(\lambda))$ and

$$\ker\left((\Delta_X - \rho^2 + \lambda^2) : C^\infty(X) \to C^\infty(X)\right) \ .$$

For all these facts see [23], [18], or [45].

Let $V \subset \partial X$ and $U \subset X$ be such that $\mathrm{clo}(U) \cap V = \emptyset$, where we take the closure of $U$ in $X \cup \partial X$. It is not difficult to show (see Lemma 3.2, (ii) below) that for $\mathrm{Re}(\lambda) > 0$ and $\phi \in C^{-\omega}(\partial X, V(\lambda))$, $\mathrm{supp}\,\phi \subset V$, the restriction of the Poisson transform $P_\lambda \phi$ to $U$ belongs to $L^2(U)$. If $\phi \in {}^\Gamma C^{-\omega}(\Lambda, V(\lambda))$, then $P_\lambda \phi$ is $\Gamma$-invariant, and therefore it descends to an eigenfunction of $\Delta_Y$ in $L^2(Y)$.

Since $Y$ is complete, the operator $\Delta_Y$ is self-adjoint with domain $\{f \in L^2(Y) | \Delta_Y f \in L^2(Y)\}$. In particular $\Delta_Y$ cannot have nontrivial eigenvectors in $L^2(Y)$ to eigenvalues with nontrivial imaginary part. Since $\mathrm{Im}(\lambda^2) \neq 0$ we conclude that $P_\lambda \phi = 0$ and hence $\phi = 0$ by injectivity of the Poisson transform. $\qquad\square$

### 2.5 The scattering matrix.

In this subsection we study the scattering matrix $\hat{S}_\lambda$ defined by (5). Our investigation is based on a detailed knowledge on the Knapp-Stein intertwining operators $\hat{J}_\lambda$ ([26], [50, Ch.10]). In order to fix our normalization conventions let us first give a definition of $\hat{J}_\lambda$ for $\mathrm{Re}(\lambda) < 0$.

*Definition* 2.15. Consider $f \in C^\infty(\partial X, V(\lambda))$ as a right $P$-equivariant function on $G$ with values in $V_\lambda$. Then $\hat{J}_\lambda : C^\infty(\partial X, V(\lambda)) \to C^\infty(\partial X, V(-\lambda))$ is defined by the convergent integral

$$(\hat{J}_\lambda f)(g) := \int_{\bar{N}} f(gw\bar{n}) d\bar{n} \ .$$

For $\mathrm{Re}(\lambda) \geq 0$ the operator $\hat{J}_\lambda$ is defined by meromorphic continuation ([26], [50, Ch.10]). It extends continuously to $C^{-\sharp}(\partial X, V(\lambda))$, $\sharp \in \{\infty, \omega\}$, and this extension coincides with the adjoint of $\hat{J}_\lambda$ acting on $C^\sharp(\partial X, V(\lambda))$.

To be precise at this point we must show that $\hat{J}_\lambda$ is meromorphic in the sense of subsection 2.2.

LEMMA 2.16. *For* $\mathrm{Re}(\lambda) < 0$ *and* $\sharp \in \{\pm\omega, \pm\infty\}$ *the intertwining operators form a holomorphic family of continuous maps* $\hat{J}_\lambda \in \mathrm{Hom}(C^\sharp(\partial X, V(\lambda)),$ $C^\sharp(\partial X, V(-\lambda)))$. *This family admits a meromorphic continuation to all of* $\mathbf{C}$ *having at most first order poles.*

*Proof.* For $\sharp = \infty$ this is as shown in [8]. The case $\sharp = -\infty$ follows by duality from Lemma 2.3.

We now show the lemma for $\sharp = \omega$. It is sufficient to prove the first assertion. The proof of the meromorphic continuation [8] (it is based on equation (17) ) for $\sharp = \infty$ applies equally well to $\sharp = \omega$.

Let $X_i$, $i = 1, \ldots, \dim(\mathbf{k})$, be an orthonormal base of $\mathbf{k}$. For any multi-index $r = (i_1, \ldots, i_{\dim(\mathbf{k})})$ we set

$$X_r = \prod_{l=1}^{\dim(\mathbf{k})} X_l^{i_l}, \ \ |r| = \sum_{l=1}^{\dim(\mathbf{k})} i_l, \ \ r! := \prod_{l=1}^{\dim(\mathbf{k})} i_l! \ ,$$

and for $f \in C^\omega(K)$ we define the seminorm

$$\|f\|_r = \sup_{k \in K} |f(X_r k)| \ .$$

For any $f \in C^\omega(K)$ there exists an $R > 0$ such that

$$\|f\|_R := \sum_r \frac{R^{|r|}}{r!} \|f\|_r < \infty \ . \tag{11}$$

Consider the Banach space $\mathcal{H}_R(K) := \{f \in C^\omega(K) | \ \|f\|_R < \infty\}$. If $0 < R' < R$, then we have an inclusion $\mathcal{H}_R(K) \hookrightarrow \mathcal{H}_{R'}(K)$. As a topological vector space,

$$C^\omega(K) = \varinjlim \mathcal{H}_R(K) \ .$$

Let $\lambda \in \mathbf{C}$ and $\text{Re}(\lambda) < 0$. As usual, we identify $C^\omega(\partial X, V(\lambda))$ with $C^\omega(K)^M$. We extend $f \in C^\omega(K)$ to a function $f_\lambda$ on $G$ by

$$f_\lambda(kan) = f(k)a^{\lambda - \rho} \ . \tag{12}$$

If $f \in C^\omega(\partial X, V(\lambda))$, then

$$\hat{J}_\lambda(f)(k) := \int_{\bar{N}} f_\lambda(kw\bar{n})d\bar{n} \ .$$

For any $\lambda_0 \in \mathbf{a}_{\mathbf{C}}^*$ with $\text{Re}(\lambda_0) < 0$ and $\delta > 0$ we can find an $\varepsilon > 0$ such that for $|\lambda - \lambda_0| < \varepsilon$

$$\int_{\bar{N}} |a(\bar{n})^{\lambda_0 - \rho} - a(\bar{n})^{\lambda - \rho}| d\bar{n} < \delta$$

holds. We then have

$$
\begin{aligned}
\|\hat{J}_{\lambda_0} f - \hat{J}_\lambda f\|_r &= \sup_{k \in K} | \int_{\bar{N}} f_{\lambda_0}(X_r kw\bar{n}) - f_\lambda(X_r kw\bar{n})d\bar{n}| \\
&\leq \sup_{k \in K} \int_{\bar{N}} |f(X_r kw\kappa(\bar{n}))||a(\bar{n})^{\lambda_0 - \rho} - a(\bar{n})^{\lambda - \rho}| d\bar{n} \\
&\leq \|f\|_r \int_{\bar{N}} |a(\bar{n})^{\lambda_0 - \rho} - a(\bar{n})^{\lambda - \rho}| d\bar{n} \\
&\leq \delta \|f\|_r \ .
\end{aligned}
$$

If $f \in \mathcal{H}_R(K)$, then we conclude

$$\|\hat{J}_{\lambda_0} f - \hat{J}_\lambda f\|_R \leq \delta \|f\|_R \ .$$

This proves continuity of the family of intertwining operators for $\sharp = \omega$ and $\text{Re}(\lambda) < 0$. Holomorphy now follows from [50, Lemma 10.1.3], and Lemma 2.1. This proves the lemma for $\sharp = \omega$. In the case $\sharp = -\omega$ we argue by duality using Lemma 2.3. $\qquad\square$

There is a meromorphic function $P(\lambda)$ such that the following functional equation holds:

$$\hat{J}_\lambda \circ \hat{J}_{-\lambda} = \frac{\text{id}}{P(\lambda)} \ . \tag{13}$$

The function $P(\lambda)$ is called the Plancherel density.

Let $1_\lambda \in C^\omega(\partial X, V(\lambda))$ be the unique $K$-invariant section normalized such that $1_\lambda(1) = 1$ (when viewed as a function on $G$ with values in $V_\lambda \cong \mathbf{C}$).

*Definition* 2.17. The *meromorphic function* $c(\lambda)$ is defined by

$$\hat{J}_\lambda 1_\lambda = c(-\lambda) 1_{-\lambda} \ .$$

We define the *normalized intertwining operator* by

$$J_\lambda := c(-\lambda)^{-1} \hat{J}_\lambda \ .$$

By (13),

$$(14) \qquad \frac{1}{P(\lambda)} = c(\lambda) c(-\lambda) \ .$$

This implies the following meromorphic identity:

$$(15) \qquad J_\lambda \circ J_{-\lambda} = \mathrm{id} \ .$$

Now we turn to the investigation of the scattering matrix $\hat{S}_\lambda$. Unfortunately, we are not able to show directly its meromorphic continuation as an operator on real analytic or hyperfunction sections of $V_B$. In [8] we studied the scattering matrix $\hat{S}_\lambda$ as an operator acting on smooth or distribution sections of $V_B(\lambda)$. First we recall some of these results.

LEMMA 2.18. (i) *For* $\mathrm{Re}(\lambda) > \delta_\Gamma$ *the scattering matrix defined by*

$$\hat{S}_\lambda := \mathrm{res}_{-\lambda} \circ \hat{J}_\lambda \circ \mathrm{ext}_\lambda : C^{-\infty}(B, V_B(\lambda)) \to C^{-\infty}(B, V_B(-\lambda))$$

*forms a meromorphic family of continuous maps, as does the normalized scattering matrix* $S_\lambda := c(-\lambda)^{-1} \hat{S}_\lambda$.

(ii) *The families* $\hat{S}_\lambda$, $S_\lambda$ *admit meromorphic continuations to all of* $\mathbf{C}$. *The singularities of* $\hat{S}_\lambda$ *and* $S_\lambda$ *in the region* $\mathrm{Re}(\lambda) < 0$ *are at most finite-dimensional.*

(iii) *The following functional equation holds*: $S_\lambda \circ S_{-\lambda} = \mathrm{id}$.

(iv) *The adjoint*

$$\hat{S}'_\lambda : C^\infty(B, V_B(\lambda)) \to C^\infty(B, V_B(-\lambda))$$

*coincides with the restriction of* $\hat{S}_\lambda$ *to* $C^\infty(B, V_B(\lambda))$. *In particular, this restriction defines a meromorphic family of continuous maps*

$$\hat{S}_\lambda : C^\infty(B, V_B(\lambda)) \to C^\infty(B, V_B(-\lambda)) \ .$$

(v) *The extension map*

$$\mathrm{ext}_\lambda : C^{-\infty}(B, V_B(\lambda)) \to C^{-\infty}(\partial X, V(\lambda))$$

*has a meromorphic continuation to all of* $\mathbf{C}$ *with at most finite-dimensional singularities. Moreover*, $\mathrm{ext}_\lambda$ *has values in* ${}^\Gamma C^{-\infty}(\partial X, V(\lambda))$.

(vi) *The equations*

$$\begin{aligned}
\mathrm{res}_\lambda \circ \mathrm{ext}_\lambda &= \mathrm{id}, \\
\hat{S}_\lambda &= \mathrm{res}_{-\lambda} \circ \hat{J}_\lambda \circ \mathrm{ext}_\lambda, \\
\mathrm{ext}_\lambda &= J_{-\lambda} \circ \mathrm{ext}_{-\lambda} \circ S_\lambda
\end{aligned}$$

*extend as meromorphic identities.*

*Proof.* All these assertions are shown in [8]. The main point is the meromorphy of the scattering matrix. The remaining assertions are easy consequences.

A meromorphic continuation of the scattering matrix was already obtained in several previous papers, e.g. in [34] in case $\delta_\Gamma < 0$, and in [40], [30], [31], [38] in the general case. Since these papers use different conventions for continuity we can only deduce meromorphy of matrix coefficients of the scattering matrix. But we can employ Cauchy's integral formula in order to obtain the meromorphy of the scattering matrix in the sense of the present paper. □

Let $W \subset \partial X$ be a closed set with nonempty interior and consider the space

(16) $$\mathcal{F}_\lambda := \{ f \in C^{-\omega}(\partial X, V(\lambda)) | f_{|W} \in C^\omega(W, V(\lambda)) \} .$$

We equip $\mathcal{F}_\lambda$ with the weakest topology such that the maps $\mathcal{F}_\lambda \hookrightarrow C^{-\omega}(\partial X, V(\lambda))$ and $\mathcal{F}_\lambda \to C^\omega(W, V(\lambda))$ are continuous. Using multiplication by suitable powers of a $K$-invariant volume form we can identify $\mathcal{F}_\lambda$, $\lambda \in \mathbf{C}$, with the fixed space $\mathcal{F}_0$. We employ this identification in order to speak of holomorphic families of vectors $f_\lambda \in \mathcal{F}_\lambda$.

Let $F \subset \mathrm{int}(W)$ be closed. By $\mathrm{res}_F$ we denote the restriction of sections of $V(\lambda)$ to $F$.

LEMMA 2.19. *The composition* $\mathrm{res}_F \circ \hat{J}_\lambda : \mathcal{F}_\lambda \to C^\omega(F, V(-\lambda))$ *is well defined and depends meromorphically on* $\lambda$.

*Proof.* In order to show that the composition is well defined we must show that if $f \in \mathcal{F}_\lambda$, then $\hat{J}_\lambda(f)$ is real analytic on $F$.

As a first step we reduce the proof of the lemma to the case $\mathrm{Re}(\lambda) < 0$ using the translation principle. In order to write down the appropriate formulas we identify $C^{-\omega}(\partial X, V(\lambda)) \cong C^{-\omega}(K)^M$. The image of $\mathcal{F}_\lambda$ under this identification is independent of $\lambda$ and will be denoted by $\mathcal{F}$.

By [50, Thm. 10.1.5], there are nonvanishing polynomial maps $b : \mathbf{a}^*_{\mathbf{C}} \to \mathbf{C}$ and $D : \mathbf{a}^*_{\mathbf{C}} \to \mathcal{U}(\mathbf{g})^K$, such that

(17) $$b(\lambda)\hat{J}_\lambda = \hat{J}_{\lambda-4\rho} \circ \pi^{\lambda-4\rho}(D(\lambda)) .$$

The family of differential operators $\pi^{\lambda-4\rho}(D(\lambda)) : \mathcal{F} \to \mathcal{F}$ is a holomorphic family of continuous maps. Thus a proof of the lemma for $\mathrm{Re}(\lambda) < \mu$ also implies a proof for $\mathrm{Re}(\lambda) < \mu + 4\rho$. It is therefore sufficient to study $\hat{J}_\lambda$ in its domain of convergence $\mathrm{Re}(\lambda) < 0$.

We are going to decompose $\mathrm{res}_F \circ \hat{J}_\lambda$ into a diagonal part $\hat{J}_\lambda^1$ and an off-diagonal part $\hat{J}_\lambda^2$. Let $F_1$ be closed such that $F \subset \mathrm{int}(F_1) \subset F_1 \subset \mathrm{int}(W)$. Then there is a compact set $V \subset \bar{N}$ such that $Fw\kappa(\bar{N} \setminus V)M \subset F_1 M$. Let $U \subset \mathrm{int}(V)$ be such that $Fw\kappa(\bar{N} \setminus U)M \subset WM$, and choose a cut-off function $\chi \in C^\infty(\bar{N})$ such that $\chi_{|U} = 1$ and $\chi_{|\bar{N} \setminus V} = 0$. If $f$ belongs to the dense subspace $C^\omega(\partial X) \subset \mathcal{F}$, then we set

$$J_\lambda^1(f)(k) := \int_{\bar{N}} f_\lambda(kw\bar{n})(1 - \chi(\bar{n}))d\bar{n} \ ,$$

$$J_\lambda^2(f)(k) := \int_{\bar{N}} f_\lambda(kw\bar{n})\chi(\bar{n})d\bar{n}, \quad k \in FM \ ,$$

where we employ the notation $f_\lambda$ introduced in (12). We have to show that these operators extend to continuous operators from $\mathcal{F}$ to $C^\omega(FM)^M$, and that this extension depends meromorphically on $\lambda$.

For any multi-index $r$,

$$
\begin{aligned}
\|J_\lambda^1(f)\|_{r,F} &\stackrel{\mathrm{def}}{=} \sup_{k \in FM} |\int_{\bar{N}} f_\lambda(X_r kw\bar{n})(1 - \chi(\bar{n}))d\bar{n}| \\
&\leq \sup_{k \in WM} |f(X_r k)| \int_{\bar{N}} a(\bar{n})^{\mathrm{Re}\lambda-\rho}d\bar{n} \\
&\leq C\|f\|_{r,W} \ ,
\end{aligned}
$$

where $C$ does not depend on $r$ and $f$. Recall the definition of the $R$-norm (11). If $R > 0$ is sufficiently small (depending on $f$), then $\|J_\lambda^1(f)\|_{R,F} \leq C\|f\|_{R,W}$. Thus for fixed $\lambda$, $\mathrm{Re}(\lambda) < 0$, there is a continuous extension $J_\lambda^1 : \mathcal{F} \to C^\omega(FM)^M$. An argument similar to the one given in the proof of Lemma 2.16 shows that $J_\lambda^1$ depends holomorphically on $\lambda$.

Let $\psi : K/M \setminus M \to \bar{N}$ be the inverse of the analytic diffeomorphism $\bar{n} \mapsto w\kappa(\bar{n})M$. If $f \in C^\omega(\partial X)$ and $x \in FM$, then we can write

$$
\begin{aligned}
J_\lambda^2(f)(x) &= \int_{\bar{N}} f_\lambda(xw\bar{n})\chi(\bar{n})d\bar{n} \\
&= \int_{\bar{N}} f(xw\kappa(\bar{n}))a(\bar{n})^{\lambda-\rho}\chi(\bar{n})d\bar{n} \\
&= \int_K f(xk)a(\psi(k))^{\lambda+\rho}\chi(\psi(k))dk \ .
\end{aligned}
$$

Let $\mathcal{G}$ be defined in the same way as $\mathcal{F}$ (see (16)), but with $W$ replaced by $w\kappa(\bar{N} \setminus \mathrm{int}(U)) \supset \mathrm{supp}(1 - \chi \circ \psi)$. There is a continuous injection of $\mathcal{F}$ into the space of analytic functions from $FM$ to $\mathcal{G}$ given by $f \mapsto a_f$, $a_f(x) := f(x) \in \mathcal{G}$, $x \in FM$. The multiplication by $\chi \circ \psi$ is a continuous operator $m_\chi : \mathcal{G} \to C^{-\omega}(\mathrm{supp}(\chi \circ \psi))$. Since the restriction of $a(\psi(.))^{\lambda+\rho}$ to $\mathrm{supp}(\chi \circ \psi)$ is analytic, it defines a continuous functional on $C^{-\omega}(\mathrm{supp}(\chi \circ \psi))$. Now we obtain a continuous extension of $J_\lambda^2$ to $\mathcal{F}$ by

$$(18) \qquad J_\lambda^2(f)(x) := \langle a(\psi(.))^{\lambda+\rho}, m_\chi(a_f(x)) \rangle \ .$$

It is now easy to verify that the map $\mathbf{a}_\mathbf{C}^* \ni \lambda \mapsto J_\lambda^2 \in \mathrm{Hom}(\mathcal{F}, C^\omega(F))$ is continuous and that for fixed $f$ and $x$ the function $\lambda \mapsto J_\lambda^2(f)(x)$ is holomorphic. Thus by Lemma 2.1 the family of maps $J_\lambda^2$ is holomorphic. □

LEMMA 2.20. (i) *The restrictions of $\hat{S}_\lambda$, $S_\lambda$ to real analytic sections of $V_B(\lambda)$ form meromorphic families of continuous maps*

$$\hat{S}_\lambda, \ S_\lambda \in \mathrm{Hom}(C^\omega(B, V_B(\lambda)), C^\omega(B, V_B(-\lambda))) \ .$$

(ii) *The singularities of $\hat{S}_\lambda$ and $S_\lambda$ in the region $\mathrm{Re}(\lambda) < 0$ are at most finite-dimensional.*

(iii) *The adjoint*

$$\hat{S}_\lambda' : C^{-\omega}(B, V_B(\lambda)) \to C^{-\omega}(B, V_B(-\lambda))$$

*of*

$$\hat{S}_\lambda : C^\omega(B, V_B(\lambda)) \to C^\omega(B, V_B(-\omega))$$

*provides a continuous extension of $\hat{S}_\lambda$ to $C^{-\omega}(B, V_B(\lambda))$. The same holds true for $S_\lambda$.*

(iv) *The following functional equation holds on $C^{-\omega}(B, V_B(-\lambda))$ : $S_\lambda \circ S_{-\lambda}$ = id.*

*Proof.* (i) and (ii): Let $p : \Omega \to B$ be the projection. We choose compact subsets $F \subset \mathrm{int}(W) \subset W \subset \Omega$ such that $p(F) = B$ and define $\mathcal{F}_\lambda$ as above. It follows from Lemma 2.18 (v), that $\mathrm{ext}_\lambda : C^\omega(B, V_B(\lambda)) \to C^{-\omega}(\partial X, V(\lambda))$ is a meromorphic family of continuous maps with at most finite-dimensional singularities. Since $\mathrm{res}_\lambda \circ \mathrm{ext}_\lambda = \mathrm{id}$ we have a holomorphic family of continuous maps $\mathrm{res}_W \circ \mathrm{ext}_\lambda : C^\omega(B, V_B(\lambda)) \to C^\omega(W, V(\lambda))$. Thus $\mathrm{ext}_\lambda : C^\omega(B, V_B(\lambda)) \to \mathcal{F}_\lambda$ is meromorphic with at most finite-dimensional singularities. We can identify $C^\omega(B, V_B(\lambda))$ with a closed subspace of $C^\omega(F, V(\lambda))$. With this identification, $\mathrm{res}_F = \mathrm{res}_\lambda$. Assertions (i) and (ii) now follow from Lemma 2.19.

(iii) and (iv): The adjoint $\hat{S}_\lambda'$ is the continuous extension of $\hat{S}_\lambda$ defined on distribution sections. The assertions now follow from Lemma 2.18 (iv) and (iii). □

2.6. *Extension of hyperfunctions and the embedding trick.*

In this subsection we want to construct the meromorphic continuation of the extension of hyperfunction sections of $V_B(\lambda)$.

Lemma 2.21. *If $\delta_\Gamma < 0$, then*

$$\mathrm{ext}_\lambda \in \mathrm{Hom}(C^{-\omega}(B, V_B(\lambda)), C^{-\omega}(\partial X, V(\lambda)))$$

*(initially defined for $\mathrm{Re}(\lambda) > \delta_\Gamma$) admits a meromorphic continuation to all of $\mathbf{C}$ with at most finite-dimensional singularities. Moreover, $\mathrm{ext}_\lambda$ has values in ${}^\Gamma C^{-\omega}(\partial X, V(\lambda))$.*

*Proof.* If $-\mathrm{Re}(\lambda) > \delta_\Gamma$, then using Lemma 2.20 (i), (iii) we define $\widetilde{\mathrm{ext}}_\lambda := J_{-\lambda} \circ \mathrm{ext}_{-\lambda} \circ S_\lambda$. We have by Lemma 2.20 (iv):

$$\mathrm{res}_\lambda \circ \widetilde{\mathrm{ext}}_\lambda = \mathrm{res}_\lambda \circ J_{-\lambda} \circ \mathrm{ext}_{-\lambda} \circ S_\lambda = S_{-\lambda} \circ S_\lambda = \mathrm{id} \ .$$

Since for $0 < \mathrm{Re}(\lambda) < -\delta_\Gamma$ and $\mathrm{Im}(\lambda) \neq 0$ the restriction map $\mathrm{res}_\lambda$ is injective by Lemma 2.13 and $\mathrm{res}_\lambda \circ \mathrm{ext}_\lambda = \mathrm{res}_\lambda \circ \widetilde{\mathrm{ext}}_\lambda = \mathrm{id}$ we conclude that $\mathrm{ext}_\lambda = \widetilde{\mathrm{ext}}_\lambda$. Thus $\widetilde{\mathrm{ext}}_\lambda$ is a meromorphic continuation of $\mathrm{ext}_\lambda$ (defined on hyperfunction sections).

Note that the family $\mathrm{ext}_\lambda$ defined on hyperfunction sections as constructed above is just the continuous extension of the previously obtained family $\mathrm{ext}_\lambda$ defined on distribution sections. Since $\mathrm{ext}_\lambda$ considered as an operator on distribution sections has at most finite-dimensional singularities by Lemma 2.18 (v), its continuous extension to hyperfunction sections has the same singularities which are, in particular, finite-dimensional. $\qquad \square$

We now show how to drop the assumption $\delta_\Gamma < 0$ using the embedding trick.

Proposition 2.22. *Lemma 2.21 holds true without the assumption $\delta_\Gamma < 0$.*

*Proof.* We realize $\mathrm{SO}(1, n)_0$ as the group of automorphisms of $\mathbf{R}^{n+1}$ fixing the quadratic form given by $\mathrm{diag}(-1, \underbrace{1, \ldots, 1}_{n\times})$. We view $\mathrm{SO}(1, n)_0$ as the subgroup of $\mathrm{SO}(1, n+1)_0$ which fixes the last element of the standard base of $\mathbf{R}^{n+2}$. This embedding is compatible with the standard Cartan involution $g \mapsto {}^t g^{-1}$.

We set $G^n := \mathrm{SO}(1, n)_0$. Then we have a sequence $G^n \hookrightarrow G^{n+1}$, $n \geq 2$, of embeddings inducing embeddings of the corresponding Iwasawa constituents $K^n \hookrightarrow K^{n+1}$, $N^n \hookrightarrow N^{n+1}$, $M^n \hookrightarrow M^{n+1}$, and $A^n \overset{\cong}{\hookrightarrow} A^{n+1}$. In particular we identify the Lie algebras $\mathbf{a}^n$ of $A^n$ in a compatible way with $\mathbf{R}$.

The embedding $G^n \hookrightarrow G^{n+1}$ induces a totally geodesic embedding $X^n \hookrightarrow X^{n+1}$ and an embedding of the geodesic boundaries $\partial X^n \hookrightarrow \partial X^{n+1}$. All these embeddings are equivariant with respect to the action of $G^n$.

If we view the discrete subgroup $\Gamma \subset G^n$ as a subgroup of $G^{n+1}$, then it is still convex cocompact. By $\Omega^n$, $\Omega^{n+1}$, $\Lambda^n$, $\Lambda^{n+1}$ we denote the corresponding domains of discontinuity and limit sets. Under the embedding $\partial X^n \hookrightarrow \partial X^{n+1}$ the limit set $\Lambda^n$ is identified with $\Lambda^{n+1}$. Moreover, we have an embedding $\Omega^n \hookrightarrow \Omega^{n+1}$ inducing the embedding of compact quotients $B^n \hookrightarrow B^{n+1}$.

The exponent of $\Gamma$ now depends on $n$ and is denoted by $\delta_\Gamma^n$. We have the relation $\delta_\Gamma^{n+1} = \delta_\Gamma^n - \frac{1}{2}$. Thus for $m$ sufficiently large, $\delta_\Gamma^{n+m} < 0$.

Let $\text{ext}_\lambda^n$ denote the extension map associated to $\Gamma \subset G^n$. The aim of the following discussion is to show how the meromorphic continuation $\text{ext}_\lambda^{n+1}$ leads to the continuation of $\text{ext}_\lambda^n$.

Let $P^n := M^n A^n N^n$, $V(\lambda)^n := G^n \times_{P^n} V_\lambda$, and $V_{B^n}(\lambda) = \Gamma \backslash V(\lambda)^n_{|\Omega^n}$. The representation $V_\lambda$ of $P^{n+1}$ restricts to the representation $V_{\lambda - \frac{1}{2}}$ of $P^n$. This induces isomorphisms of bundles

$$V(\lambda)^{n+1}_{|\partial X^n} \cong V(\lambda - \frac{1}{2})^n, \quad V_{B^{n+1}}(\lambda)_{|B^n} \cong V_{B^n}(\lambda - \frac{1}{2}) \ .$$

Let

$$i^* : C^\omega(B^{n+1}, V_{B^{n+1}}(\lambda)) \to C^\omega(B^n, V_{B^n}(\lambda - \frac{1}{2}))$$

and

$$j^* : C^\omega(\partial X^{n+1}, V(\lambda)^{n+1}) \to C^\omega(\partial X^n, V(\lambda - \frac{1}{2})^n)$$

denote the maps given by restriction of sections. Note that $j^*$ is $G^n$-equivariant. The adjoint maps define the push forward of hyperfunction sections

$$i_* : C^{-\omega}(B^n, V_{B^n}(\lambda)) \quad \to \quad C^{-\omega}(B^{n+1}, V_{B^{n+1}}(\lambda - \frac{1}{2})) \ ,$$

$$j_* : C^{-\omega}(\partial X^n, V(\lambda)^n) \quad \to \quad C^{-\omega}(\partial X^{n+1}, V(\lambda - \frac{1}{2})^{n+1}) \ .$$

If $\phi \in C^{-\omega}(B^n, V_{B^n}(\lambda))$, then the push forward $i_*\phi$ has support in $B^n \subset B^{n+1}$. Since $\text{res}_{\lambda - \frac{1}{2}}^{n+1} \circ \text{ext}_{\lambda - \frac{1}{2}}^{n+1} = \text{id}$,

$$(19) \qquad \text{supp}(\text{ext}_{\lambda - \frac{1}{2}}^{n+1} \circ i_*)(\phi) \subset \Lambda^{n+1} \cup \Omega^n = \partial X^n \ .$$

We are now going to continue $\text{ext}_\lambda^n$ using $i_*$, $\text{ext}_{\lambda - \frac{1}{2}}^{n+1}$ and a left inverse of $j_*$. As on previous occasions we trivialize $V(\lambda)^{n+1}$ by powers of a $K^{n+1}$-invariant volume form. We thus identify $C^{-\omega}(\partial X^{n+1}, V(\lambda)^{n+1})$ with $C^{-\omega}(\partial X^{n+1})$ for all $\lambda \in \mathbf{C}$. Let $U \subset \partial X^{n+1}$ be a closed tubular neighbourhood of $\partial X$ and fix an analytic diffeomorphism $T : [-1, 1] \times \partial X^n \overset{\cong}{\to} U$. Then we can define a continuous extension $t : C^\omega(\partial X^n) \to C^\omega(U)$ by $(T^* \circ t)f(r, x) := f(x)$, $r \in [-1, 1]$, $x \in \partial X^n$. Let $t' : C^{-\omega}(U) \to C^{-\omega}(\partial X^n)$ be the adjoint of $t$. Then $t' \circ j_* = \text{id}$. Because of (19) we can define

$$\widetilde{\text{ext}}_\lambda^n \phi := (t' \circ \text{ext}_{\lambda - \frac{1}{2}}^{n+1} \circ i_*)(\phi) \ .$$

Then

$$\widetilde{\text{ext}}_\lambda^n \in \text{Hom}(C^{-\omega}(B^n, V_{B^n}(\lambda)), C^{-\omega}(\partial X^n, V(\lambda)^n))$$

is a meromorphic family of continuous maps with at most finite-dimensional singularities.

In order to prove that $\widetilde{\mathrm{ext}}_\lambda^n$ provides the desired meromorphic continuation it remains to show that it coincides with $\mathrm{ext}_\lambda^n$ in the region $\mathrm{Re}(\lambda) > \delta_\Gamma^n$. If $\mathrm{Re}(\lambda) > \delta_\Gamma^n$, then $\mathrm{Re}(\lambda) - \frac{1}{2} > \delta_\Gamma^{n+1}$, and the push-down maps $\pi_{*,-\lambda}^n$, $\pi_{*,-\lambda+\frac{1}{2}}^{n+1}$ are defined. It is easy to see from the definition of the push-down that

$$i^* \circ \pi_{*,-\lambda+\frac{1}{2}}^{n+1} = \pi_{*,-\lambda}^n \circ j^* \ .$$

Taking adjoints we obtain $\mathrm{ext}_{\lambda-\frac{1}{2}}^{n+1} \circ i_* = j_* \circ \mathrm{ext}_\lambda^n$. Therefore,

$$\widetilde{\mathrm{ext}}_\lambda^n = t' \circ \mathrm{ext}_{\lambda-\frac{1}{2}}^{n+1} \circ i_* = t' \circ j_* \circ \mathrm{ext}_\lambda^n = \mathrm{ext}_\lambda^n \ .$$

It follows by meromorphy that $\mathrm{im}(\widetilde{\mathrm{ext}}_\lambda^n)$ consists of $\Gamma$-invariant sections for all $\lambda \in \mathbf{C}$. This finishes the proof of the proposition. $\qquad\square$

LEMMA 2.23. *The identities* 2.18 (vi), *are valid on hyperfunction sections.*

*Proof.* This follows by continuity. $\qquad\square$

## 3. Green's formula and applications

### 3.1. *Asymptotics of Poisson transforms.*

In this subsection we recall facts about the asymptotic behaviour of the Poisson transform $P_\lambda \phi \in C^\infty(X)$ near $\partial X$, where $\phi \in C^{-\omega}(\partial X, V(\lambda))$.

LEMMA 3.1. *Let* $\partial X = U \cup Q$, *where* $U$ *is open and* $Q = \partial X \setminus U$. *Then for any closed subset* $F \subset U$ *the composition* $\mathrm{res}_F \circ \hat{J}_\lambda : C^{-\omega}(Q, V(\lambda)) \to C^\omega(F, V(-\lambda))$ *extends to a holomorphic family of continuous maps.*

*Proof.* Choose a closed set $W$ with nonempty interior such that $F \subset W \subset U$. Then $C^{-\omega}(Q, V(\lambda)) \subset \mathcal{F}_\lambda$ (see (16)), and we see by Lemma 2.19 that $\mathrm{res}_F \circ \hat{J}_\lambda : C^{-\omega}(Q, V(\lambda)) \to C^\omega(F, V(-\lambda))$ is a meromorphic family of continuous maps. We argue that it is in fact holomorphic. In the course of the proof of Lemma 2.19 we have constructed a decomposition $\mathrm{res}_F \circ \hat{J}_\lambda = \hat{J}_\lambda^1 + \hat{J}_\lambda^2$, valid for $\mathrm{Re}(\lambda) < 0$, where $\hat{J}_\lambda^2$ admits a holomorphic continuation to all of $\mathbf{C}$ by (18). Since $\hat{J}_{\lambda \mid C^{-\omega}(Q, V(\lambda))}^1 \equiv 0$, we conclude that $\mathrm{res}_F \circ \hat{J}_\lambda = \hat{J}_\lambda^2$ is holomorphic. $\qquad\square$

LEMMA 3.2.   (i) *Let* $f \in C^\infty(\partial X, V(\lambda))$. *If* $\mathrm{Re}(\lambda) > 0$, *then there exists* $\varepsilon > 0$ *such that for* $a \to \infty$, *we have*

$$(P_\lambda f)(ka) = a^{\lambda-\rho} c(\lambda) f(k) + O(a^{\lambda-\rho-\varepsilon})$$

*uniformly in* $k \in K$. *If* $\mathrm{Re}(\lambda) = 0$, $\lambda \neq 0$, *then*

$$(P_\lambda f)(ka) = a^{\lambda-\rho} c(\lambda) f(k) + a^{-\lambda-\rho} \hat{J}_\lambda f(k) + O(a^{-\rho-\varepsilon})$$

*uniformly in* $k \in K$.

(ii) *Let $\partial X = U \cup Q$, where $U$ is open and $Q := \partial X \setminus U$. Let $f \in C^{-\omega}(Q, V(\lambda))$. Then there exist smooth functions $\psi_n$ on $U$ such that*

$$(20) \qquad (P_\lambda f)(ka) = a^{-(\lambda+\rho)}(\hat{J}_\lambda f)(k) + \sum_{n \geq 1} a^{-(\lambda+\rho)-n} \psi_n(k) , \quad k \in U .$$

*The series converges uniformly for $a \gg 0$ and $kM$ in compact subsets of $U$. In particular, there exists $\varepsilon > 0$ such that for $a \to \infty$,*

$$(P_\lambda f)(ka) = a^{-\lambda-\rho}(\hat{J}_\lambda f)(k) + O(a^{-\lambda-\rho-\varepsilon})$$

*uniformly as $kM$ varies in compact subsets of $U$.*

(iii) *Let $U, Q$ be as in* (ii) *and $f \in C^{-\omega}(\partial X, V(\lambda))$ such that $\mathrm{res}_U f \in C^\infty(U, V(\lambda))$. If $\mathrm{Re}(\lambda) > 0$, then there exists $\varepsilon > 0$ such that for $a \to \infty$,*

$$(P_\lambda f)(ka) = a^{\lambda-\rho} c(\lambda) f(k) + O(a^{\lambda-\rho-\varepsilon})$$

*uniformly as $kM$ varies in compact subsets of $U$. If $\mathrm{Re}(\lambda) = 0, \lambda \neq 0$, then there exists $\varepsilon > 0$ such that for $a \to \infty$,*

$$(P_\lambda f)(ka) = a^{\lambda-\rho} c(\lambda) f(k) + a^{-\lambda-\rho}(\hat{J}_\lambda f)(k) + O(a^{-\rho-\varepsilon})$$

*uniformly as $kM$ varies in compact subsets of $U$.*

*The asymptotic expansions can be differentiated with respect to $a$.*

*Proof.* Assertion (i) is a simple consequence of the general results concerning the asymptotics of eigenfunctions on symmetric spaces [2, Thms. 3.5 and 3.6], combined with the limit formulas for the Poisson transform (see [45, Thm. 5.1.4], [49, Thm. 5.3.4]):

$$\lim_{a \to \infty} a^{\rho-\lambda}(P_\lambda f)(ka) = c(\lambda) f(k), \qquad \mathrm{Re}(\lambda) > 0 ,$$
$$\lim_{a \to \infty} a^{\rho+\lambda}(P_\lambda f)(ka) = (\hat{J}_\lambda f)(k), \qquad \mathrm{Re}(\lambda) < 0 .$$

For $f \in C(\partial X, V(\lambda))$, $\mathrm{supp}(f) \subset Q$, assertion (ii) follows from Theorem 4.8 in [3]. However, the proof given in that paper works as well for $f \in C^{-\omega}(Q, V(\lambda))$.

(iii) is a consequence of (i) and (ii). Indeed, let $W, W_1 \subset U$ be compact subsets such that $W \subset \mathrm{int}(W_1)$. Let $\chi \in C_c^\infty(U)$ be such that $\chi_{|W_1} \equiv 1$. Then we can write $f = \chi f + (1-\chi)f$, where $\chi f$ is smooth and $\mathrm{supp}(1-\chi)f \subset \partial X \setminus \mathrm{int}(W_1)$. We now apply (i) to $\chi f$ and (ii) to $(1-\chi)f$ for $kM \in W$. $\quad\square$

### 3.2. *An orthogonality result.*

The following facts were shown in [8, §6].

LEMMA 3.3. *Let $F \subset X \cup \Omega$ be a fundamental domain for $\Gamma$. There exists a cut-off function $\chi \in C^\infty(X)$ such that:*

(i) *There is a finite subset $L \subset \Gamma$ such that $\operatorname{supp}(\chi) \subset \bigcup_{g \in L} gF$;*

(ii) $\sum_{g \in \Gamma} g^* \chi = 1$;

(iii) *For any $i \in \mathbf{N}$ and compact $W \subset \Omega$, $\sup_{k \in WM, a \in A_+} a \, |\nabla^i \chi(ka)| < \infty$;*

(iv) *$\chi$ extends continuously to $\partial X$ and the restriction $\chi_\infty$ of this extension to $\Omega$ is a smooth function with compact support.*

Let $\phi \in {}^\Gamma C^{-\omega}(\Lambda, V(\lambda))$. Then by Lemma 3.1 the composition $\operatorname{res}_{-\lambda} \circ \hat{J}_\lambda(\phi) \in {}^\Gamma C^\omega(\Omega, V(-\lambda)) = C^\omega(B, V_B(-\lambda))$ is well-defined for all $\lambda \in \mathbf{C}$.

PROPOSITION 3.4. *If $\operatorname{Re}(\lambda) \geq 0$, $\phi \in {}^\Gamma C^{-\omega}(\Lambda, V(\lambda))$, and $f \in {}^\Gamma C^{-\omega}(\partial X, V(\lambda))$, then*

$$\langle \operatorname{res}_{-\lambda} \circ \hat{J}_\lambda(\phi), \operatorname{res}_\lambda(f) \rangle = 0 \ .$$

*Proof.* At first we need the following:

LEMMA 3.5. *The space*

$${}^\Gamma C_\Omega^{-\omega}(\partial X, V(\lambda)) := \{ f \in {}^\Gamma C^{-\omega}(\partial X, V(\lambda)) \mid f_{|\Omega} \in C^\infty(\Omega, V(\lambda)) \}$$

*is dense in ${}^\Gamma C^{-\omega}(\partial X, V(\lambda))$.*

*Proof.* By Proposition 2.22 the family $\operatorname{ext}_\mu$ has an at most finite-dimensional singularity at $\mu = \lambda$. Thus there is a finite-dimensional subspace $W \subset C^\omega(B, V_B(-\lambda))$ such that $(\operatorname{ext}_\lambda)_{|W^\perp} : W^\perp \to {}^\Gamma C^{-\omega}(\partial X, V(\lambda))$ is a well-defined continuous map, where $W^\perp := \{ \phi \in C^\omega(B, V_B(\lambda)) \mid \langle \phi, W \rangle = \{0\} \}$. Since $C^\infty(B, V_B(\lambda)) \subset C^{-\omega}(B, V_B(\lambda))$ is dense we can choose a complement $\tilde{W} \subset C^\infty(B, V_B(\lambda))$ such that $C^{-\omega}(B, V_B(\lambda)) = W^\perp \oplus \tilde{W}$.

Let $f \in {}^\Gamma C^{-\omega}(\partial X, V(\lambda))$. Then we can write $\operatorname{res}_\lambda f = g = g^\perp \oplus \tilde{g}$, $g^\perp \in W^\perp$, $\tilde{g} \in \tilde{W}$. Now $\operatorname{res}_\lambda(f - \operatorname{ext}_\lambda(g^\perp)) = \tilde{g}$. It follows that $f - \operatorname{ext}_\lambda(g^\perp) \in {}^\Gamma C_\Omega^{-\omega}(\partial X, V(\lambda))$. Let now $g_i \in C^\infty(B, V_B(\lambda))$ be a sequence such that $\lim_{i \to \infty} g_i = g$. Then we can decompose $g_i = g_i^\perp + \tilde{g}_i$. The sections $g_i^\perp$ are smooth since $g_i$ and $\tilde{g}_i \in \tilde{W}$ are so. It follows that $\operatorname{ext}_\lambda(g_i^\perp) \in {}^\Gamma C_\Omega^{-\omega}(\partial X, V(\lambda))$. By continuity of $(\operatorname{ext}_\lambda)_{|W^\perp}$ we have $\operatorname{ext}_\lambda(g^\perp) = \lim_{i \to \infty} \operatorname{ext}_\lambda(g_i^\perp)$. The assertion of the lemma now follows from $f = f - \operatorname{ext}_\lambda(g^\perp) + \lim_{i \to \infty} \operatorname{ext}_\lambda(g_i^\perp)$. $\square$

By Lemma 3.5 it suffices to prove the proposition for $f \in {}^\Gamma C_\Omega^{-\omega}(\partial X, V(\lambda))$. By Lemma 2.13 we can assume that $\lambda \in [0, \infty) \cup \iota\mathbf{R}$.

Let $A := \Delta_X - \rho^2 + \lambda^2$, $P := P_\lambda$ be the Poisson transform, and let $\chi$ be the cut-off function given by Lemma 3.3. By $B_R \subset X$ we denote the metric ball of radius $R$ centered at the origin represented by $K$.

Note that $APf = 0$ and $AP\phi = 0$ (see subsection 2.4). Integration by parts gives

$$
\begin{aligned}
(21) \quad 0 &= \langle \chi AP\phi, Pf \rangle_{L^2(B_R)} - \langle \chi P\phi, APf \rangle_{L^2(B_R)} \\
&= \langle A\chi P\phi, Pf \rangle_{L^2(B_R)} - \langle \chi P\phi, APf \rangle_{L^2(B_R)} - \langle [A, \chi] P\phi, Pf \rangle_{L^2(B_R)} \\
&= -\langle \nabla_n \chi P\phi, Pf \rangle_{L^2(\partial B_R)} + \langle \chi P\phi, \nabla_n Pf \rangle_{L^2(\partial B_R)} \\
&\quad - \langle [A, \chi] P\phi, Pf \rangle_{L^2(B_R)} ,
\end{aligned}
$$

where $n$ is the exterior unit normal vector field at $\partial B_R$. For the following discussion we distinguish between the three cases:

 (i) $\lambda \in (0, \infty)$,

 (ii) $\mathrm{Re}(\lambda) = 0$, $\lambda \neq 0$, and

 (iii) $\lambda = 0$.

We first consider the case $\lambda > 0$. Lemma 3.3 (iii) implies that $|[A, \chi] P\phi Pf|$ is integrable over all of $X$, and by Lemma 3.3 (ii), Lemma 3.2 (ii), (iii) and the $\Gamma$-invariance of $f, \phi$,

$$
\langle [A, \chi] P\phi, Pf \rangle_{L^2(X)} = \langle \sum_{\gamma \in \Gamma} \gamma^*([A, \chi]) P\phi, Pf \rangle_{L^2(Y)} = 0 .
$$

Taking the limit $R \to \infty$ in (21), and using Lemma 3.2 (i), (ii) and Lemma 3.3, we obtain

$$
\begin{aligned}
0 &= (\lambda + \rho) \int_K \chi_\infty(k)(\hat{J}_\lambda \phi)(k) c(\lambda) f(k) dk \\
&\quad + (\lambda - \rho) \int_K \chi_\infty(k)(\hat{J}_\lambda \phi)(k) c(\lambda) f(k) dk \\
&= 2\lambda c(\lambda) \int_K \chi_\infty(k)(\hat{J}_\lambda \phi)(k) f(k) dk \\
&= 2\lambda c(\lambda) \langle \mathrm{res}_{-\lambda} \circ \hat{J}_\lambda(\phi), \mathrm{res}_\lambda(f) \rangle .
\end{aligned}
$$

This is the assertion of the proposition for $\lambda > 0$ since $c(\lambda) \neq 0$ (see [19, Ch. IV, Thm. 6.14]).

Now we discuss the case $\mathrm{Re}(\lambda) = 0$ and $\lambda \neq 0$. By Lemma 3.2 (i), the function $P(f)(ka)$ now has two leading terms. Instead of taking the limit $R \to \infty$ in (21) we apply $\lim_{r \to \infty} \frac{1}{r} \int_0^r dR$. Again we have

$$
\lim_{r \to \infty} \frac{1}{r} \int_0^r \langle [A, \chi] P\phi, Pf \rangle_{L^2(B_R)} dR = 0 .
$$

The leading term $a^{-\lambda - \rho} \hat{J}_\lambda f(k)$ of $P(f)(ka)$ does not contribute to the limit since

$$
-2\lambda \lim_{r \to \infty} \frac{1}{r} \int_0^r R^{-2\lambda - 2\rho} \langle \chi(R) \hat{J}_\lambda \phi, \hat{J}_\lambda f \rangle_{L^2(\partial X)} dR = 0 .
$$

The contribution of the term $a^{\lambda-\rho}c(\lambda)f(k)$ leads to

$$0 = 2\lambda c(\lambda)\langle\mathrm{res}_{-\lambda} \circ \hat{J}_\lambda(\phi), \mathrm{res}_\lambda(f)\rangle$$

as in the case $\lambda > 0$. The proposition again follows since $c(\lambda) \neq 0$.

Now we consider the last case $\lambda = 0$. The function $c(\lambda)$ has a first order pole at $\lambda = 0$ with residue $c_1 \neq 0$ (see [19, Ch. IV, Thm. 6.14]). We redo the computation (21) using the Poisson transform $P = P_\mu$ at $\mu$ in a neighbourhood of 0. Of course in general $P_\mu(f), P_\mu(\phi)$ are not $\Gamma$-invariant except for $\mu = 0$. Nevertheless,

$$\lim_{R\to\infty}\langle[A,\chi]P_{1/R}\phi, P_{1/R}f\rangle_{L^2(B_R)} = 0$$

by the theorem of Lebesgue about dominated convergence. Moreover,

$$
\begin{aligned}
- \lim_{R\to\infty}&(\langle\nabla_n\chi P_{1/R}\phi, P_{1/R}f\rangle_{L^2(\partial B_R)} - \langle\chi P_{1/R}\phi, \nabla_n P_{1/R}f\rangle_{L^2(\partial B_R)})\frac{e^{2\rho R}}{\mathrm{vol}(\partial B_R)}\\
&= \lim_{R\to\infty}((1/R+\rho)\int_K \chi_\infty(k)(\hat{J}_{1/R}\phi)(k)c(1/R)f(k)dk\\
&\qquad +(1/R-\rho)\int_K \chi_\infty(k)(\hat{J}_{1/R}\phi)(k)c(1/R)f(k)dk)\\
&= 2\int_K \chi_\infty(k)(\hat{J}_0\phi)(k)\lim_{R\to\infty}\frac{c(1/R)}{R}f(k)dk\\
&= 2c_1\langle\mathrm{res}_0 \circ \hat{J}_0(\phi), \mathrm{res}_0(f)\rangle \ .
\end{aligned}
$$

It follows from (21) that the latter pairing vanishes. This finishes the proof of the proposition in the last case $\lambda = 0$. □

### 3.3. *Miscellaneous results.*

The following lemma is very similar to the general result [3, Thm. 4.1].

LEMMA 3.6. *Let $U \subset \partial X$ be a nonempty open subset. Assume that $\lambda \in \mathbf{C}$ satisfies $\lambda \notin -\rho - \mathbf{N}_0 \cup -\mathbf{N}$. If $\phi \in C^{-\omega}(\partial X, V(\lambda))$ satisfies $\mathrm{res}_U(\phi) = 0$ and $\mathrm{res}_U \circ \hat{J}_\lambda(\phi) = 0$, then $\phi = 0$ (note that $(\hat{J}_\lambda\phi)_{|U}$ is defined even if $\hat{J}_\lambda$ has a pole).*

*Proof.* We reduce the proof to the case $\mathrm{Re}(\lambda) \geq 0$ by replacing $\phi$ by $\hat{J}_\lambda(\phi)$, if $\mathrm{Re}(\lambda) < 0$. We can do so because $\lambda \notin -\rho - \mathbf{N}_0 \cup -\mathbf{N}$ and thus $\hat{J}_\lambda$ and $\hat{J}_{-\lambda}$ are regular and bijective (see Lemma 4.13 below).

We now assume that $\mathrm{Re}(\lambda) \geq 0$. Then the Poisson transform $P_\lambda$ is a bijection between $C^{-\omega}(\partial X, V(\lambda))$ and $\ker(D)$, where $D := \Delta_X - \rho^2 + \lambda^2$. Since $D$ is elliptic with real analytic coefficients $P_\lambda\phi$ is a real analytic function.

We argue by contradiction and assume that $\phi \not\equiv 0$. Without loss of generality we can assume that $M \in U$. Since $P_\lambda\phi$ is real analytic and not identically zero the expansion (20) has nontrivial terms. Let $m$ be the smallest integer such that $\psi_m \neq 0$ near $M$, where $\psi_0 := \hat{J}_\lambda \phi$. The contradiction will be obtained by showing that $m = 0$.

With respect to the coordinates $k, a$ the operator $D$ has the form $D = D_0 + a^{-\alpha_1} R(a, k)$, where $D_0$ is a constant coefficient operator on $A$ and $R$ is a differential operator with coefficients which remain bounded if $a \to \infty$ (see [19, Ch. IV, §5, (8)]). Moreover, it is known that $D_0$ coincides with the $\bar{N}$-radial part of $D$.

We consider the $\bar{N}$-invariant function $f \in C^\infty(X)$ defined by $f(\bar{n}a) := a^{-(\lambda+\rho+m)}$. Since $D$ annihilates the asymptotic expansion (20) we have $Df = D_0 f = 0$. On the other hand, $f$ satisfies $(\Delta_X - \rho^2 + (\lambda+m)^2)f = 0$. Hence $(\lambda + m)^2 = \lambda^2$. Since $\mathrm{Re}(\lambda) \geq 0$ we conclude that $m = 0$. $\qquad\square$

The following corollary is an immediate consequence of Lemma 3.6 and Lemma 3.1.

COROLLARY 3.7. *If $\lambda \notin -\rho - \mathbf{N}_0 \cup -\mathbf{N}$, then*

$$\mathrm{res}_\Omega \circ \hat{J}_\lambda : C^{-\omega}(\Lambda, V(\lambda)) \to C^\omega(\Omega, V(-\lambda))$$

*is injective.*

LEMMA 3.8. *If $\lambda > 0$, then the order of a pole of $\mathrm{ext}_\mu$ at $\lambda$ is at most 1.*

*Proof.* Let $f_\mu \in C^\omega(B, V_B(\mu))$, $\mu \in \mathbf{C}$, be a holomorphic family such that $\mathrm{ext}_\mu(f_\mu)$ has a pole of order $l \geq 1$ at $\mu = \lambda$, $\lambda > 0$. We assume that $l \geq 2$ and argue by contradiction. Let $0 \neq \phi \in {}^\Gamma C^{-\omega}(\partial X, V(\lambda))$ be the leading singular part of $\mathrm{ext}_\mu(f_\mu)$ at $\mu = \lambda$. Because of $\mathrm{res}_\mu \circ \mathrm{ext}_\mu = \mathrm{id}$ we have $\mathrm{res}_\lambda \phi = 0$ and hence $\phi \in {}^\Gamma C^{-\omega}(\Lambda, V(\lambda))$.

By Lemma 3.2 for any compact subset $F \subset \Omega$ there exist constants $C_1, C_2, C_3$ such that for $a \gg 0$, $k \in FM$, $\mu$ near $\lambda$, $0 < \mu < \lambda$,

$$|P_\lambda \phi(ka)| \leq C_1 a^{-\lambda-\rho}$$

and

$$\begin{aligned}
(22) \qquad |(\mu - \lambda)^l P_\mu \mathrm{ext}_\mu(f_\mu)(ka)| &\leq C_2((\lambda-\mu)^l a^{\mu-\rho} + a^{-\rho}) \\
&\leq C_3(1 + \log a)^{-l} a^{\lambda-\rho} \ .
\end{aligned}$$

In particular, $P_\lambda \phi \in L^2(Y)$. Since $l \geq 2$ we obtain by Lebesgue's theorem of dominated convergence

$$\begin{aligned}
||P_\lambda \phi||^2 &= \langle \lim_{\substack{\mu \to \lambda \\ \mu < \lambda}} (\mu - \lambda)^l P_\mu \mathrm{ext}_\mu(f_\mu), P_\lambda \phi \rangle_{L^2(Y)} \\
&= \lim_{\substack{\mu \to \lambda \\ \mu < \lambda}} (\mu - \lambda)^l \langle P_\mu \mathrm{ext}_\mu(f_\mu), P_\lambda \phi \rangle_{L^2(Y)} \ .
\end{aligned}$$

On the other hand the estimates (22) allow partial integration, and we obtain for $\mu < \lambda$,

$$\begin{aligned}
\langle P_\mu \mathrm{ext}_\mu(f_\mu), P_\lambda \phi \rangle &= \frac{1}{\lambda^2 - \mu^2} \langle (\Delta_Y - \rho^2 + \lambda^2) P_\mu \mathrm{ext}_\mu(f_\mu), P_\lambda \phi \rangle \\
&= \frac{1}{\lambda^2 - \mu^2} \langle P_\mu \mathrm{ext}_\mu(f_\mu), (\Delta_Y - \rho^2 + \lambda^2) P_\lambda \phi \rangle \\
&= 0 \ .
\end{aligned}$$

Hence $\|P_\lambda \phi\|_{L^2(Y)} = 0$, which is a contradiction. We conclude that $l = 1$. $\quad\square$

For $\lambda \in \imath\mathbf{R}$ there is a conjugate-linear pairing $V_\lambda \otimes V_\lambda \to V_{-\rho}$ and hence a natural $L^2$-scalar product $C^\infty(B, V_B(\lambda)) \otimes C^\infty(B, V_B(\lambda)) \to \mathbf{C}$. Let $L^2(B, V_B(\lambda))$ be the associated Hilbert space. Using Lemma 2.18 (iv), we see that the adjoint $S_\lambda^*$ with respect to this Hilbert space structure is just $S_{-\lambda}$.

LEMMA 3.9.   *If* $\mathrm{Re}(\lambda) = 0$, *then* $S_\lambda$ *is regular and unitary.*

*Proof.* Assume that $S_\lambda$ is regular. If $f \in C^\infty(B, V_B(\lambda))$, then by Lemma 2.18 (iii),

$$\|S_\lambda f\|_{L^2(B, V_B(\lambda))}^2 = \langle S_{-\lambda} \circ S_\lambda f, f \rangle_{L^2(B, V_B(\lambda))} = \|f\|_{L^2(B, V_B(\lambda))}^2 \ .$$

Thus $S_\lambda$ is unitary.

Meromorphy of $S_\lambda$ and $\|S_\lambda\| = 1$ for all $\lambda \in \imath\mathbf{R}$ with $S_\lambda$ regular imply regularity of $S_\lambda$ for all $\lambda \in \imath\mathbf{R}$. $\quad\square$

LEMMA 3.10.   *If* $\mathrm{Re}(\lambda) = 0$, $\lambda \neq 0$, *then*

$$\mathrm{ext}_\lambda : C^{-\omega}(B, V_B(\lambda)) \to {}^\Gamma C^{-\omega}(\partial X, V(\lambda))$$

*is regular.*

*Proof.* By Lemma 3.9 the scattering matrix $S_\lambda$ is regular for $\mathrm{Re}(\lambda) = 0$. Since $c(-\lambda)$ is regular for $\mathrm{Re}(\lambda) = 0$, $\lambda \neq 0$, the same holds true for $\hat{S}_\lambda$. Recall that $\hat{S}_\lambda = \mathrm{res}_{-\lambda} \circ \hat{J}_\lambda \circ \mathrm{ext}_\lambda$. Corollary 3.7 and the fact that the leading singular part of $\mathrm{ext}_\lambda$ maps to hyperfunctions supported on the limit set $\Lambda$ show that a pole of $\mathrm{ext}_\lambda$ would necessarily imply a singularity of $\hat{S}_\lambda$. Thus $\mathrm{ext}_\lambda$ is regular for $\mathrm{Re}(\lambda) = 0$, $\lambda \neq 0$. $\quad\square$

LEMMA 3.11.   (i) *At* $\lambda = 0$ *the family* $\mathrm{ext}_\lambda$ *has at most a first order pole.*

(ii) $\mathrm{ext}_0$ *is regular if and only if* $S_0 = \mathrm{id}$.

(iii) *The operator* $S_0 - \mathrm{id}$ *is finite-dimensional.*

*Proof.* Since $S_0$ is regular and unitary by Lemma 3.9 and $c(\lambda)$ has a first order pole at $\lambda = 0$ the unnormalized scattering matrix $\hat{S}_\lambda$ has a first order pole at $\lambda = 0$, too. The principal series representation of $G$ on $C^\infty(\partial X, V(0))$ is irreducible ([25, Cor. 14.30] or [20, Ch.VI, Thm. 3.6]). Thus $J_0 = \mathrm{id}$ by our choice of normalization.

We expand $\hat{J}_\lambda = (\mathrm{res}_{\lambda=0} c(-\lambda))\,\lambda^{-1} + J^0_\lambda$ and $\mathrm{ext}_\lambda = \sum_{k<0} \lambda^k \mathrm{ext}^k + \mathrm{ext}^0_\lambda$, where $J^0_\lambda$ and $\mathrm{ext}^0_\lambda$ depend holomorphically on $\lambda$ near $\lambda = 0$. Since for $k < 0$ we have $\mathrm{supp}(\mathrm{ext}^k(f)) \subset \Lambda$ for all $f \in C^{-\omega}(B, V_B(0))$ Corollary 3.7 implies that $\mathrm{res}_\Omega \circ J^0_0 \circ \mathrm{ext}^k = 0$ if and only if $\mathrm{ext}^k = 0$. Since $\hat{S}_\lambda$ has a first order pole it follows that $\mathrm{ext}^k = 0$ for $k \leq -2$, hence (i). We now consider the residue of $\hat{S}_\lambda$.

$$
\begin{aligned}
\mathrm{res}_{\lambda=0} c(-\lambda) S_0 &= \mathrm{res}_{\lambda=0} \hat{S}_\lambda \\
&= \mathrm{res}_{\lambda=0} \mathrm{res}_{-\lambda} \circ \hat{J}_\lambda \circ \mathrm{ext}_\lambda \\
&= \mathrm{res}_{\lambda=0} c(-\lambda)\mathrm{id} + \mathrm{res}_0 \circ J^0_0 \circ \mathrm{ext}^{-1} \ .
\end{aligned}
$$

We conclude that $\mathrm{res}_{\lambda=0} c(-\lambda)(S_0 - \mathrm{id}) = \mathrm{res}_0 \circ J^0_0 \circ \mathrm{ext}^{-1}$. The right-hand side of this equation is a finite-dimensional operator, and (iii) follows. Again by Corollary 3.7 we conclude that $\mathrm{ext}^{-1} = 0$ if and only if $S_0 = \mathrm{id}$. This is assertion (ii). $\qquad\square$

LEMMA 3.12.    *If $\mathrm{Re}(\lambda) \geq 0$, then the residue of $\mathrm{ext}_\mu : C^{-\omega}(B, V_B(\mu)) \to {}^\Gamma C^{-\omega}(\partial X, V(\mu))$ at $\mu = \lambda$ spans ${}^\Gamma C^{-\omega}(\Lambda, V(\lambda))$. In particular, $\mathrm{ext}_\mu$ is regular at $\lambda$ if and only if ${}^\Gamma C^{-\omega}(\Lambda, V(\lambda)) = 0$.*

*Proof.* By Lemmas 3.8, 3.11 and 3.10 for $\mathrm{Re}(\mu) \geq 0$ the family $\mathrm{ext}_\mu$ has poles of at most first order. Thus $\mathrm{res}_{\mu=\lambda} \mathrm{ext}_\mu$ has values in ${}^\Gamma C^{-\omega}(\Lambda, V(\lambda))$. Let $d := \dim \mathrm{im}\, \mathrm{res}_{\mu=\lambda} \mathrm{ext}_\mu \leq \dim {}^\Gamma C^{-\omega}(\Lambda, V(\lambda))$. The meromorphic identity $\mathrm{res}_\mu \circ \mathrm{ext}_\mu = \mathrm{id}$ now implies that $d \geq \dim \mathrm{coker}\, \mathrm{res}_\lambda$. By Proposition 3.4 and Corollary 3.7,

$$
\begin{aligned}
\dim {}^\Gamma C^{-\omega}(\Lambda, V(\lambda)) &\geq d \geq \dim \mathrm{coker}\, \mathrm{res}_\lambda \\
&= \dim(\mathrm{im}\, \mathrm{res}_\lambda)^\perp \geq \dim {}^\Gamma C^{-\omega}(\Lambda, V(\lambda)) \ .
\end{aligned}
$$

It follows that $d = \dim {}^\Gamma C^{-\omega}(\Lambda, V(\lambda))$. This proves the lemma. $\qquad\square$

COROLLARY 3.13.    *The set $\mathcal{F} := \{\lambda \in \mathbf{C}\,|\,\mathrm{Re}(\lambda) \geq 0, {}^\Gamma C^{-\omega}(\Lambda, V(\lambda)) \neq 0\}$ is finite and contained in $[0, \delta_\Gamma]$. Moreover $\mathcal{F}$ is the set of singularities of $\mathrm{ext}_\lambda$ for $\mathrm{Re}(\lambda) \geq 0$. If $\lambda \in \mathcal{F}$ then $\mathrm{ext}_\mu$ has a first order pole at $\mu = \lambda$.*

*Proof.* We combine Lemmas 2.13, 3.10, and 3.12. $\qquad\square$

# 4. Cohomology

## 4.1. *Hyperfunctions with parameters.*

In this subsection we compare different definitions of holomorphic families of hyperfunctions. Let $M$ be a real analytic manifold and $K \subset M$ be a compact subset. Then the space $C^{-\omega}(K)$ of hyperfunctions on $M$ with support in $K$ has a natural Fréchet topology as the strong dual of the space of germs of real analytic functions at $K$. Thus we can consider holomorphic families of hyperfunctions in $C^{-\omega}(K)$ as in subsection 2.2. If $U \subset \mathbf{C}$, then by $\mathcal{O}C^{-\omega}(K)(U)$ we denote the space of holomorphic functions from $U$ with values in $C^{-\omega}(K)$.

Hyperfunctions on $M$ form a flabby sheaf $\mathcal{B}_M$. Thus we can consider the space of hyperfunctions $\mathcal{B}_M(V)$ on an open subset $V \subset M$. This space does not carry a natural topology. In order to define holomorphic families of hyperfunctions on $V$ we have to follow a different approach. Let $M_c$ be a complex neighbourhood of $M$ and $\mathcal{O}_{M_c}$ be the sheaf of holomorphic functions on $M_c$. Then by definition $\mathcal{B}_M = \mathcal{H}^n_M(\mathcal{O}_{M_c})$, where $n = \dim(M)$ and $\mathcal{H}^*_M(\mathcal{O}_{M_c})$ denotes the relative cohomology sheaf [27, p. 192ff]. It is a theorem (see [27, p. 219]) that $\{\phi \in \mathcal{B}_M(M) \,|\, \mathrm{supp}(\phi) \subset K\} = C^{-\omega}(K)$ in a natural manner. If $V \subset M$ is open, then a hyperfunction $f \in \mathcal{B}_M(V)$ can be represented (in a non-unique way) as a locally finite sum of hyperfunctions with compact support contained in $V$.

Using the approach via relative cohomology we can define the sheaf of holomorphic families of hyperfunctions as follows. Consider the embedding $\mathbf{C} \times M \hookrightarrow \mathbf{C} \times M_c$. Then we define the sheaf $\mathcal{O}\mathcal{B}_M := \mathcal{H}^n_{\mathbf{C} \times M}(\mathcal{O}_{\mathbf{C} \times M_c})$ on $\mathbf{C} \times M$. If $U \subset \mathbf{C}$, $V \subset M$ are open, then $\mathcal{O}\mathcal{B}_M(V)(U) := \mathcal{O}\mathcal{B}_M(U \times V)$ is by definition the space of holomorphic families of hyperfunctions on $V$ parametrized by $U$.

An important consequence of this definition is the following property. Let $T \subset \mathbf{C}$ be Stein. Then the prescription

$$M \supset V \mapsto \mathcal{O}\mathcal{B}_M(V)(T)$$

defines a flabby sheaf on $M$ (see [33, §5]).

Let $\bar{\partial}$ be the Cauchy-Riemann operator acting on the first variable of $\mathbf{C} \times M$. Then we have an alternative description of the sheaf $\mathcal{O}\mathcal{B}_M$ as the solution sheaf of $\bar{\partial}$; i.e., $\mathcal{O}\mathcal{B}_M(W) = \{\phi \in \mathcal{B}_{\mathbf{C} \times M}(W) \,|\, \bar{\partial}\phi = 0\}$ (see [27, p. 308]). If $K \subset M$ is compact, then we can consider

$$\mathcal{O}\mathcal{B}_M(K)(U) := \{\phi \in \mathcal{O}\mathcal{B}_M(M)(U) \,|\, \mathrm{supp}(\phi) \subset U \times K\} \,.$$

LEMMA 4.1. *There is a natural isomorphism* $\mathcal{O}C^{-\omega}(K)(U) \cong \mathcal{O}\mathcal{B}_M(K)(U)$.

*Proof.* In this proof we employ the description of $\mathcal{O}\mathcal{B}_M$ as the solution sheaf of the partial $\bar{\partial}$-equation. First we define a map $\Phi : \mathcal{O}C^{-\omega}(K)(U) \to \mathcal{O}\mathcal{B}_M(K)(U)$. Let $\{U_\alpha\}$ be a locally finite covering of $U$ and let $\{\chi_\alpha\}$ be a subordinated partition of unity. If $f \in \mathcal{O}C^{-\omega}(K)(U)$, then $\chi_\alpha f$ can be considered as an analytic functional on $\bar{U} \times K$ with support in $\mathrm{supp}(\chi_\alpha) \times K$ as follows: For $\phi \in C^\omega(\bar{U} \times K)$ we set

$$\langle \chi_\alpha f, \phi \rangle := \int_U \chi_\alpha(\lambda)\langle f_\lambda, \phi_\lambda \rangle d\lambda \,,$$

where $\phi_\lambda := \phi(\lambda, .) \in C^\omega(K)$.

Then we define $\Phi(f) \in \mathcal{B}_{\mathbf{C} \times M}(U \times K)$ to be the hyperfunction represented by the locally finite sum $\sum_\alpha \chi_\alpha f \in \mathcal{B}_{\mathbf{C} \times M}(U \times M)$ of hyperfunctions with compact support. Since the functions $\chi_\alpha$ form a partition of unity on $U$ and $f_\lambda$ depends holomorphically on $\lambda$ it is easy to see that $\bar{\partial} \sum_\alpha \chi_\alpha f = 0$. Thus indeed $\Phi(f) \in \mathcal{OB}_M(U \times K)$.

We now construct the inverse $\Psi : \mathcal{OB}_M(K)(U) \rightarrow \mathcal{O}C^{-\omega}(K)(U)$. Let $F \in \mathcal{OB}_M(K)(U)$. Since $F \in \mathcal{B}_{\mathbf{C} \times M}(U \times M)$ it satisfies $\bar{\partial}F = 0$. Therefore the specialization $F(\lambda, .) \in \mathcal{B}_M(M)$ is defined for all $\lambda \in U$. Clearly $\operatorname{supp}(F(\lambda, .)) \subset K$ and we define $\Psi(f)_\lambda := F(\lambda, .)$. We must show that $\Psi$ is well-defined, i.e. that $U \ni \lambda \mapsto \Psi(F)_\lambda \in C^{-\omega}(K)$ is a holomorphic function.

From the arguments of the proof of [22, Thm. 4.4], there exists a local elliptic operator $J$ on $\mathbf{C} \times M$ of possible infinite order acting on the second variable and a function $u \in C^\infty(U \times M)$ satisfying $\bar{\partial}u = 0$ such that $Ju = F$. By [22, Thm. 4.3], we can view $U \ni \lambda \mapsto u_\lambda := u(\lambda, .)$ as a holomorphic function on $U$ with values in $C^\infty(M)$.

Let $\chi \in C_c^\infty(M)$ be a cut-off function with $\chi \equiv 1$ in a neighbourhood of $K$. Since $J$ defines a continuous operator from $C^{-\omega}(\operatorname{supp}(\chi))$ into itself the function $\lambda \mapsto J(\chi u_\lambda)$ is holomorphic on $U$ with values in $C^{-\omega}(\operatorname{supp}(\chi))$. We can write $J(\chi u_\lambda) = F(\lambda, .) + F_1(\lambda)$, where $\operatorname{supp}(F_1)(\lambda) \subset \operatorname{supp}(d\chi)$. Since $K$ and $\operatorname{supp}(d\chi)$ are separated we have a continuous decomposition $C^{-\omega}(K \cup \operatorname{supp}(d\chi)) = C^{-\omega}(K) \oplus C^{-\omega}(\operatorname{supp}(d\chi))$. Thus $U \ni \lambda \mapsto F(\lambda, .)$ is a holomorphic function with values in $C^{-\omega}(K)$.

It remains to show that $\Phi$ and $\Psi$ are inverse to each other. Let $F \in \mathcal{OB}_M(K)(U)$. If $\chi \in C_c^\infty(U)$ is considered as a function on $U \times M$, then since $\bar{\partial}F = 0$ the product $\chi F \in \mathcal{B}_{\mathbf{C} \times M}(U \times M)$ is well-defined. Since $\operatorname{supp}(\chi F) \subset \operatorname{supp}(\chi) \times K$ is compact, $\chi F \in C^{-\omega}(\operatorname{supp}(\chi) \times K)$. Using the partition of unity $\{\chi_\alpha\}$ introduced above we can write $F$ as a locally finite sum of analytic functionals $F = \sum_\alpha \chi_\alpha F$.

Now we show $\Phi \circ \Psi = \operatorname{id}$. We have $\Psi(F)_\lambda = F(\lambda, .)$. Then $\Phi \circ \Psi(F)$ is given by $\sum_\alpha \chi_\alpha F = F$. In order to prove that $\Psi \circ \Phi = \operatorname{id}$, note that $\Phi(f)(\lambda, .) = \sum_\alpha \chi_\alpha(\lambda) f_\lambda = f_\lambda$. This finishes the proof of Lemma 4.1. $\qquad\square$

Fix $\lambda \in \mathbf{C}$. If $V \subset M$ is open, then we define

$$\mathcal{O}_\lambda C^{-\omega}(V) := \varinjlim \mathcal{OB}_M(V)(U) ,$$

where $U$ runs over all open neighbourhoods of $\lambda$.

LEMMA 4.2. *The prescription $M \supset V \mapsto \mathcal{O}_\lambda C^{-\omega}(V)$ defines a flabby sheaf $\mathcal{O}_\lambda \mathcal{B}_M$ on $M$.*

*Proof.* There exists a fundamental sequence of Stein neighbourhoods $\{T_n\}$ of $\lambda$. We can write $\mathcal{O}_\lambda C^{-\omega}(V) = \varinjlim_n \mathcal{OB}_M(V)(T_n)$. Recall that $M \supset V \mapsto$ $\mathcal{OB}_M(V)(T_n)$ defines a flabby sheaf on $M$. A direct limit of flabby sheaves is again a flabby sheaf. $\qquad\square$

For later reference we need the following result. Let $\phi \in C^\omega(M)$ be such that 0 is a regular value. Then $N := \{\phi = 0\} \subset M$ is a real analytic submanifold of codimension one. Let $i : N \hookrightarrow M$ denote the embedding. Then there is a natural morphism of sheaves $i_*\mathcal{O}_\lambda\mathcal{B}_N \to \mathcal{O}_\lambda\mathcal{B}_M$. Multiplication by $\phi$ induces a morphism $\phi : \mathcal{O}_\lambda\mathcal{B}_M \to \mathcal{O}_\lambda\mathcal{B}_M$.

LEMMA 4.3. *The sequence of sheaves*

$$(23) \qquad 0 \to i_*\mathcal{O}_\lambda\mathcal{B}_N \to \mathcal{O}_\lambda\mathcal{B}_M \xrightarrow{\phi} \mathcal{O}_\lambda\mathcal{B}_M \to 0$$

*is exact.*

*Proof.* Let $M_c$ be a complex neighbourhood of $M$ such that $\phi$ extends holomorphically to $M_c$ and 0 remains to be a regular value of this extension. Then $N_c := \{\phi = 0\}$ is a complex neighbourhood of $N$. We use $i$ to denote the embedding $N_c \hookrightarrow M_c$, too. There is a natural exact sequence of sheaves

$$(24) \qquad 0 \to \mathcal{O}_{\mathbf{C} \times M_c} \xrightarrow{\phi} \mathcal{O}_{\mathbf{C} \times M_c} \to i_*\mathcal{O}_{\mathbf{C} \times N_c} \to 0 \ .$$

Since $N_c \subset M_c$ is a locally closed submanifold, $i_* = i_!$ is an exact functor (see [24, Prop. 2.5.4]). It follows that

$$\mathcal{H}^*_{\mathbf{C} \times M}(i_*\mathcal{O}_{\mathbf{C} \times N_c}) = i_*\mathcal{H}^*_{\mathbf{C} \times (M \cap N_c)}(\mathcal{O}_{\mathbf{C} \times N_c}) = i_*\mathcal{H}^*_{\mathbf{C} \times N}(\mathcal{O}_{\mathbf{C} \times N_c}) \ .$$

By [27, p. 26], for $p \neq n$ and $q \neq n-1$ we have $\mathcal{H}^p_{\mathbf{C} \times M}(\mathcal{O}_{\mathbf{C} \times M_c}) = 0$ and $\mathcal{H}^q_{\mathbf{C} \times N}(\mathcal{O}_{\mathbf{C} \times N_c}) = 0$. The long exact sequence in cohomology associated with (24) gives

$$0 \to i_*\mathcal{H}^{n-1}_{\mathbf{C} \times N}(\mathcal{O}_{\mathbf{C} \times N_c}) \to \mathcal{H}^n_{\mathbf{C} \times M}(\mathcal{O}_{\mathbf{C} \times M_c}) \xrightarrow{\phi} \mathcal{H}^n_{\mathbf{C} \times M}(\mathcal{O}_{\mathbf{C} \times M_c}) \to 0 \ ,$$

or equivalently

$$0 \to i_*\mathcal{OB}_N \to \mathcal{OB}_M \xrightarrow{\phi} \mathcal{OB}_M \to 0 \ .$$

Taking germs at $\lambda \in \mathbf{C}$ we obtain (23) and conclude exactness. $\qquad\square$

### 4.2. *Acyclic resolutions.*

The main goal of this paper is to compute the cohomology groups $H^*(\Gamma, \mathcal{M})$, where the coefficients $\mathcal{M}$ are certain $\Gamma$-modules (i.e. complex representations of $\Gamma$). A $\Gamma$-module $\mathcal{M}$ is called $\Gamma$-acyclic if and only if $H^p(\Gamma, \mathcal{M}) = 0$ for all $p \geq 1$. A $\Gamma$-acyclic resolution of $\mathcal{M}$ is a complex

$$\mathcal{C}^\cdot : 0 \to C^0 \to C^1 \to C^2 \to \dots \ ,$$

of $\Gamma$-modules with $H^0(\mathcal{C}^\cdot) = \mathcal{M}$, $H^p(\mathcal{C}^\cdot) = 0$, $p \geq 1$, where all $C^i$ are $\Gamma$-acyclic. Standard homological algebra gives an isomorphism

$$H^*(\Gamma, \mathcal{M}) = H^*(^\Gamma \mathcal{C}^\cdot) \ ,$$

where

$$^{\Gamma}\mathcal{C}^{\cdot} : 0 \to {}^{\Gamma}C^0 \to {}^{\Gamma}C^1 \to {}^{\Gamma}C^2 \to \dots \; .$$

Using suitable trivializations of the families $V(\lambda)$, $V_B(\lambda)$ we can carry over the results of subsection 4.1 to hyperfunction sections of $V(\lambda)$ and $V_B(\lambda)$.

If $V \subset \partial X$ is open, then let $\mathcal{O}_\lambda C^{-\omega}(V)$ denote the space of germs at $\lambda$ of holomorphic families of hyperfunction sections of $V(\lambda)$ on $V$. By Lemma 4.2 the prescription $\partial X \supset V \mapsto \mathcal{O}_\lambda C^{-\omega}(V)$ defines a flabby sheaf. Let $\mathrm{res}_\Omega : \mathcal{O}_\lambda C^{-\omega}(\partial X) \to \mathcal{O}_\lambda C^{-\omega}(\Omega)$ be the restriction. Then we define $\mathcal{O}_\lambda C^{-\omega}(\Lambda)$ := ker $\mathrm{res}_\Omega$.

Note that $\Lambda$ and $\partial X$ are compact. By Lemma 4.1 we can identify the spaces $\mathcal{O}_\lambda C^{-\omega}(\Lambda)$ and $\mathcal{O}_\lambda C^{-\omega}(\partial X)$ with the spaces of germs at $\lambda$ of holomorphic functions $\mu \to f_\mu \in C^{-\omega}(\Lambda, V(\mu))$ and $\mu \to f_\mu \in C^{-\omega}(\partial X, V(\mu))$ in the sense of subsection 2.2.

LEMMA 4.4.  *The following complex of $\Gamma$-modules*

(25)        $$0 \to \mathcal{O}_\lambda C^{-\omega}(\Lambda) \to \mathcal{O}_\lambda C^{-\omega}(\partial X) \overset{\mathrm{res}_\Omega}{\to} \mathcal{O}_\lambda C^{-\omega}(\Omega) \to 0$$

*is exact.*

*Proof.* The restriction $\mathrm{res}_\Omega$ is surjective since the sheaf $V \subset \partial X \mapsto \mathcal{O}_\lambda C^{-\omega}(V)$ is flabby. Exactness at the other places holds by definition. □

LEMMA 4.5.  *The $\Gamma$-module $\mathcal{O}_\lambda C^{-\omega}(\Omega)$ is $\Gamma$-acyclic.*

*Proof.* The group $\Gamma$ acts freely and properly on $\Omega$. The $\Gamma$-module $\mathcal{O}_\lambda C^{-\omega}(\Omega)$ is the space of sections of a $\Gamma$-equivariant flabby sheaf on $\Omega$. Thus the assertion can be shown by repeating the arguments of the proof of [4, Lemma 2.6]. □

The proof of the following proposition will occupy the remainder of this subsection.

PROPOSITION 4.6.  *If $\lambda \notin -\rho - \mathbf{N}_0$, then $\mathcal{O}_\lambda C^{-\omega}(\partial X)$ is $\Gamma$-acyclic.*

*Proof.* We will show that $H^p(\Gamma, \mathcal{O}_\lambda C^{-\omega}(\partial X)) = 0$ using a suitable acyclic resolution. If $U \subset \mathbf{C}$ is open, then define

$$\mathcal{O}C^\infty(X)(U) := \{f \in C^\infty(U \times X) \,|\, \bar{\partial} f = 0\} \; ,$$

where $\bar{\partial}$ is the partial Cauchy-Riemann operator acting on the first variable. By [22, Thm. 4.3], $\mathcal{O}C^\infty(X)(U)$ is the space of holomorphic functions from $U$ to the Fréchet space $C^\infty(X)$. We define

$$\mathcal{O}_\lambda C^\infty(X) := \varinjlim \mathcal{O}C^\infty(X)(U) \; ,$$

where $U$ runs over all neighbourhoods of $\lambda$.

We define the operator $A : \mathcal{O}C^\infty(X)(U) \to \mathcal{O}C^\infty(X)(U)$ by $(Af)_\mu = (\Delta_X - \rho^2 + \mu^2)f_\mu$, where $f_\mu = f(\mu, .)$. We use the same symbol $A$ in order to denote the induced operator on $\mathcal{O}_\lambda C^\infty(X)$.

Note that the Poisson transform $P_\mu$ comes as a holomorphic family of continuous maps $P_\mu : C^{-\omega}(\partial X, V(\mu)) \to C^\infty(X)$. Viewing $\mathcal{O}_\lambda C^{-\omega}(\partial X)$ as a space of germs at $\lambda$ of holomorphic functions with values in a topological vector space we can define a Poisson transform $P : \mathcal{O}_\lambda C^{-\omega}(\partial X) \to \mathcal{O}_\lambda C^\infty(X)$ by $(Pf)_\mu = P_\mu f_\mu$. Since $(\Delta_X - \rho^2 + \mu^2)P_\mu f_\mu = 0$ we have $A \circ P = 0$. We need to know that $P : \mathcal{O}_\lambda C^{-\omega}(\partial X) \to \{f \in \mathcal{O}_\lambda C^\infty(X) \mid Af = 0\}$ is surjective.

PROPOSITION 4.7. *Let $U \subset \mathbf{C}$ be open such that $U \cap -\rho - \mathbf{N}_0 = \emptyset$. Let $\mu \mapsto f_\mu$, $\mu \in U$, be a continuous (smooth, holomorphic) family of eigenfunctions, i.e., $Af = 0$. Then there exists a continuous (smooth, holomorphic) family of hyperfunctions $\beta f$, $(\beta f)_\mu \in C^{-\omega}(\partial X, V(\mu))$, such that $P(\beta f) = f$.*

*Proof.* In order to prove the proposition one may take any of the existing proofs of the pointwise surjectivity of $P_\mu$ ([18], [23], [46]). One has to control the dependence on $\mu$ of the construction of the inverse map $\beta_\mu$. We prefer to follow Helgason's proof [18] because it is the most elementary one.

Let $\hat{K}_M$ denote the set of equivalence classes of irreducible representations of $K$ with a nontrivial $M$-fixed vector. Let $\delta_p$ be the class represented by the representation of $K$ on the space of homogeneous harmonic polynomials of degree $p$ on $\mathbf{R}^n$. Then $\hat{K}_M = \{\delta_p \mid p \in \mathbf{N}_0\}$. Hence there is a decomposition

$$(26) \qquad f_\mu = \sum_{p=0}^\infty f_{\mu,p} \ ,$$

where $f_{\mu,p} \in \ker(\Delta_X - \rho^2 + \mu^2)$ transforms according to the $K$-type $\delta_p$. For $\mu \notin -\rho - \mathbf{N}_0$, any $K$-finite eigenfunction is the Poisson transform of a unique $K$-finite function on $\partial X$ (see e.g. [20, Ch.III, Thm. 6.1]); i.e.,

$$f_{\mu,p} = P_\mu \varphi_{\mu,p} \qquad \text{for a certain } \varphi_{\mu,p} \in C^\infty(\partial X) \ .$$

Here and in the following we identify $C^*(\partial X, V(\mu))$ with $C^*(\partial X)$ using a $K$-invariant volume form. Let $\psi_p$ be the unique $M$-spherical function on $\partial X$ of $K$-type $\delta_p$, and set $\Phi_{\mu,p} := P_\mu \psi_p$. Then we have ([18, Lemma 4.2])

$$(27) \qquad f_{\mu,p}(ka) = \varphi_{\mu,p}(k) \cdot \Phi_{\mu,p}(a) \ , \qquad k \in K, \ a \in A \ .$$

We want to show that for a continuous family $f$ the series $\mu \mapsto \beta_\mu f_\mu := \sum_p \varphi_{\mu,p}$ converges in the space of continuous functions from $U$ to $C^{-\omega}(\partial X)$. The resulting limit $\beta f$ then satisfies $P(\beta f) = f$. In fact, we will show the stronger result that $\beta : \{f \in C(U, C^\infty(X)) \mid Af = 0\} \to C(U, C^{-\omega}(\partial X))$ is continuous.

Let $\varphi \in C^{-\omega}(\partial X)$ and $\varphi = \sum_p \varphi_p$ be its decomposition with respect to $K$-types. For $0 < R < 1$ we define seminorms $||.||_R$ on $C^{-\omega}(\partial X)$ by

$$||\varphi||_R^2 := \sum_{p=0}^{\infty} R^{2p} ||\varphi_p||_{L^2(\partial X)}^2 \ .$$

It is not difficult to see that the seminorms $||.||_R$, $0 < R < 1$, define the topology on $C^{-\omega}(\partial X)$ (compare [18, Prop. 5.2]). Thus the topology on $C(U, C^{-\omega}(\partial X))$ is given by the seminorms

$$||\varphi||_{W,R} := \sup_{\mu \in W} ||\varphi_\mu||_R \ , \qquad W \subset U \text{ compact}, \ 0 < R < 1 \ .$$

The function $\Phi_{\mu,p}$ can be represented in terms of the hypergeometric function $F = {}_2F_1$ as follows ([18, p. 341]):

$$\Phi_{\mu,p}(a) = c_p(\mu)(1-r^2)^{\rho-\frac{1}{2}} r^p F(\mu + \frac{1}{2}, -\mu + \frac{1}{2}, p + \rho + \frac{1}{2}, \frac{r^2}{r^2-1}) \ ,$$

where $a \in [1, \infty)$, $r = \frac{a-1}{a+1}$ and $c_p(\mu) := \frac{\Gamma(\mu+\rho+p)}{\Gamma(\mu+\rho)} \frac{\Gamma(\rho+\frac{1}{2})}{\Gamma(\rho+\frac{1}{2}+p)}$.

Now fix $R$ and a compact set $W \subset U$. The following two estimates are crucial:

(i) There exist constants $C_1$ and $C_2$ such that for all $p \in \mathbf{N}_0$ and $\mu \in W$

$$\frac{1}{|c_p(\mu)|} \leq C_1 \frac{|\Gamma(\mu+\rho)|}{\Gamma(\rho+\frac{1}{2})}(1+p)^{\frac{1}{2}-\operatorname{Re}(\mu)} \leq C_2(1+p)^{\frac{1}{2}-\operatorname{Re}(\mu)} \ .$$

In particular,

$$\sup_{\mu \in W} \frac{1}{|c_p(\mu)|} \leq C(1+p)^{2k} \tag{28}$$

for some $C > 0$, $k \in \mathbf{N}_0$.

(ii) Set $y := \frac{R^2}{R^2-1}$. Then there exists $P \in \mathbf{N}_0$ such that for all $p \geq P$, $\mu \in W$,

$$|F(\mu + \frac{1}{2}, -\mu + \frac{1}{2}, p + \rho + \frac{1}{2}, y)| \geq \frac{1}{2} \ . \tag{29}$$

Assertion (i) follows from $\lim_{x \to \infty} x^\alpha \frac{\Gamma(x)}{\Gamma(x+\alpha)} = 1$ for real $x$ and $\alpha \in \mathbf{C}$ (see [11, p. 47]). In order to verify (ii) we estimate

$$\begin{aligned} \varrho_1(\mu, p, y) &:= F(\mu + \frac{1}{2}, -\mu + \frac{1}{2}, p + \rho + \frac{1}{2}, y) - 1 \\ &= F(-\mu + \frac{1}{2}, \mu + \frac{1}{2}, p + \rho + \frac{1}{2}, y) - 1 \end{aligned}$$

for $\mu$ varying in the compact set $\tilde{W} := (W \cup -W) \cap \{\mathrm{Re}(\mu) \geq 0\}$. By [11, p. 76] we have for $\mu \in \tilde{W}$ and $p > \sup_{\mu \in \tilde{W}}(\mathrm{Re}(\mu) - \rho)$,

$$
\begin{aligned}
|\varrho_1(\mu, p, y)| &\leq (|y| + 1)^{\mathrm{Re}(\mu) + 1}|\mu^2 - \frac{1}{2}| \cdot \\
&\quad \cdot \frac{\Gamma(\mathrm{Re}(\mu) + \frac{1}{2})}{|\Gamma(\mu + \frac{1}{2})|} \frac{\Gamma(p + \rho - \mathrm{Re}(\mu))}{|\Gamma(p + \rho - \mu)|} \frac{1}{p + \rho + \frac{1}{2}} \\
&\leq C\frac{1}{p + \rho + \frac{1}{2}} \ .
\end{aligned}
$$

If $P > \sup_{\mu \in \tilde{W}}(\mathrm{Re}(\mu) - \rho)$ is large enough, then $|\varrho_1(\mu, p, y)| \leq \frac{1}{2}$ for all $p \geq P$ and $\mu \in \tilde{W}$. Assertion (ii) follows.

Now we are able to estimate $\|\beta f\|_{W,R}$ for a continuous family of eigenfunctions $f$. Equations (26) and (27) imply that $Q_1$ defined by

$$
Q_1(f) := \sup_{\mu \in W}(\sum_{p=0}^{P-1} R^{2p}\|\varphi_{\mu,p}\|_{L^2(\partial X)}^2)^{\frac{1}{2}}
$$

is a continuous seminorm on $\{f \in C(U, C^\infty(X)) \mid Af = 0\}$. Using (i), (ii) and the fact that $\delta_p(C_K) = p(p + n - 2)\mathrm{id}$, where $C_K$ is the Casimir operator of $K$, we estimate the remainder

$$
\begin{aligned}
\sum_{p=P}^{\infty} &R^{2p}\|\varphi_{\mu,p}\|_{L^2(\partial X)}^2 \\
&\overset{(29)}{\leq} 4\sum_{p=P}^{\infty} R^{2p}|F(\mu + \frac{1}{2}, -\mu + \frac{1}{2}, p + \rho + \frac{1}{2}, \frac{R^2}{R^2 - 1})|^2|\,\|\varphi_{\mu,p}\|_{L^2(\partial X)}^2 \\
&\overset{(28)}{\leq} \frac{4C}{(1 - R^2)^{2\rho - 1}}\sum_{p=0}^{\infty}(1 + p)^{4k} \cdot \\
&\qquad \cdot \|c_p(\mu)(1 - R^2)^{\rho - \frac{1}{2}}R^p F(\mu + \frac{1}{2}, -\mu + \frac{1}{2}, p + \rho + \frac{1}{2}, \frac{R^2}{R^2 - 1})\varphi_{\mu,p}\|_{L^2(\partial X)}^2 \\
&\leq C'\int_K |(1 + C_K)^k f_\mu(ka)|^2 dk \ , \qquad a = \frac{1 + R}{1 - R} \ .
\end{aligned}
$$

We obtain

$$
\|\beta f\|_{W,R} \leq Q_1(f) + C'^{\frac{1}{2}}Q_2(f) \ ,
$$

where $Q_2$ is the continuous seminorm on $C(U, C^\infty(X))$ given by

$$
Q_2(f) := \sup_{\mu \in W}\left(\int_K |(1 + C_K)^k f_\mu(ka)|^2 dk\right)^{\frac{1}{2}} \ .
$$

Continuity of $\beta : \{f \in C(U, C^\infty(X)) \mid Af = 0\} \to C(U, C^{-\omega}(\partial X))$ follows.

It remains to discuss smooth and holomorphic families of eigenfunctions $f$. Let $\bar{\partial}$ be the Cauchy-Riemann operator acting on the $\mu$-variable. Let $f \in C^\infty(U, C^\infty(X))$ be such that $Af = 0$. Then $\bar{\partial}^k f$ is again a smooth family of eigenfunctions. Since $\Phi_{\mu,p}$ is holomorphic in $\mu$ it follows from (27) that

$$\beta(\bar{\partial}^l f) = \sum_{p=0}^{\infty} \bar{\partial}^l \varphi_{.,p}, \quad \text{for all } l \in \mathbf{N}_0 .$$

Hence $\sum_{p=0}^{\infty} \bar{\partial}^l \varphi_{.,p}$ converges for all $l \in \mathbf{N}_0$, and we conclude that

$$(30) \qquad \bar{\partial}^k(\beta f) = \bar{\partial}^k \sum_{p=0}^{\infty} \varphi_{.,p} = \sum_{p=0}^{\infty} \bar{\partial}^k \varphi_{.,p} = \beta(\bar{\partial}^k f).$$

Thus, if $f$ is holomorphic, then so is $\beta f$. If $f$ is smooth, then the right-hand side of (30) is continuous for all $k$. Hence, by elliptic regularity $\beta f$ has to be smooth. This completes the proof of the proposition. □

LEMMA 4.8. *If $\lambda \notin -\rho - \mathbf{N}_0$, then the complex*

$$0 \to \mathcal{O}_\lambda C^{-\omega}(\partial X) \xrightarrow{P} \mathcal{O}_\lambda C^\infty(X) \xrightarrow{A} \mathcal{O}_\lambda C^\infty(X) \to 0$$

*is exact.*

*Proof.* Since $P_\mu$ is injective for $\mu \notin -\rho - \mathbf{N}_0$ the map $P$ is injective. Exactness in the middle follows from Proposition 4.7. It remains to show that $A$ is surjective. This turns out to be quite complicated. We will need the following technical result.

PROPOSITION 4.9. *Let $(M, g)$ be a connected, noncompact, real analytic, Riemannian manifold, and let $\Delta_M$ be the associated Laplace operator. Let $r \in \mathbf{C}$ and $U \subset \mathbf{C}$ be open. Then the operator $D : C^\infty(U \times M) \to C^\infty(U \times M)$ given by $(Df)(\mu, m) := ((\Delta_M + r + \mu^2) f)(\mu, m)$ is surjective.*

*Proof.* We consider the adjoint operator ${}^tD : C_c^{-\infty}(U \times M) \to C_c^{-\infty}(U \times M)$. To prove surjectivity of $D$ it is sufficient to show that ${}^tD$ is injective and has closed range.

First we prove that ${}^tD$ is injective. Let $f \in C_c^{-\infty}(U \times M)$ satisfy ${}^tDf = 0$. We consider $f$ as a hyperfunction in $\mathcal{B}_{U \times M}(U \times M)$. The hyperplanes $\{\mu = \text{const}\}$ are non-characteristic for ${}^tD$. By a theorem of Sato [44], $f(\mu, m)$ contains the variable $m \in M$ as a real analytic parameter. Since $\text{supp}(f)$ is compact we conclude by [22, Thm. 1.5], that $f = 0$. This shows injectivity of ${}^tD$.

We now prove that the range of ${}^tD$ is closed. Let $f_i \in C_c^{-\infty}(U \times M)$ be a sequence such that ${}^tDf_i =: h_i$ converges to $h \in C_c^{-\infty}(U \times M)$. We have to find $f \in C_c^{-\infty}(U \times M)$ such that ${}^tDf = h$.

There is a compact subset $T \times W \subset U \times M$ such that $\text{supp}(h_i) \subset T \times W$ for all $i$. Let $\chi \in C_c^\infty(U)$ such that $\chi_{|T} \equiv 1$. Since ${}^tD$ is of order zero with respect to the first variable we have ${}^t D(\chi f_i) = \chi h_i = h_i$. Replacing $f_i$ by $\chi f_i$ and enlarging $T$ we can and will assume that $\mathrm{supp}(f_i) \subset T \times M$. We will show that $\lim_{i \to \infty} f_i =: f$ exists.

By enlarging $W$ we can assume that any connected component of the closure of $M \setminus W$ is noncompact. We claim that $\mathrm{supp}(f_i) \subset T \times W$ for all $i$. We will argue by contradiction. Let $p : \mathbf{C} \times M \to M$ denote the projection. Let $V$ be a connected component of $M \setminus W$ and assume that $(f_i)_{|U \times V} \neq 0$. Then by the unique continuation result [22, Thm. 1.5], $p(\mathrm{supp}(f_i)) \supset V$. This is impossible since $\mathrm{supp}(f_i)$ and thus $p(\mathrm{supp}(f_i))$ are compact. It follows that $(f_i)_{|U \times V} = 0$ for all connected components $V$ of $M \setminus W$.

Let $V \subset M$ be an open subset with smooth boundary $\partial V$ such that $W \subset V$. We choose a real nonnegative potential $0 \neq F \in C_c^\infty(V)$ with $\mathrm{supp}(F) \subset V \setminus W$. For small $t$ we consider the Dirichlet extension $\Delta_t$ of the operator $\Delta_M + tF$ on $L^2(V)$. Then $\Delta_t$ is self-adjoint with eigenvalues $0 < \kappa_1(t) < \kappa_2(t) \leq \kappa_3(t) \leq \ldots$ (counted with multiplicity) such that $\kappa_j(t) \to \infty$ as $j \to \infty$ for fixed $t$. Furthermore, the eigenvalues $\kappa_i$ are continuous, strictly increasing functions of $t$. We employ this fact in order to choose for any $\delta \in U$ a number $t_\delta$ such that $-(r + \delta^2) \notin \mathrm{spec}(\Delta_{t_\delta})$.

Let $H^s(V)$ be the scale of Sobolev spaces associated with $\Delta_0$. For $s \geq 0$ we have $H^s(V) := \mathrm{dom}(1 + \Delta_0)^{s/2} = \mathrm{dom}(1 + \Delta_t)^{s/2}$. If $s < 0$, then $H^s(V) := H^{-s}(V)'$. By elliptic regularity $\bigcap_{s \in \mathbf{R}} H^s(V) \subset C^\infty(V)$ and dually $C_c^{-\infty}(V) \subset \bigcup_{s \in \mathbf{R}} H^s(V)$. If $\chi \in C_c^\infty(V)$, then multiplication by $\chi$ defines a continuous map $\chi : C^\infty(V) \to \bigcap_{s \in \mathbf{R}} H^s(V)$.

If $d \notin \mathrm{spec}(\Delta_{t_\delta})$, then there is a neighbourhood $Z_\delta \subset \mathbf{C} \setminus \mathrm{spec}(\Delta_{t_\delta})$ of $d$ such that $Z_\delta \ni a \mapsto (\Delta_{t_\delta} - a)^{-1} : \bigcap_{s \in \mathbf{R}} H^s(V) \to \bigcap_{s \in \mathbf{R}} H^s(V)$ is a holomorphic family of continuous isomorphisms. We choose open neighbourhoods $U_\delta$ of $\delta$ such that $-(r + \mu^2) \in Z_\delta$ for all $\mu \in U_\delta$.

We choose a countable set of points $\delta_j \in U$ such that $\{U_{\delta_j}\}$ is a locally finite cover of $U$. Let $\{\xi_j\}$ be an associated smooth partition of unity. Set $U_j := U_{\delta_j}$, $t_j := t_{\delta_j}$ and choose $\chi \in C_c^\infty(V) \subset C_c^\infty(M)$ with $\chi_{|W} \equiv 1$. Then we define continuous maps

$$L_j : C^\infty(U_j \times M) \to C^\infty(U_j, \bigcap_{s \in \mathbf{R}} H^s(V)) \subset C^\infty(U_j \times V)$$

by $L_j(\phi)(\mu) := (\Delta_{t_j} + r + \mu^2)^{-1} \chi \phi_\mu$, where $\phi_\mu(m) = \phi(\mu, m)$. Note that $\xi_j h_i \in C_c^{-\infty}(U_j \times V)$. Thus we can define the distributions $\tilde{f_{i,j}} \in C_c^{-\infty}(U_j \times M)$ by $\langle \tilde{f_{i,j}}, \phi \rangle := \langle \xi_j h_i, L_j(\phi) \rangle$, for all $\phi \in C^\infty(U_j \times V)$; i.e., $\tilde{f_{i,j}} := {}^t L_j(\xi_j h_i)$.

We claim that $f_{i,j} = \xi_j f_i$. Indeed, we have for all $\phi \in C^\infty(U_j \times M)$

$$
\begin{aligned}
\langle \tilde{f_{i,j}}, \phi \rangle &= \langle \xi_j h_i, L_j(\phi) \rangle \\
&= \langle \xi_j {}^t D f_i, L_j(\phi) \rangle \\
&= \langle {}^t D(\xi_j f_i), L_j(\phi) \rangle
\end{aligned}
$$

$$
\begin{aligned}
&= \langle ({}^t D + t_j F)(\xi_j f_i), L_j(\phi) \rangle \\
&= \langle \xi_j f_i, (\Delta_{t_j} + r + \mu^2) L_j(\phi) \rangle \\
&= \langle \xi_j f_i, \chi \phi \rangle \\
&= \langle \xi_j f_i, \phi \rangle \; .
\end{aligned}
$$

Consider the set $I := \{ j \in \mathbf{N} \mid U_j \cap T \neq \emptyset \}$. By local finiteness of the cover $\{U_j\}$ and compactness of $T$ the set $I$ is finite. By the claim above we have $\tilde{f_{i,j}} = 0$ for $j \notin I$ and $f_i = \sum_{j \in I} \tilde{f_{i,j}} = \sum_{j \in I} {}^t L_j(\xi_j h_i)$. Since $h_i \to h$ as $i \to \infty$ and ${}^t L_j$ is continuous we conclude that $\lim_{i \to \infty} f_i =: f$ exists. This finishes the proof of Proposition 4.9. □

We now finish the proof of Lemma 4.8, proving surjectivity of

$$
A : \mathcal{O}_\lambda C^\infty(X) \to \mathcal{O}_\lambda C^\infty(X) \; .
$$

Let $[f] \in \mathcal{O}_\lambda C^\infty(X)$ be represented by $f \in C^\infty(U \times X)$, where $U$ is a precompact neighbourhood of $\lambda$ with smooth boundary such that $U \cap -\rho - \mathbf{N}_0 = \emptyset$ and $\bar{\partial} f = 0$. By Proposition 4.9 we can find $h \in C^\infty(U \times X)$ with $Ah = f$. We have $0 = \bar{\partial} Ah = A\bar{\partial} h$. Thus $\bar{\partial} h$ is a smooth family of eigenfunctions. By Proposition 4.7 there exists a smooth family $\beta(\bar{\partial} h) \in C^\infty(U, C^{-\omega}(\partial X))$ such that $P(\beta(\bar{\partial} h)) = \bar{\partial} h$. Here we have identified $C^{-\omega}(\partial X, V(\mu))$ with $C^{-\omega}(\partial X)$, as usual.

Let $\Delta_U$ be the Dirichlet extension of $\bar{\partial} \bar{\partial}^*$ on $L^2(U)$. Then $\Delta_U^{-1}$ is bounded and can be considered as a continuous operator $\Delta_U^{-1} : C_c^\infty(U) \to C^\infty(U)$. Using that $C_{(c)}^\infty(U, C^{-\omega}(\partial X)) = C_{(c)}^\infty(U) \hat{\otimes} C^{-\omega}(\partial X)$ we obtain an induced operator $Q = \Delta_U^{-1} \otimes \mathrm{id}$. Let $W \subset \bar{W} \subset U$ be a neighbourhood of $\lambda$ and choose $\chi \in C_c^\infty(U)$ such that $\chi_{|W} \equiv 1$. We define $k := \bar{\partial}^* Q(\chi \beta(\bar{\partial} h))$. Then $(\bar{\partial} k)_{|W} = (\beta(\bar{\partial} h))_{|W}$.

If we define $\phi \in C^\infty(W \times X)$ by $\phi := (h - Pk)_{|W \times X}$, then

$$
\bar{\partial} \phi = \bar{\partial}(h - Pk)_{|W \times X} = \bar{\partial} h_{|W \times X} - (P\bar{\partial})k_{|W \times X} = \bar{\partial} h_{|W \times X} - (P\beta \bar{\partial} h)_{|W \times X} = 0
$$

and $A\phi = A(h - Pk)_{|W \times X} = Ah_{|W \times X} = f_{|W \times X}$. We let $[\phi] \in \mathcal{O}_\lambda C^\infty(X)$ be represented by $\phi$ and obtain $A[\phi] = [f]$. This proves surjectivity of $A$ and thus Lemma 4.8. □

LEMMA 4.10. *The $\Gamma$-module $\mathcal{O}_\lambda C^\infty(X)$ is acyclic.*

*Proof.* If $V \subset X$ is open, then let $\mathcal{O}_\lambda C^\infty(V)$ be the space of germs at $\lambda$ of holomorphic functions with values in $C^\infty(V)$. The prescription $X \supset V \mapsto \mathcal{O}_\lambda C^\infty(V)$ defines a soft sheaf on $X$. Since $\Gamma$ acts freely and properly on $X$, and the $\Gamma$-module $\mathcal{O}_\lambda C^\infty(X)$ is the space of sections of a $\Gamma$-equivariant soft sheaf on $X$, we can argue as in the proof of [4, Lemma 2.4]. □

We return to the proof of Proposition 4.6. By Lemma 4.8 and Lemma 4.10

$$0 \to \mathcal{O}_\lambda C^\infty(X) \xrightarrow{A} \mathcal{O}_\lambda C^\infty(X) \to 0$$

is a $\Gamma$-acyclic resolution of $\mathcal{O}_\lambda C^{-\omega}(\partial X)$.

Since $\mathcal{O}_\lambda C^\infty(X) = \mathcal{O}_{-\lambda} C^\infty(X)$ we can assume that $\mathrm{Re}(\lambda) \geq 0$. Let $\mathcal{O}_\lambda C^\infty(Y) := {}^\Gamma \mathcal{O}_\lambda C^\infty(X)$ be the space of germs at $\lambda$ of holomorphic functions with values in $C^\infty(Y)$ and $A_Y : \mathcal{O}_\lambda C^\infty(Y) \to \mathcal{O}_\lambda C^\infty(Y)$ be the operator induced by $A$. Then

$$
\begin{aligned}
H^0(\Gamma, \mathcal{O}_\lambda C^{-\omega}(\partial X)) &= \ker(A_Y) \,, \\
H^1(\Gamma, \mathcal{O}_\lambda C^{-\omega}(\partial X)) &= \operatorname{coker}(A_Y) \,, \\
H^p(\Gamma, \mathcal{O}_\lambda C^{-\omega}(\partial X)) &= 0, \quad \text{for all} \quad p \geq 2 \,.
\end{aligned}
$$

To prove Proposition 4.6 it remains to show that $\operatorname{coker}(A_Y) = 0$ or, equivalently, that $A_Y$ is surjective. The proof of surjectivity of $A_Y$ is similar to the proof of surjectivity of $A$.

Let $[f] \in \mathcal{O}_\lambda C^\infty(Y)$ be represented by $f \in C^\infty(U \times Y)$ satisfying $\bar\partial f = 0$. By Proposition 4.9 we can find $h \in C^\infty(U \times Y)$ with $A_Y h = f$. We have $0 = \bar\partial A_Y h = A_Y \bar\partial h$. By $p^* : C^\infty(Y) \hookrightarrow C^\infty(X)$ we denote the inclusion given by the pull-back associated to the covering projection $p : X \to Y$. Then $A p^* \bar\partial h = 0$ and $\beta_\mu(p^* \bar\partial h)_\mu \in {}^\Gamma C^{-\omega}(\partial X, V(\mu))$ is defined. The family $\mu \mapsto \mathrm{res}_\mu \circ \beta_\mu(p^* \bar\partial h)_\mu \in C^{-\omega}(B, V_B(\mu))$ is a smooth family of analytic functionals.

Using a suitable holomorphic trivialization of the family of bundles $V_B(\mu)$ we can identify $C^{-\omega}(B, V_B(\mu))$ with $C^{-\omega}(B)$. We consider $\mathrm{res} \circ \beta(p^* \bar\partial h) \in C^\infty(U, C^{-\omega}(B))$. As in the proof of Lemma 4.8 we solve the $\bar\partial$-problem $\bar\partial k = \mathrm{res} \circ \beta(p^* \bar\partial h)$ for $k \in C^\infty(W, C^{-\omega}(B))$, where $W \subset \bar W \subset U$ is a smaller neighbourhood of $\lambda$.

Replacing $k$ by $k - k(\lambda)$, we can and will assume that $k(\lambda) = 0$. By our assumption $\mathrm{Re}(\lambda) \geq 0$ and Corollary 3.13, $\mathrm{ext}_\mu$ is regular at $\mu = \lambda$ or it has a pole of at most first order. Thus $\mathrm{ext}(k)$ is smooth on $W \setminus \{\lambda\}$ and bounded on $W$.

We define the family $W \ni \mu \mapsto \phi_\mu := h_\mu - P_\mu \circ \mathrm{ext}_\mu(k_\mu) \in C^\infty(Y)$. This family is bounded on $W$, and on $W \setminus \{\lambda\}$ it satisfies $\bar\partial \phi = 0$. We conclude that $\phi$ is a holomorphic function from $W$ to $C^\infty(Y)$. Moreover, $A_Y \phi = f_{|W \times Y}$. If we define $[\phi] \in \mathcal{O}_\lambda C^\infty(Y)$ to be the element which is represented by $\phi$, then $A_Y[\phi] = [f]$. This proves that $A_Y$ is surjective. The proof of Proposition 4.6 is now complete.

### 4.3. *Computation of* $H^*(\Gamma, \mathcal{O}_\lambda C^{-\omega}(\Lambda))$.

Let $\mathcal{O}_\lambda C^{-\omega}(B)$ be the space of germs at $\lambda$ of holomorphic families $\mu \mapsto f_\mu \in C^{-\omega}(B, V_B(\mu))$ and $\mathrm{res} : {}^\Gamma \mathcal{O}_\lambda C^{-\omega}(\partial X) \to \mathcal{O}_\lambda C^{-\omega}(B)$ be defined by $\mathrm{res}(f)_\mu := \mathrm{res}_\mu f_\mu$.

LEMMA 4.11.  *If* $\lambda \notin -\rho - \mathbf{N}_0$, *then* $H^p(\Gamma, \mathcal{O}_\lambda C^{-\omega}(\Lambda)) = 0$ *for all* $p \neq 1$ *and*

$$(31) \qquad H^1(\Gamma, \mathcal{O}_\lambda C^{-\omega}(\Lambda)) \cong \text{coker}(\text{res}) \ .$$

*Proof.* By Lemmas 4.4, 4.5, and Proposition 4.6 the complex

$$0 \to \mathcal{O}_\lambda C^{-\omega}(\partial X) \overset{\text{res}_\Omega}{\to} \mathcal{O}_\lambda C^{-\omega}(\Omega) \to 0$$

is a $\Gamma$-acyclic resolution of $\mathcal{O}_\lambda C^{-\omega}(\Lambda)$. If we employ the identifications

$$^\Gamma C^{-\omega}(\Omega, V(\mu)) = C^{-\omega}(B, V_B(\mu)) \ ,$$

then the corresponding complex of $\Gamma$-invariant vectors can be written in the form

$$0 \to {}^\Gamma\mathcal{O}_\lambda C^{-\omega}(\partial X) \overset{\text{res}}{\to} \mathcal{O}_\lambda C^{-\omega}(B) \to 0 \ .$$

We conclude that

$$\begin{aligned}
H^0(\Gamma, \mathcal{O}_\lambda C^{-\omega}(\Lambda)) &= \text{ker}(\text{res}) \ , \\
H^1(\Gamma, \mathcal{O}_\lambda C^{-\omega}(\Lambda)) &= \text{coker}(\text{res}) \ , \\
H^p(\Gamma, \mathcal{O}_\lambda C^{-\omega}(\Lambda)) &= 0, \quad \text{for all } p \geq 2 \ .
\end{aligned}$$

Since $\text{res}_\mu$ is injective for generic $\mu \in \mathbf{C}$ we have $\text{ker}(\text{res}) = 0$ (see [8, Prop. 6.11]). $\qquad \square$

The aim of this subsection is to express $H^p(\Gamma, \mathcal{O}_\lambda C^{-\omega}(\Lambda))$ in terms of spectral and topological data as eigenvalues of the Laplacian, behaviour of the scattering matrix and cohomology of $\Gamma$ with values in finite-dimensional representations. We start with the case $\text{Re}(\lambda) \geq 0$.

PROPOSITION 4.12.  *If* $\text{Re}(\lambda) \geq 0$, *then* $H^p(\Gamma, \mathcal{O}_\lambda C^{-\omega}(\Lambda)) = 0$ *for all* $p \neq 1$ *and*

$$\begin{aligned}
H^1(\Gamma, \mathcal{O}_\lambda C^{-\omega}(\Lambda)) &\cong {}^\Gamma C^{-\omega}(\Lambda, V(\lambda)) \\
&\cong \begin{cases} \dim \text{ker}_{L^2}(\Delta_Y - \rho^2 + \lambda^2) & \text{Re}(\lambda) \geq 0, \lambda \neq 0 \\ \dim \text{ker}(S_0 + \text{id}) & \lambda = 0 \ . \end{cases}
\end{aligned}$$

*In particular, if* $H^1(\Gamma, \mathcal{O}_\lambda C^{-\omega}(\Lambda)) \neq 0$, *then* $\lambda \in \mathcal{F} \subset [0, \delta_\Gamma]$ *(see Corollary 3.13).*

*Proof.* By Lemma 4.11 it remains to compute $\text{coker}(\text{res})$. Since $\text{ext}_\lambda$ has at most first order poles for $\text{Re}(\lambda) \geq 0$ and $\text{res} \circ \text{ext} = \text{id}$ we can define $ev : \mathcal{O}_\lambda C^{-\omega}(B) \to {}^\Gamma C^{-\omega}(\Lambda, V(\lambda))$ by $ev(f) := \text{res}_{\mu=\lambda} \text{ext}_\mu(f_\lambda)$.

We claim that the sequence

$$(32) \qquad 0 \to {}^\Gamma\mathcal{O}_\lambda C^{-\omega}(\partial X) \overset{\text{res}}{\to} \mathcal{O}_\lambda C^{-\omega}(B) \overset{ev}{\to} {}^\Gamma C^{-\omega}(\Lambda, V(\lambda)) \to 0$$

is exact. It is clear that res is injective. Surjectivity of $ev$ follows from Lemma 3.12. Let $f \in \ker(ev)$. Then $\text{ext}(f) \in {}^{\Gamma}\mathcal{O}_{\lambda}C^{-\omega}(\partial X)$ and res $\circ$ $\text{ext}(f) = f$. It remains to show $ev \circ \text{res} = 0$. Let $\phi \in {}^{\Gamma}\mathcal{O}_{\lambda}C^{-\omega}(\partial X)$. Then $\text{ext}_{\mu} \circ \text{res}_{\mu}(\phi_{\mu}) = \phi_{\mu}$, for all $\mu \neq \lambda$, $\mu$ close to $\lambda$. It follows $ev \circ \text{res}(\phi) = \lim_{\mu \to \lambda}(\mu - \lambda)\text{ext}_{\mu} \circ \text{res}_{\mu}(\phi_{\mu}) = \lim_{\mu \to \lambda}(\mu - \lambda)\phi_{\mu} = 0$. This proves the claim.

Exactness of (32) immediately implies $\text{coker}(\text{res}) \cong {}^{\Gamma}C^{-\omega}(\Lambda, V(\lambda))$. If $\text{Re}(\lambda) > 0$, then the Poisson transform provides an isomorphism between ${}^{\Gamma}C^{-\omega}(\Lambda, V(\lambda))$ and $\ker_{L^2}(\Delta_Y - \rho^2 + \lambda^2)$. In fact, that $P_{\lambda}$ maps ${}^{\Gamma}C^{-\omega}(\Lambda, V(\lambda))$ injectively into $\ker_{L^2}(\Delta_Y - \rho^2 + \lambda^2)$ was already observed in the proof of Lemma 2.13. For the surjectivity of $P_{\lambda}$ see e.g. [8, Prop. 9.2] or [6, Lemma 2.1]. If $\text{Re}(\lambda) = 0$, $\lambda \neq 0$, then ${}^{\Gamma}C^{-\omega}(\Lambda, V(\lambda))$ as well as $\ker_{L^2}(\Delta_Y - \rho^2 + \lambda^2)$ are trivial (see Corollary 3.13 and [28], respectively). For a derivation of the latter fact in the framework of the present paper we refer to [8].

We finish the proof of the proposition by showing that ${}^{\Gamma}C^{-\omega}(\Lambda, V(0)) \cong \text{im}(S_0 - \text{id})$. From the proof of Lemma 3.11 we recall the equation $S_0 - \text{id} = C \, \text{res}_0 \circ J_0^0 \circ \text{res}_{\mu=0}\text{ext}_{\mu}$, where $C$ is some nonzero constant and $J_0^0$ is the constant term in the Laurent expansion of $\hat{J}_{\mu}$ at 0. Now by Lemma 3.12 we have $\text{im}(\text{res}_{\mu=0}\text{ext}_{\mu}) = {}^{\Gamma}C^{-\omega}(\Lambda, V(0))$, and $\text{res}_0 \circ J_0^0$ is injective by Lemma 3.7. Thus $\text{res}_0 \circ J_0^0$ provides the desired isomorphism. $\qquad\square$

In order to compute $H^*(\Gamma, \mathcal{O}_{\lambda}C^{-\omega}(\Lambda))$ for $\text{Re}(\lambda) < 0$ we need detailed information about the singularities of the intertwining operators.

LEMMA 4.13. (i) *If $\lambda \in -\rho - \mathbf{N}_0$, then the range of the Knapp-Stein intertwining operator $\hat{J}_{\lambda}$ is an irreducible finite-dimensional representation $F_{\lambda}$ of $G$.*

(ii) *If $\text{Re}(\lambda) < 0$ and $\lambda \notin -\rho - \mathbf{N}_0$, then $\hat{J}_{\lambda}$ is an isomorphism.*

(iii) *$\hat{J}_{\lambda}$ has a pole if and only if $\lambda \in \mathbf{N}_0$, and this pole is of first order.*

(iv) *We consider the following renormalized versions of the intertwining operators: If $\lambda \in \mathbf{N}$, then $J_{\mu}^{\lambda} := (\mu - \lambda)\hat{J}_{\mu}$. If $\text{Re}(\lambda) > 0$ and $\lambda \notin \mathbf{N}$, then $J_{\mu}^{\lambda} := \hat{J}_{\mu}$ (thus $J_{\mu}^{\lambda}$ is regular at $\mu = \lambda$).*

*If $\lambda \in \rho + \mathbf{N}_0$, then $\text{im}(J_{\lambda}^{\lambda}) = \ker(\hat{J}_{-\lambda})$ and $\ker(J_{\lambda}^{\lambda}) = \text{im}(\hat{J}_{-\lambda})$. If $\text{Re}(\lambda) > 0$, $\lambda \notin \rho + \mathbf{N}_0$, then $J_{\lambda}^{\lambda}$ is an isomorphism.*

(v) *If $\lambda \in -\rho - \mathbf{N}_0$, then the meromorphic family $J_{\mu}^{\lambda}$ is defined by $J_{\mu}^{\lambda} := (\mu - \lambda)^{-1}\hat{J}_{\mu}$. If $\text{Re}(\lambda) < 0$ and $\lambda \notin -\rho - \mathbf{N}_0$, then $J_{\mu}^{\lambda} := \hat{J}_{\mu}$.*

*If $\text{Re}(\lambda) < 0$, then there exists a holomorphic function $q^{\lambda}(\mu)$ which is defined in a small neighbourhood of $\lambda$ such that $q^{\lambda}(\lambda) \neq 0$ and $J_{\mu}^{\lambda} \circ J_{-\mu}^{-\lambda} = q^{\lambda}(\mu)\text{id}$.*

*Proof.* The lemma is a consequence of several well-known facts concerning intertwining operators and reducibility of principal series representations.

By [20, Ch.VI, Thm. 3.6], the principal series representation of $G$ on $C^{-\omega}(\partial X, V(\lambda))$ is (topologically) reducible if and only if $\lambda \in \pm(\rho + \mathbf{N}_0)$. If $\lambda \in \rho + \mathbf{N}_0$, then by the Cartan-Helgason Theorem [19, Ch. V, Thm. 4.1], and Casselman's Frobenius reciprocity [49, 3.8.2], $C^{-\omega}(\partial X, V(\lambda))$ contains a finite-dimensional irreducible submodule $F_{-\lambda}$. If $\mathrm{Re}(\lambda) < 0$, then by [49, 5.4.1.(2)], the range of $\hat{J}_\lambda$ is exactly the unique irreducible submodule of $C^{-\omega}(\partial X, V(-\lambda))$ (the "Langlands quotient"). If $\lambda \in -\rho - \mathbf{N}_0$, then this submodule is $F_\lambda$. This proves (i). If $\mathrm{Re}(\lambda) < 0$, $\lambda \notin -\rho - \mathbf{N}_0$, then $C^{-\omega}(\partial X, V(\pm\lambda))$ are irreducible and $\hat{J}_\lambda$ is an isomorphism. This is assertion (ii).

(iii) follows from [26, Thm. 3 and Prop 4.4] (see also [20, Ch. II, Thm. 5.4]). In order to apply Prop. 4.4 (loc. cit.) we must know that $P(0) = 0$. One can either employ irreducibility of $C^{-\omega}(\partial X, V(0))$ and [26, 7.1], or the explicit formulas for $c(\lambda)$ (and the relation (14) $P(\lambda)^{-1} = c(\lambda)c(-\lambda)$) given in [19, Ch. IV, Thm. 6.14].

We now consider (iv) and (v). If $\mathrm{Re}(\lambda) < 0$, $\lambda \notin -\rho - \mathbf{N}_0$, then by the nonvanishing result [20, Ch. II, Prop. 5.7], and the irreducibility of $C^{-\omega}(\partial X, V(\pm\lambda))$ both $J_{\pm\mu}^{\pm\lambda}$ are isomorphisms for $\mu$ in a neighbourhood of $\lambda$. Thus $J_\mu^\lambda \circ J_{-\mu}^{-\lambda} = q^\lambda(\mu)\mathrm{id}$ for some nowhere-vanishing local holomorphic function $q^\lambda$.

Let now $\lambda \in -\rho - \mathbf{N}_0$. By (13),

$$(33) \qquad \hat{J}_\mu \circ J_{-\mu}^{-\lambda} = J_{-\mu}^{-\lambda} \circ \hat{J}_\mu = r^\lambda(\mu)\mathrm{id} \ ,$$

where

$$r^\lambda(\mu) = \begin{cases} \frac{1}{P(\mu)} & n \equiv 0(2) \\ \frac{\mu-\lambda}{P(\mu)} & n \equiv 1(2) \end{cases} \ .$$

Since $J_{-\lambda}^{-\lambda}$ is regular and $\hat{J}_\lambda$ is not surjective by (i), we obtain $r^\lambda(\lambda) = 0$. We conclude that $P$ has a pole at $\lambda$ if $n \equiv 0(2)$, and that $P(\lambda) \neq 0$ if $n \equiv 1(2)$. By [26, Sec. 12], the Plancherel density $P$ has at most simple poles if $n \equiv 0(2)$, and it is holomorphic if $n \equiv 1(2)$. Thus the zero of $r^\lambda$ at $\mu = \lambda$ is simple. Now (iv) follows from (33), and

$$q^\lambda(\mu) := \begin{cases} \frac{1}{(\mu-\lambda)P(\mu)} & n \equiv 0(2) \\ \frac{1}{P(\mu)} & n \equiv 1(2) \end{cases}$$

is holomorphic and nonvanishing in a neighbourhood of $\lambda$. This finishes the proof of (v). $\qquad\square$

Fix $\lambda \in -\rho - \mathbf{N}_0$. We define the evaluation map $b : \mathcal{O}_\lambda C^{-\omega}(\partial X) \to F_\lambda$ by $b(\phi) := \hat{J}_\lambda \phi_\lambda$. Set $\mathcal{O}_\lambda^0 C^{-\omega}(\partial X) := \ker b$. Let $\mathrm{res}^0$ be the restriction of $\mathrm{res} : {}^\Gamma\mathcal{O}_\lambda C^{-\omega}(\partial X) \to \mathcal{O}_\lambda C^{-\omega}(B)$ to ${}^\Gamma\mathcal{O}_\lambda^0 C^{-\omega}(\partial X)$.

In order to make our notation more uniform we set for $\mathrm{Re}(\lambda) < 0$, $\lambda \not\in -\rho - \mathbf{N}_0$:

$$F_\lambda := 0, b := 0, \mathrm{res}^0 := \mathrm{res}, \mathcal{O}_\lambda^0 C^{-\omega}(\partial X) := \mathcal{O}_\lambda C^{-\omega}(\partial X) \ .$$

Then for all $\lambda \in \mathbf{C}$ with $\mathrm{Re}(\lambda) < 0$ we have an exact complex of $G$-modules

$$(34) \qquad 0 \to \mathcal{O}_\lambda^0 C^{-\omega}(\partial X) \to \mathcal{O}_\lambda C^{-\omega}(\partial X) \xrightarrow{b} F_\lambda \to 0 \ .$$

There is a well-defined map $J^\lambda : \mathcal{O}_\lambda^0 C^{-\omega}(\partial X) \to \mathcal{O}_{-\lambda} C^{-\omega}(\partial X)$ given by $(J^\lambda f)_{-\mu} := J_\mu^\lambda f_\mu$. It follows from Lemma 4.13 (v), that it is an isomorphism. Thus, by Proposition 4.6 the $\Gamma$-module $\mathcal{O}_\lambda^0 C^{-\omega}(\partial X)$ is acyclic.

LEMMA 4.14. *If* $\mathrm{Re}(\lambda) < 0$, *then for all* $p \geq 1$

$$(35) \qquad H^p(\Gamma, \mathcal{O}_\lambda C^{-\omega}(\partial X)) \cong H^p(\Gamma, F_\lambda) \ .$$

*Proof.* This is an immediate consequence of (34) and the $\Gamma$-acyclicity of $\mathcal{O}_\lambda^0 C^{-\omega}(\partial X)$. $\qquad\blacksquare$

Next we introduce a regularized scattering matrix. Let $\mathcal{O}_{-\lambda}^0 C^{-\omega}(B)$ $:= \ker(ev)$. Then by exactness of (32) $\mathcal{O}_{-\lambda}^0 C^{-\omega}(B) = \mathrm{im}(\mathrm{res})$. We define $S^{-\lambda} : \mathcal{O}_{-\lambda}^0 C^{-\omega}(B) \to \mathcal{O}_\lambda C^{-\omega}(B)$ by $S^{-\lambda} := \mathrm{res} \circ J^{-\lambda} \circ \mathrm{ext}$. Note that $S_\mu^{-\lambda}$ is regular at $\lambda = \mu$ since $\mathrm{ext}_\mu$ is so when restricted to $\ker(ev)$.

LEMMA 4.15. *If* $\mathrm{Re}(\lambda) < 0$, *then* $\mathrm{coker}(\mathrm{res}^0) = \mathrm{coker}(S^{-\lambda})$.

*Proof.* The maps

$$J^\lambda : {}^\Gamma\mathcal{O}_\lambda^0 C^{-\omega}(\partial X) \to {}^\Gamma\mathcal{O}_{-\lambda} C^{-\omega}(\partial X)$$

and

$$\mathrm{res} : {}^\Gamma\mathcal{O}_{-\lambda} C^{-\omega}(\partial X) \to \mathcal{O}_{-\lambda}^0 C^{-\omega}(B)$$

are isomorphisms. The inverse of $\mathrm{res} \circ J^\lambda$ is given by $\frac{1}{q^\lambda} J^{-\lambda} \circ \mathrm{ext}$ (see Lemma 4.13 (v)). Thus $J^{-\lambda} \circ \mathrm{ext}(\mathcal{O}_{-\lambda}^0 C^{-\omega}(B)) = {}^\Gamma\mathcal{O}_\lambda^0 C^{-\omega}(\partial X)$. We conclude that $\mathrm{coker}(\mathrm{res}^0) = \mathrm{coker}(S^{-\lambda})$. $\qquad\blacksquare$

PROPOSITION 4.16. *If* $\mathrm{Re}(\lambda) < 0$, *then*

(i) $H^0(\Gamma, \mathcal{O}_\lambda C^{-\omega}(\Lambda)) = 0$ ,

(ii) $\dim H^1(\Gamma, \mathcal{O}_\lambda C^{-\omega}(\Lambda)) = \dim \mathrm{coker}(\mathrm{res}) + \dim H^1(\Gamma, F_\lambda)$

$$= \dim \mathrm{coker}(S^{-\lambda}) + \dim H^1(\Gamma, F_\lambda) - \dim H^0(\Gamma, F_\lambda) \ ,$$

(iii) $H^p(\Gamma, \mathcal{O}_\lambda C^{-\omega}(\Lambda)) \cong H^p(\Gamma, F_\lambda)$ *for all* $p \geq 2$ .

*In particular*, dim $H^*(\Gamma, \mathcal{O}_\lambda C^{-\omega}(\Lambda)) < \infty$ *and*

$$\chi(\Gamma, \mathcal{O}_\lambda C^{-\omega}(\Lambda)) = \chi(\Gamma, F_\lambda) - \dim \operatorname{coker}(S^{-\lambda}) \ .$$

*Proof.* According to Lemma 4.14 the long exact cohomology sequence associated to (25) reads as follows:

$$(36) \quad 0 \to {}^\Gamma\mathcal{O}_\lambda C^{-\omega}(\partial X) \overset{\mathrm{res}}{\to} \mathcal{O}_\lambda C^{-\omega}(B) \quad \overset{\delta}{\to} \quad H^1(\Gamma, \mathcal{O}_\lambda C^{-\omega}(\Lambda))$$
$$\overset{b_1}{\to} \quad H^1(\Gamma, F_\lambda) \to 0 \ ,$$
$$(37) \qquad\qquad 0 \to H^p(\Gamma, \mathcal{O}_\lambda C^{-\omega}(\Lambda)) \quad \overset{b_p}{\to} \quad H^p(\Gamma, F_\lambda) \to 0 \ , \quad p \geq 2 \ .$$

This implies (i), (iii) and the first equation of (ii). The short exact sequence of complexes

$$
\begin{array}{ccccccccc}
0 & \to & {}^\Gamma\mathcal{O}_\lambda^0 C^{-\omega}(\partial X) & \overset{\mathrm{res}^0}{\to} & \mathcal{O}_\lambda C^{-\omega}(B) & \to & \operatorname{coker}(\mathrm{res}^0) & \to & 0 \\
 & & \downarrow \mathrm{i} & & \downarrow & & \downarrow & & \\
0 & \to & {}^\Gamma\mathcal{O}_\lambda C^{-\omega}(\partial X) & \overset{\mathrm{res}}{\to} & \mathcal{O}_\lambda C^{-\omega}(B) & \to & \operatorname{coker}(\mathrm{res}) & \to & 0
\end{array}
$$

induces the exact sequence

$$0 \to \operatorname{coker}(i) \to \operatorname{coker}(\mathrm{res}^0) \to \operatorname{coker}(\mathrm{res}) \to 0 \ .$$

Since $\mathcal{O}_\lambda^0 C^{-\omega}(\partial X)$ is acyclic we find by (34) that $\operatorname{coker}(i) \cong H^0(\Gamma, F_\lambda)$. Combining this with (36) and Lemma 4.15 we obtain the remaining assertions. $\square$

### 4.4. *The $\Gamma$-modules $\mathcal{O}_{(\lambda,k)} C^{-\omega}(\Lambda)$.*

*Definition* 4.17. For any $k \in \mathbf{N}$ we define the $\Gamma$-module $\mathcal{O}_{(\lambda,k)} C^{-\omega}(\Lambda)$ as the quotient:

$$(38) \qquad 0 \to \mathcal{O}_\lambda C^{-\omega}(\Lambda) \overset{L_\lambda^k}{\to} \mathcal{O}_\lambda C^{-\omega}(\Lambda) \to \mathcal{O}_{(\lambda,k)} C^{-\omega}(\Lambda) \to 0 \ ,$$

where $L_\lambda^k$ is defined by $(L_\lambda^k f)_\mu := (\mu - \lambda)^k f_\mu$. For any $\lambda \in \mathbf{C}$ define $k(\lambda)$ := $\operatorname{Ord}_{\mu=\lambda}\operatorname{ext}_\mu + \varepsilon(\lambda)$, where $\varepsilon(\lambda) = 0$ if $\lambda \not\in -\rho - \mathbf{N}_0$, and $\varepsilon(\lambda) = 1$ elsewhere. Here $\operatorname{Ord}_{\mu=\lambda}$ denotes the (positive) order of a pole at $\mu = \lambda$, if there is one, and zero otherwise.

There are isomorphisms of $\Gamma$-modules $\mathcal{O}_{(\lambda,1)} C^{-\omega}(\Lambda) \cong C^{-\omega}(\Lambda, V(\lambda))$. If $\operatorname{Re}(\lambda) \geq 0$, then by Corollary 3.13 we have $k(\lambda) \leq 1$.

We consider the operators

$$(L_\lambda^k)_p : H^p(\Gamma, \mathcal{O}_\lambda C^{-\omega}(\Lambda)) \to H^p(\Gamma, \mathcal{O}_\lambda C^{-\omega}(\Lambda))$$

induced by $L_\lambda^k$.

LEMMA 4.18.  *If $p \neq 1$, then $(L_\lambda)_p = 0$. If $k \geq k(\lambda)$, then $(L_\lambda^k)_1 = 0$.*

*Proof.* Because of the triviality of $H^0(\Gamma, \mathcal{O}_\lambda C^{-\omega}(\Lambda))$ it is enough to consider the case $p > 0$. If $\lambda \notin -\rho - \mathbf{N}_0$, then $H^p(\Gamma, \mathcal{O}_\lambda C^{-\omega}(\Lambda)) = 0$ for all $p \neq 1$ by Lemma 4.11; thus $(L_\lambda)_p = 0$. Let $\lambda \in -\rho - \mathbf{N}_0$. We give $F_\lambda$ the structure of a $\mathbf{C}[L_\lambda]$-module setting $L_\lambda v := 0$, $v \in F_\lambda$. Then $b$ becomes a morphism of $\Gamma$- and $\mathbf{C}[L_\lambda]$-modules. Thus (37) is an isomorphism of $\mathbf{C}[L_\lambda]$-modules. Hence $(L_\lambda)_p = 0$ for $p \geq 2$.

Consider (25) as an exact sequence of $\Gamma$- and $\mathbf{C}[L_\lambda]$-modules. Then (31) and (36) become sequences of $\mathbf{C}[L_\lambda]$-modules. Let $\lambda \in \mathbf{C}$ with $\lambda \notin -\rho - \mathbf{N}_0$. If $[f] \in H^1(\Gamma, \mathcal{O}_\lambda C^{-\omega}(\Lambda)) \cong \mathrm{coker}(\mathrm{res})$ is represented by $f \in \mathcal{O}_\lambda C^{-\omega}(B)$, then $\mathrm{ext} \circ L_\lambda^k f =: g \in {}^\Gamma \mathcal{O}_\lambda C^{-\omega}(\partial X)$ exists for $k \geq \mathrm{Ord}_{\mu=\lambda}\mathrm{ext}_\mu = k(\lambda)$. It follows that $(L_\lambda^k)_1[f] = [L_\lambda^k f] = [\mathrm{res}(g)] = 0$ for $k \geq k(\lambda)$.

Let now $\lambda \in -\rho - \mathbf{N}_0$. Consider $\phi \in H^1(\Gamma, \mathcal{O}_\lambda C^{-\omega}(\Lambda))$. We employ the sequence (36). Since $(L_\lambda)_1$ acts trivially on $H^1(\Gamma, F_\lambda)$ we have $b_1 \circ (L_\lambda)_1(\phi) = (L_\lambda)_1 \circ b_1(\phi) = 0$. Thus $(L_\lambda)_1\phi = \delta(f)$ for some $f \in \mathcal{O}_\lambda C^{-\omega}(B)$. Suppose that $k \geq \mathrm{Ord}_{\mu=\lambda}\mathrm{ext}_\mu + 1 = k(\lambda)$. Putting $g := \mathrm{ext} \circ (L_\lambda^{k-1})_1(f)$ we obtain $(L_\lambda^k)_1\phi = \delta \circ \mathrm{res}(g) = 0$. This finishes the proof of the lemma. $\square$

We recall the definition of the first derived Euler characteristic $\chi_1(\Gamma, V)$ of a $\Gamma$-module $V$:

$$\chi_1(\Gamma, V) := \sum_{p=0}^{\infty} (-1)^p p \, \dim H^p(\Gamma, V) \ .$$

The following proposition contains the first three assertions of Theorem 1.3 and the fact that equation (3) implies equation (4).

PROPOSITION 4.19. *Let $\lambda \in \mathbf{C}$. Then*

(i) $\dim H^*(\Gamma, \mathcal{O}_{(\lambda,k)} C^{-\omega}(\Lambda)) < \infty$. *In particular,* $\dim H^*(\Gamma, C^{-\omega}(\Lambda, V(\lambda))) < \infty$.

(ii) $\chi(\Gamma, \mathcal{O}_{(\lambda,k)} C^{-\omega}(\Lambda)) = 0$.

 *If $k \geq k(\lambda)$, then*

(iii) $\dim H^*(\Gamma, \mathcal{O}_{(\lambda,k+1)} C^{-\omega}(\Lambda)) = \dim H^*(\Gamma, \mathcal{O}_{(\lambda,k)} C^{-\omega}(\Lambda))$.

(iv) $\chi_1(\Gamma, \mathcal{O}_{(\lambda,k)} C^{-\omega}(\Lambda)) = \chi(\Gamma, \mathcal{O}_\lambda C^{-\omega}(\Lambda))$.

*Proof.* Assertions (i) and (ii) follow from Proposition 4.16 and the long exact cohomology sequence associated to (38). If $k \geq k(\lambda)$, then by Lemma 4.18 this long exact sequence splits into short exact sequences ($p = 0, 1, \ldots$),

$$0 \to H^p(\Gamma, \mathcal{O}_\lambda C^{-\omega}(\Lambda)) \to H^p(\Gamma, \mathcal{O}_{(\lambda,k)} C^{-\omega}(\Lambda)) \to H^{p+1}(\Gamma, \mathcal{O}_\lambda C^{-\omega}(\Lambda)) \to 0 \ .$$

Now assertions (iii) and (iv) follow, too. $\square$

PROPOSITION 4.20.   *If $\lambda \in \mathbf{C}$ and $k \geq k(\lambda)$, then*

$$\left.\begin{array}{l} \chi(\Gamma, \mathcal{O}_\lambda C^{-\omega}(\Lambda)) \\ \chi_1(\Gamma, \mathcal{O}_{(\lambda,k)} C^{-\omega}(\Lambda)) \end{array}\right\} = \left\{ \begin{array}{ll} -\dim \ker_{L^2}(\Delta_Y - \rho^2 + \lambda^2) & \operatorname{Re}(\lambda) \geq 0, \quad \lambda \neq 0 \\ -\dim \ker(S_0 + \operatorname{id}) & \lambda = 0 \\ -\dim \operatorname{coker}(S^{-\lambda}) + \chi(\Gamma, F_\lambda) & \operatorname{Re}(\lambda) < 0 \ . \end{array} \right.$$

*Proof.* The proposition is an immediate consequence of Propositions 4.12, 4.16 and 4.19. □

In fact, a closer study of the long exact cohomology sequence associated to (38) in combination with Propositions 4.12 and 4.16 gives more precise information about the cohomology groups $H^*(\Gamma, \mathcal{O}_{(\lambda,k)} C^{-\omega}(\Lambda))$.

COROLLARY 4.21.   *If $\operatorname{Re}(\lambda) \geq 0$, then set $F_\lambda := 0$. For all $\lambda \in \mathbf{C}$ and $k \in \mathbf{N}$ we have the following isomorphisms and exact sequences:*

$$H^0(\Gamma, \mathcal{O}_{(\lambda,k)} C^{-\omega}(\Lambda)) \quad \cong \quad \ker(L_\lambda^k)_1 \ ,$$

$$0 \to \operatorname{coker}(L_\lambda^k)_1 \to H^1(\Gamma, \mathcal{O}_{(\lambda,k)} C^{-\omega}(\Lambda)) \quad \to \quad H^2(\Gamma, F_\lambda) \to 0 \ ,$$

$$0 \to H^p(\Gamma, F_\lambda) \to H^p(\Gamma, \mathcal{O}_{(\lambda,k)} C^{-\omega}(\Lambda)) \quad \to \quad H^{p+1}(\Gamma, F_\lambda) \to 0 \ , \quad p \geq 2 \ .$$

*In particular,*

$$\dim H^1(\Gamma, \mathcal{O}_{(\lambda,k)} C^{-\omega}(\Lambda)) = \dim H^0(\Gamma, \mathcal{O}_{(\lambda,k)} C^{-\omega}(\Lambda)) + \dim H^2(\Gamma, F_\lambda)$$

*and if $k \geq k(\lambda)$, then*

$$\dim H^0(\Gamma, \mathcal{O}_{(\lambda,k)} C^{-\omega}(\Lambda)) = \dim \operatorname{coker}(\operatorname{res}) + \dim H^1(\Gamma, F_\lambda)$$

$$= \left\{ \begin{array}{ll} \dim {}^\Gamma C^{-\omega}(\Lambda, V(\lambda)) = \dim \ker_{L^2}(\Delta_Y - \rho^2 + \lambda^2) & \operatorname{Re}(\lambda) \geq 0, \quad \lambda \neq 0 \\ \dim {}^\Gamma C^{-\omega}(\Lambda, V(\lambda)) = \dim \ker(S_0 + \operatorname{id}) & \lambda = 0 \\ \dim \operatorname{coker}(S^{-\lambda}) + \dim H^1(\Gamma, F_\lambda) - \dim H^0(\Gamma, F_\lambda) & \operatorname{Re}(\lambda) < 0 \ . \end{array} \right.$$

We conclude this subsection with a generalization of Lemma 3.12 to the case $\lambda \notin -\rho - \mathbf{N}_0$. For $k \geq k(\lambda)$ the singular part of ext at $\lambda$ defines a map $\operatorname{ext}_\lambda^{<0} : \mathcal{O}_\lambda C^{-\omega}(B) \to {}^\Gamma \mathcal{O}_{(\lambda,k)} C^{-\omega}(\Lambda)$ by

$$\operatorname{ext}_\lambda^{<0}(f) := \operatorname{ext} \circ L_\lambda^k(f) \mod L_\lambda^k(\mathcal{O}_\lambda C^{-\omega}(\partial X)) \ .$$

PROPOSITION 4.22.   *If $\lambda \notin -\rho - \mathbf{N}_0$ and $k \geq k(\lambda)$, then*

$$\operatorname{ext}_\lambda^{<0} : \mathcal{O}_\lambda C^{-\omega}(B) \to {}^\Gamma \mathcal{O}_{(\lambda,k)} C^{-\omega}(\Lambda)$$

*is surjective. Moreover,*

$${}^\Gamma C^{-\omega}(\Lambda, V(\lambda)) = \left\{ \operatorname{res}_{\mu=\lambda}(\operatorname{ext}_\mu(f_\mu)) \mid \begin{array}{l} f \in \mathcal{O}_\lambda C^{-\omega}(B) \text{ such that } \operatorname{ext}_\mu(f_\mu) \\ \text{has a pole of first order at } \mu = \lambda \end{array} \right\} \ .$$

*In particular, $\operatorname{ext}_\mu$ is regular at $\lambda$ if and only if ${}^\Gamma C^{-\omega}(\Lambda, V(\lambda)) = 0$.*

*Proof.* $\text{ext}_\lambda^{<0}$ factorizes over $\text{coker}(\text{res})$, and since $\text{res} \circ \text{ext} = \text{id}$, this factorization is injective. By Corollary 4.21, $\dim \text{coker}(\text{res}) = \dim \, {}^\Gamma\mathcal{O}_{(\lambda,k)}C^{-\omega}(\Lambda)$. This implies surjectivity of $\text{ext}_\lambda^{<0}$. Again by Corollary 4.21 it follows that

$$\dim \, {}^\Gamma C^{-\omega}(\Lambda, V(\lambda)) = \dim \ker(L_\lambda : \text{coker}(\text{res}) \to \text{coker}(\text{res})) \ .$$

The remaining assertions of the proposition are now obvious. $\qquad\square$

*Remark.* Since in general for $\lambda \in -\rho - \mathbf{N}_0$ we have $\dim H^1(\Gamma, F_\lambda) \neq 0$, the map $\text{ext}_\lambda^{<0}$ is not surjective in view of the formula

$$\dim \, {}^\Gamma\mathcal{O}_{(\lambda,k)}C^{-\omega}(\Lambda) = \dim \text{ coker}(\text{res}) + \dim H^1(\Gamma, F_\lambda) \ .$$

## 5. The singularities of the Selberg zeta function

### 5.1. *The embedding trick.*

Let $Z_S(\lambda)$ denote the Selberg zeta function associated to $\Gamma$ introduced in subsection 1.1. For $\dim(X)$ even the spectral description of its singularities was worked out in Patterson-Perry [38]. This description simplifies considerably if we assume that $\delta_\Gamma < 0$.

We are going to prove the remaining assertion (iv), equation (3), of Theorem 1.3 in two steps:

(i) First we employ the embedding trick (which was already used in the proof of Proposition 2.22) in order to show in Corollary 5.5 that the equality (3) under the additional assumptions $\delta_\Gamma < 0$ and $\dim(X) \equiv 0(2)$ implies (3) in general.

(ii) Then we prove (3) under the additional assumptions $\delta_\Gamma < 0$ and $\dim(X) \equiv 0(2)$.
   In the present subsection we are concerned with step (i).

We adopt the notation $G^n$, $\partial X^n$, etc. as introduced in the proof of Proposition 2.22.

**PROPOSITION 5.1.**

$$(39) \qquad \chi(\Gamma, \mathcal{O}_\lambda C^{-\omega}(\Lambda^n)) = \chi(\Gamma, \mathcal{O}_{\lambda-\frac{1}{2}}C^{-\omega}(\Lambda^{n+1})) - \chi(\Gamma, \mathcal{O}_{\lambda+\frac{1}{2}}C^{-\omega}(\Lambda^{n+1})) \ .$$

*Proof.* We will construct an exact sequence of $\Gamma$-modules

$$(40) \qquad 0 \to \mathcal{O}_\lambda C^{-\omega}(\Lambda^n) \xrightarrow{i_*} \mathcal{O}_{\lambda-\frac{1}{2}}C^{-\omega}(\Lambda^{n+1}) \xrightarrow{\Phi} \mathcal{O}_{\lambda+\frac{1}{2}}C^{-\omega}(\Lambda^{n+1}) \to 0 \ .$$

Here $i_*$ is the push forward associated to the embedding $\partial X^n \hookrightarrow X^{n+1}$ (see the proof of 2.22). The map $\Phi$ is multiplication by a $G^n$-invariant real analytic section $\phi$ of $V(\rho^{n+1}+1)^{n+1} \to \partial X^{n+1}$ which vanishes on $\partial X^n$ of first order. We will construct $\phi$ in Lemma 5.2 and show exactness of (40) in Lemma 5.3 below. Then Proposition 5.1 will follow in view of the long exact sequence in cohomology associated to (40) and the fact that the cohomology groups are finite-dimensional.

LEMMA 5.2.    *There exists a real analytic function $\phi \in C^\omega(G^{n+1})$ which satisfies*:

(i)  $\{\phi = 0\} = G^n M^{n+1} A^{n+1} N^{n+1}$,

(ii)  $\phi(hg) = \phi(g)$, *for all $h \in G^n, g \in G^{n+1}$*,

(iii)  $\phi(gman) = a\phi(g)$, *for all $man \in M^{n+1}A^{n+1}N^{n+1}$*.

(iv)  *If $\phi$ is considered as a section of $V(\rho^{n+1}+1)^{n+1}$, then it vanishes on $\partial X^n$ of first order.*

*Proof.*    We consider the standard representation $\pi$ of $G^{n+1} = SO(1, n+1)_0$ on $\mathbf{R}^{n+2}$. Let $\{e_i \,|\, i = 0, \ldots, n+1\}$ be the standard base and $\langle .,. \rangle$ be the $G^{n+1}$-invariant bilinear form $\mathrm{diag}(-1, 1, \ldots, 1)$. We specify the Iwasawa decomposition by $K^{n+1} := \{g \in G^{n+1} \,|\, ge_0 = e_0\} \cong SO(n+1)$ and

$$A^+ := \left\{ \begin{pmatrix} \cosh(t) & \sinh(t) & 0 \\ \sinh(t) & \cosh(t) & 0 \\ 0 & 0 & \mathrm{id}_{\mathbf{R}^n} \end{pmatrix} \,\Big|\, t > 0 \right\}.$$

Then we define $\phi(g) := \langle e_{n+1}, \pi(g)(e_0 + e_1)\rangle$. Since $e_{n+1}$ is $G^n$-invariant, $e_0 + e_1$ is $M^{n+1}$-invariant and the highest weight vector with respect to $A^+$, the assertions (ii) and (iii) follow. Thus we can consider $\phi$ as a $G^n$-invariant section of $V(\rho^{n+1}+1)^{n+1}$. In the usual trivialization by the $K^{n+1}$-invariant section it is simply given by the height function $\phi(x) = x_{n+1}$, $x = (x_1, \ldots, x_{n+1}) \in S^n = \partial X^{n+1}$. This is easily seen using the $K^{n+1}$-equivariant embedding $S^n \ni x \mapsto (1, x) \in \mathbf{R}^{n+2}$. We conclude that $\{x \in \partial X^{n+1} \,|\, \phi(x) = 0\} = \{x \in \partial X^{n+1} \,|\, x_{n+1} = 0\} = \partial X^n$, and that $\phi$ vanishes on $\partial X^n$ of first order.    □

We now define $\Phi_\mu : C^\omega(\partial X^{n+1}, V(\mu)^{n+1}) \to C^\omega(\partial X^{n+1}, V(\mu+1)^{n+1})$ by $\Phi_\mu(f)(g) := \phi(g)f(g)$. Indeed by Lemma 5.2 (iii),

$$\Phi_\mu(f)(gman) = \phi(gman)f(gman) = a\phi(g)a^{\mu-\rho^{n+1}}f(g) = a^{\mu+1-\rho^{n+1}}\Phi_\mu(f)(g)\,,$$

for all $g \in G^{n+1}$, $man \in M^{n+1}A^{n+1}N^{n+1}$, and thus $\Phi(f)$ is a section of $V(\mu+1)^{n+1}$. By Lemma 5.2 (ii), the map $\Phi_\mu$ is $G^n$-equivariant. $\Phi_\mu$ extends to hyperfunction sections and satisfies $\mathrm{supp}(\Phi_\mu f) \subset \mathrm{supp}(f)$, for all $f \in C^{-\omega}(\partial X^{n+1}, V(\mu)^{n+1})$.

We define $\Phi : \mathcal{O}_{\lambda-\frac{1}{2}} C^{-\omega}(\Lambda^{n+1}) \to \mathcal{O}_{\lambda+\frac{1}{2}} C^{-\omega}(\Lambda^{n+1})$ by $(\Phi f)_\mu := \Phi_\mu f_\mu$.

LEMMA 5.3.   *The sequence* (40) *is exact.*

*Proof.* Using Lemmas 4.2 and 4.3 we obtain the following diagram:

$$
\begin{array}{ccccccccc}
 & & 0 & & 0 & & 0 & & \\
 & & \downarrow & & \downarrow & & \downarrow & & \\
0 & \to & \mathcal{O}_\lambda C^{-\omega}(\partial X^n) & \stackrel{i_*}{\to} & \mathcal{O}_{\lambda-\frac{1}{2}} C^{-\omega}(\partial X^{n+1}) & \stackrel{\Phi}{\to} & \mathcal{O}_{\lambda+\frac{1}{2}} C^{-\omega}(\partial X^{n+1}) & \to & 0 \\
 & & \downarrow & & \downarrow & & \downarrow & & \\
0 & \to & \mathcal{O}_\lambda C^{-\omega}(\Omega^n) & \stackrel{i_*}{\to} & \mathcal{O}_{\lambda-\frac{1}{2}} C^{-\omega}(\Omega^{n+1}) & \stackrel{\Phi}{\to} & \mathcal{O}_{\lambda+\frac{1}{2}} C^{-\omega}(\Omega^{n+1}) & \to & 0 \\
 & & \downarrow & & \downarrow & & \downarrow & & \\
 & & 0 & & 0 & & 0 & & \\
\end{array}
,
$$

where the rows are exact, and the columns are surjective. The sequence (40) is just the complex of kernels of the columns of the diagram above, and this complex is exact by the Snake lemma. This finishes the proof of both the lemma and Proposition 5.1. □

The significance of Proposition 5.1 is that the numbers $\chi(\Gamma, \mathcal{O}_\lambda C^{-\omega}(\Lambda^n))$ behave in the same way with respect to embeddings $G^n \hookrightarrow G^{n+1}$ as the orders of the singularities of the Selberg zeta function $Z_S$. Let $\Gamma \subset G^n$ be convex cocompact. If we consider $\Gamma$ as a convex cocompact subgroup of $G^{n+1}$, then $Z_{S,n+1}(s)$ denotes the corresponding Selberg zeta function. The following elementary fact was noted in [34].

LEMMA 5.4.
$$
Z_{S,n+1}(s) = \prod_{j=0}^{\infty} Z_{S,n}(s + j + \frac{1}{2})
$$

*and, consequently,*

$$
(41) \quad \operatorname{ord}_{s=\lambda} Z_{S,n+1}(s) = \sum_{j=0}^{\infty} \operatorname{ord}_{s=\lambda+j+\frac{1}{2}} Z_{S,n}(s) \ ,
$$
$$
\operatorname{ord}_{s=\lambda} Z_{S,n}(s) = \operatorname{ord}_{s=\lambda-\frac{1}{2}} Z_{S,n+1}(s) - \operatorname{ord}_{s=\lambda+\frac{1}{2}} Z_{S,n+1}(s) \ .
$$

Combining equation (41) with Proposition 5.1 we obtain:

COROLLARY 5.5.   *Let* $\Gamma \subset G^n$ *be convex cocompact. If equation* (3) *holds true for* $\Gamma$ *viewed as a subgroup of* $G^{n+1}$, *then so does it for* $\Gamma$ *viewed as a subgroup of* $G^n$. *In particular, it is sufficient to prove* (3) *under the assumptions* $\delta_\Gamma < 0$ *and* $\dim(X) \equiv 0(2)$.

### 5.2. *Singularities and cohomology.*

In this subsection we prove Theorem 1.3 (iv), (3), under the assumptions that $n \equiv 0(2)$ and $\delta_\Gamma < 0$. In [38] the order of the singularities of the Selberg zeta function was expressed in terms of traces of residues of certain meromorphic families of operators. Our main task is to link this description with our computation of $\chi(\Gamma, \mathcal{O}_\lambda C^{-\omega}(\Lambda))$.

First observe that (3) holds true for $\text{Re}(\lambda) \geq 0$. Indeed because of our assumption $\delta_\Gamma < 0$ on the one hand the infinite product (2) defining the Selberg zeta function $Z_S(\lambda)$ converges and thus $\text{ord}_{\mu=\lambda} Z_S(\mu) = 0$. On the other hand $\chi(\Gamma, \mathcal{O}_\lambda C^{-\omega}(\Lambda)) = 0$ by Proposition 4.12.

It remains to prove (3) for $\text{Re}(\lambda) < 0$. First we give an expression for $\dim \text{coker}(S^{-\lambda} : \mathcal{O}_{-\lambda}^0 C^{-\omega}(B) \to \mathcal{O}_\lambda C^{-\omega}(B))$ in terms of the trace of the residue of the logarithmic derivative of $S^{-\lambda}$. Because of the assumption $\delta_\Gamma < 0$ we have $\mathcal{O}_{-\lambda}^0 C^{-\omega}(B) = \mathcal{O}_{-\lambda} C^{-\omega}(B)$.

In the following we review results of Patterson-Perry [38]. We fix an analytic Riemannian metric on $B$ in the canonical conformal class. This metric induces a volume form which we employ in order to identify all bundles $V_B(\lambda)$ with $B \times \mathbf{C} = V_B(\rho)$. Then the scattering matrix becomes a germ at $\lambda$ of a meromorphic family of operators $S_\mu^{-\lambda}$ on $C^*(B)$, $* = \pm\omega, \pm\infty$.

Let $\Delta_B$ be the Laplace operator on $B$ associated to the Riemannian metric. Viewed as an unbounded symmetric operator on the Hilbert space $\mathcal{H} := L^2(B)$ the sum $\Delta_B + 1$ is positive. Let $P := \sqrt{\Delta_B + 1}$ be the positive square root defined by spectral theory. Then $P := \sqrt{\Delta_B + 1}$ is a pseudodifferential operator of order 1. In particular, $P$ and its complex powers $P^\mu$ act on $C^{\pm\infty}(B)$. For $\mu$ close to $\lambda$ the scattering matrix can be factorized as

$$(42) \qquad S_{-\mu}^{-\lambda} = P^{-\mu}(\text{id} + K(-\mu))P^{-\mu} ,$$

where $K(-\mu)$ is a holomorphic family of pseudodifferential operators belonging to the $(n+1)^{\text{th}}$ Schatten class (i.e. $K(-\mu)^{n+1}$ is of trace class). The inverse $(1 + K(-\mu))^{-1}$ is a meromorphic family of operators with finite-dimensional residues.

PROPOSITION 5.6.

$$\dim \text{coker}\left(S^{-\lambda} : \mathcal{O}_{-\lambda} C^{-\omega}(B) \to \mathcal{O}_\lambda C^{-\omega}(B)\right) = \text{Tr} \, \text{res}_{\mu=-\lambda}(1 + K(\mu))^{-1} K'(\mu) ,$$

where $K'(\mu)$ denotes the derivative of $K(\mu)$ with respect to $\mu$.

*Proof.* Let $\mathcal{O}_\lambda \mathcal{H}$ denote the space of germs at $\lambda$ of holomorphic families of vectors in $\mathcal{H}$ and let $(1 + K) : \mathcal{O}_{-\lambda} \mathcal{H} \to \mathcal{O}_\lambda \mathcal{H}$ be given by $((1 + K)f)_\mu := (1 + K(-\mu))f_{-\mu}$, $f \in \mathcal{O}_{-\lambda} \mathcal{H}$.

LEMMA 5.7.

$$\dim \text{coker}(S^{-\lambda}) = \dim \text{coker}(1 + K).$$

*Proof.* For $k \in \mathbf{N}$ we define

$$\mathcal{O}_{(\lambda,k)}C^{-\omega}(B) := \operatorname{coker}(L_\lambda^k : \mathcal{O}_\lambda C^{-\omega}(B) \to \mathcal{O}_\lambda C^{-\omega}(B)) \ .$$

Let $S^{(-\lambda,k)} : \mathcal{O}_{(-\lambda,k)}C^{-\omega}(B) \to \mathcal{O}_{(\lambda,k)}C^{-\omega}(B)$ be the operator induced by $S^{-\lambda}$. Fix $k \geq k(\lambda)$. If $f \in \operatorname{im}(L_\lambda^k)$, then $g := \frac{1}{q^\lambda}S^\lambda(f) \in \mathcal{O}_{-\lambda}C^{-\omega}(B)$ is defined and $f = S^{-\lambda}(g)$ (see Lemma 4.13 (v)). It follows that $L_\lambda^k$ acts trivially on $\operatorname{coker}(S^{-\lambda})$. Hence we can identify $\operatorname{coker}(S^{-\lambda})$ with $\operatorname{coker}(S^{(-\lambda,k)})$ in the natural way. As spaces of finite Taylor series with values in a Fréchet space, the spaces $\mathcal{O}_{(\pm\lambda,k)}C^{-\omega}(B)$ are Fréchet spaces, too. By Lemma 2.20 the map $S^{(-\lambda,k)}$ is continuous with finite-dimensional cokernel. Hence it has closed range by the open mapping theorem. Thus the induced topology on $\operatorname{coker}(S^{-\lambda}) = \operatorname{coker}(S^{(-\lambda,k)})$ is Hausdorff.

Let $\mathcal{O}_\lambda C^{-\infty}(B)$ denote the space of germs at $\lambda$ of holomorphic families of distributions on $B$, and let $\mathcal{O}_{(\lambda,k)}C^{-\infty}(B)$ and $\mathcal{O}_{(\lambda,k)}\mathcal{H}$ be the corresponding quotient spaces of finite Taylor series. We define

$$P : \mathcal{O}_\lambda C^{-\infty}(B) \to \mathcal{O}_\lambda C^{-\infty}(B)$$

by $(Pf)_\mu = P^{-\mu}f_\mu$. We denote the induced operator on $\mathcal{O}_{(\lambda,k)}C^{-\infty}(B)$ by the same symbol. The composition

$$p : \mathcal{O}_{(\lambda,k)}\mathcal{H} \hookrightarrow \mathcal{O}_{(\lambda,k)}C^{-\infty}(B) \xrightarrow{P} \mathcal{O}_{(\lambda,k)}C^{-\infty}(B) \hookrightarrow \mathcal{O}_{(\lambda,k)}C^{-\omega}(B)$$

has dense range.

Now let $h \in \mathcal{O}_\lambda \mathcal{H}$. We claim that $h \in \operatorname{im}(1 + K)$ if and only if $Ph$ considered as a hyperfunction is in the image of $S^{-\lambda}$. Indeed, let $Ph = S^{-\lambda}j$ for some $j \in \mathcal{O}_\lambda C^{-\omega}(B)$. Then $-\mu \mapsto g_{-\mu} := (1 + K(\mu))^{-1}(h_\mu)$ defines a meromorphic family of vectors in $\mathcal{H}$. It is regular at $\mu = \lambda$ since $P^\mu(g_{-\mu}) = j_{-\mu}$ by (42) and $P^\mu$ is injective. Vice versa, let $h = (1 + K)g$ for some $g \in \mathcal{O}_{-\lambda}\mathcal{H}$. The family $j_{-\mu} := P^\mu(g_{-\mu})$, $\mu$ near $\lambda$, defines an element $j \in \mathcal{O}_{-\lambda}C^{-\infty}(B)$ satisfying $Ph = S^{-\lambda}j$.

The above claim implies that $p$ induces a map $p_* : \operatorname{coker}(1 + K) \to \operatorname{coker}(S^{-\lambda})$ which is injective. Moreover, it has dense range. Since $\operatorname{coker}(S^{-\lambda})$ is finite-dimensional and Hausdorff $p_*$ must be surjective. The lemma follows. $\qquad\square$

The following lemma is known in one form or another (compare [15]). For completeness we include a proof here.

LEMMA 5.8.

$$\operatorname{Tr} \operatorname{res}_{\mu=-\lambda}(1 + K(\mu))^{-1}K'(\mu) = \dim \operatorname{coker}(1 + K) \ .$$

*Proof.* Set $s := -\lambda$. Let $P(\mu)$ be the holomorphic family of finite-dimensional projections given by

$$P(\mu) := \frac{1}{2\pi i} \oint_c \frac{1}{z - 1 - K(\mu)} dz \ ,$$

where the path of integration is a small circle enclosing $0 \in \mathbf{C}$ counterclockwise and $\mu$ is close to $s$. There is a holomorphic family of invertible operators $U(\mu)$ ([43, Thm. XII.12]) such that $U(\mu)^{-1} P(\mu) U(\mu) = P(s)$. We define $T(\mu) := U(\mu)^{-1}(1 + K(\mu))U(\mu)$. Then

$$
\begin{aligned}
T(\mu)P(s) &= U(\mu)^{-1}(1 + K(\mu))P(\mu)U(\mu) \\
&= U(\mu)^{-1}P(\mu)(1 + K(\mu))U(\mu) = P(s)T(\mu) \ .
\end{aligned}
$$

Let $V := P(s)\mathcal{H}$ and $W := (1 - P(s))\mathcal{H}$. Then

$$T(\mu) = \begin{pmatrix} A(\mu) & 0 \\ 0 & B(\mu) \end{pmatrix} : \begin{matrix} V \\ \oplus \\ W \end{matrix} \rightarrow \begin{matrix} V \\ \oplus \\ W \end{matrix} \ ,$$

where $B(\mu)$ is invertible for $\mu - s$ small.

We claim that

$$\mathrm{Tr}\,\mathrm{res}_{\mu=s} T(\mu)^{-1} T'(\mu) = \mathrm{Tr}\,\mathrm{res}_{\mu=s} (1 + K(\mu))^{-1} K'(\mu) \ .$$

In fact

$$
\begin{aligned}
T'(\mu) &= U(\mu)^{-1}(1 + K(\mu))U'(\mu) + U(\mu)^{-1} K'(\mu)U(\mu) \\
&\quad - U(\mu)^{-1} U'(\mu) U(\mu)^{-1}(1 + K(\mu))U(\mu) \ .
\end{aligned}
$$

Using the facts that all singular terms of $(1 + K(\mu))^{-1}$ are finite-dimensional and that the trace is cyclic, we compute

$$
\begin{aligned}
&\mathrm{Tr}\,\mathrm{res}_{\mu=s} T(\mu)^{-1} T'(\mu) \\
=\ &\mathrm{Tr}\,\mathrm{res}_{\mu=s} U(\mu)^{-1}(1 + K(\mu))^{-1} U(\mu) \\
&[U(\mu)^{-1}(1 + K(\mu))U'(\mu) + U(\mu)^{-1} K'(\mu)U(\mu) \\
&- U(\mu)^{-1} U'(\mu) U(\mu)^{-1}(1 + K(\mu))U(\mu)] \\
=\ &\mathrm{Tr}\,\mathrm{res}_{\mu=s} \\
&[U(\mu)^{-1} U'(\mu) + U(\mu)^{-1}(1 + K(\mu))^{-1} K'(\mu)U(\mu) \\
&- U(\mu)^{-1}(1 + K(\mu))^{-1} U'(\mu) U(\mu)^{-1}(1 + K(\mu))U(\mu)] \\
=\ &\mathrm{Tr}\,\mathrm{res}_{\mu=s} (1 + K(\mu))^{-1} K'(\mu) \ .
\end{aligned}
$$

This proves the claim.

Let $T : \mathcal{O}_s \mathcal{H} \to \mathcal{O}_s \mathcal{H}$ be given by $(Tf)(\mu) = T(\mu)f(\mu)$, $f \in \mathcal{O}_s \mathcal{H}$. Then we have

$$\dim \mathrm{coker}(T) = \dim \mathrm{coker}(1 + K) \ .$$

Now

$$\operatorname{Tr} \operatorname{res}_{\mu=s} T(\mu)^{-1} T'(\mu) = \operatorname{Tr} \operatorname{res}_{\mu=s} A(\mu)^{-1} A'(\mu)$$

for the holomorphic family of operators $A(\mu)$ on the finite-dimensional space $V$. Moreover $\dim \operatorname{coker}(T) = \dim \operatorname{coker}(A)$, where $A : \mathcal{O}_s V \to \mathcal{O}_s V$ is given by $(Af)(\mu) := A(\mu)f(\mu)$, $f \in \mathcal{O}_s V$. In order to finish the proof of the proposition we must show that

$$\operatorname{Tr} \operatorname{res}_{\mu=s} A(\mu)^{-1} A'(\mu) = \dim \operatorname{coker}(A) \ .$$

Now

$$
\begin{aligned}
\operatorname{Tr} \operatorname{res}_{\mu=s} A(\mu)^{-1} A'(\mu) &= \operatorname{res}_{\mu=s} \operatorname{Tr} A(\mu)^{-1} A'(\mu) \\
&= \operatorname{res}_{\mu=s} \frac{\det(A(\mu))'}{\det(A(\mu))} \\
&= \operatorname{ord}_{\mu=s} \det(A(\mu)) \ .
\end{aligned}
$$

By Gauss's algorithm $A(\mu)$ can be transformed to a holomorphic family of diagonal matrices $\tilde{A}(\mu)$ through multiplication from the left and right with holomorphic matrix functions with invertible determinants. We have $\dim \operatorname{coker}(\tilde{A}) = \dim \operatorname{coker}(A)$, where $\tilde{A} : \mathcal{O}_s V \to \mathcal{O}_s V$ is given by $(\tilde{A}f)(\mu) := \tilde{A}(\mu)f(\mu)$, $f \in \mathcal{O}_s V$, and $\operatorname{ord}_{\mu=s} \det(\tilde{A}(\mu)) = \operatorname{ord}_{\mu=s} \det(A(\mu))$. But for holomorphic families of diagonal matrices the equation $\dim \operatorname{coker} \tilde{A} = \operatorname{ord}_{\mu=s} \det(\tilde{A}(\mu))$ is obvious. This finishes the proof of the lemma. $\qquad \square$

Proposition 5.6 follows from Lemmas 5.7 and 5.8. $\qquad \square$

We now recall the description of the singularities of the Selberg zeta function $Z_S(\lambda)$ for $\operatorname{Re}(\lambda) < 0$ given in [38]. Note that our standing hypothesis is $n \equiv 0(2)$ and $\delta_\Gamma < 0$. This simplifies things considerably because the point spectrum of $\Delta_Y$ is absent. Let $\operatorname{Re}(\lambda) < 0$ and set $n_\lambda := \operatorname{Tr} \operatorname{res}_{\mu=-\lambda}(1 + K(\mu))^{-1} K'(\mu)$. Then

$$\operatorname{ord}_{s=\lambda} Z_S(s) = \begin{cases} n_\lambda & \lambda \notin -\rho - \mathbf{N}_0 \\ n_\lambda - \chi(Y) \dim F_\lambda & \lambda \in -\rho - \mathbf{N}_0 \end{cases} .$$

Since $Y$ has the homotopy type of a finite CW-complex we have $\chi(Y, F_\lambda) = \chi(Y) \dim F_\lambda$. Equation (3) now follows from Propositions 5.6 and 4.20.

By Corollary 5.5 a proof of (3) under the assumptions $n \equiv 0(2)$ and $\delta_\Gamma < 0$ implies (3) without these assumptions. This finishes the proof of Theorem 1.3.

Universität Göttingen, Göttingen, Germany
*E-mail address*: bunke@uni-math.gwdg.de
*E-mail address*: olbrich@uni-math.gwdg.de